\documentclass[11pt,a4paper]{article}
\usepackage{multirow}
\usepackage{graphicx}
\usepackage{subfigure}
\usepackage[T1]{fontenc}
\usepackage{geometry}
\usepackage{amsbsy,amsmath,latexsym,amsfonts, epsfig, color, authblk, amssymb, graphics, bm}
\usepackage{epsf,slidesec,epic,eepic}
\usepackage{fancybox}
\usepackage{fancyhdr}
\usepackage{setspace}
\usepackage{nccmath}
\usepackage[titletoc]{appendix}
\usepackage[colorlinks = true, linkcolor = red, citecolor = blue, pagebackref = true]{hyperref}

\newtheorem{theorem}{Theorem}
\newtheorem{lemma}[theorem]{Lemma}

\newtheorem{corollary}[theorem]{Corollary}
\newtheorem{assumption}{Assumption}

\newenvironment{proof}{{\noindent \bf Proof:}}{\hfill$\Box$\medskip}

\definecolor{lred}{rgb}{1,0.8,0.8}
\definecolor{lblue}{rgb}{0.8,0.8,1}
\definecolor{dred}{rgb}{0.6,0,0}
\definecolor{dblue}{rgb}{0,0,0.5}
\definecolor{dgreen}{rgb}{0,0.5,0.5}

\title{A Rank-Corrected Procedure for Matrix Completion with Fixed Basis Coefficients}
\author{Weimin Miao\footnote{Risk Management Institute, National University of Singapore, 21 Heng Mui Keng Terrace, Singapore 119613 (rmimw@nus.edu.sg). This author's research is supported in part by Willis Research Network.},  \ \ Shaohua Pan\footnote{Department of Mathematics, South China University of Technology, Tianhe District of Guangzhou City, China (shhpan@scut.edu.cn).}  \ \ and \  Defeng Sun\footnote{Department of Mathematics and Risk Management Institute, National University of Singapore, 10 Lower Kent Ridge Road, Singapore 119076 (matsundf@nus.edu.sg). This author's research is supported in part by Academic Research Fund under grant R-146-000-149-112.}}
\date{This Version: June 22, 2015}

\begin{document}

\maketitle

\setlength{\parskip}{0.5em}

\begin{abstract}
For the problems of low-rank matrix completion, the efficiency of the widely-used nuclear norm technique may be challenged under many circumstances, especially when certain basis coefficients are fixed, for example, the low-rank correlation matrix completion in various fields such as the financial market and the low-rank density matrix completion from the quantum state tomography. To seek a solution of  high recovery quality beyond the reach of the nuclear norm, in this paper, we propose a rank-corrected procedure using a nuclear semi-norm to generate a new estimator. For this new estimator, we establish a non-asymptotic recovery error bound. More importantly, we quantify the reduction of the recovery error bound for  this rank-corrected procedure. Compared with the one obtained  for the nuclear norm penalized least squares estimator, this reduction can be substantial (around $50\%$).  We also provide necessary and sufficient conditions for rank consistency in the sense of Bach (2008). %\cite{Bac08},
Very interestingly, these conditions  are highly related  to the concept of constraint
nondegeneracy in matrix optimization. As a \hyphenation{byproduct}byproduct, our \hyphenation{results}results provide a theoretical foundation for the majorized penalty method of Gao and Sun (2010) and Gao (2010)
%\cite{GaoS10} and Gao \cite{Gao10}
for  structured low-rank matrix optimization problems. Extensive numerical experiments demonstrate that our proposed rank-corrected procedure can simultaneously achieve a high recovery accuracy and capture the low-rank structure.
%In this paper, we address low-rank matrix completion problems with fixed basis coefficients, which include the low-rank correlation matrix completion in various fields such as the financial market and the low-rank density matrix completion from the quantum state tomography. For this class of problems,
% the efficiency of the common nuclear norm penalized estimator for recovery may be challenged. Here, with a reasonable initial estimator, we propose a rank-corrected procedure to generate
% an estimator of high accuracy and low rank. For this new estimator, we establish a non-asymptotic recovery error bound and demonstrate the power of the rank-correction term on improving the recoverability. We also provide necessary
% and sufficient conditions for rank consistency in the sense of Bach \cite{Bac08}, in which the concept of constraint
% nondegeneracy in matrix optimization plays an important role. As a byproduct, our results provide a theoretical foundation for the majorized penalty method of Gao and Sun
% \cite{GaoS10} and Gao \cite{Gao10} for  structured low-rank matrix optimization problems.
 %for the positive semidefinite case, our results provide a theoretical foundation for the reweighted trace norm minimization advocated by Fazel et al. \cite{FazHB03, MohF10} for the matrix rank minimization problem, and the majorized penalty method of Gao and Sun \cite{GaoS10} and Gao \cite{Gao10} for the structured low-rank matrix optimization problem.
 \\

\end{abstract}

\section{Introduction}

 The low-rank matrix completion is to recover an unknown low-rank matrix
 from the under-sampled observations with or without noises. This problem is of considerable interest in
 many application areas, from machine learning to quantum state tomography.
 A basic idea to address a low-rank matrix completion problem is to minimize the rank of
 a matrix subject to certain constraints from observations. Since the direct minimization of rank function
 is generally NP-hard, a widely-used convex relaxation approach is to replace the rank function with the nuclear norm --- the convex envelope of the rank function over a unit ball of the spectral norm \cite{Faz02}.

 The nuclear norm technique has been observed to provide a low-rank solution in practice for a long time (see, e.g., \cite{MesP97,Mes98,Faz02}). The first remarkable theoretical characterization for the minimum rank solution via the nuclear norm minimization was
 given by Recht, Fazel and Parrilo \cite{RecFP10}, with the help of the concept of Restricted Isometric Property (RIP).
 Recognizing that the matrix completion problem does not obey the RIP, Cand{\`e}s and Recht \cite{CanR09} introduced the concept of incoherence property and proved that most low-rank matrices can be exactly recovered from a surprisingly small number of noiseless observations of randomly sampled entries via the nuclear norm minimization. The bound of the number of sampled entries was later improved to be near-optimal by Cand{\`e}s and Tao \cite{CanT10} through a counting argument. Such a bound was also
 obtained by Keshavan et al. \cite{KesMO10} for their proposed OptSpace algorithm. Later, Gross \cite{Gro11} sharpened the bound by employing a novel technique from quantum information theory developed in \cite{GroLFBE10}, in which noiseless observations were extended from entries to coefficients relative to an arbitrary basis. This technique was also adapted by Recht \cite{Rec11}, leading to a short and intelligible analysis. Besides the above results for the noiseless case, matrix completion with noise was first addressed by
 Cand{\`e}s and Plan \cite{CanP10}. More recently, nuclear norm penalized estimators for matrix completion with noise have been well studied by Koltchinskii, Lounici and Tsybakov \cite{KolLT11}, Negahban and Wainwright \cite{NegW12}, and Klopp \cite{Klo12} under different settings.  Besides the nuclear norm, estimators with other penalties  for matrix completion have also been considered in terms of recoverability in the literature, e.g., \cite{RohT11,Klo11,Kol12,SreRJ05,FoyS11}.

%
% In an alternative line, Keshavan, Montanari and Oh \cite{KesMO10, KesMO10_2} designed the OptSpace algorithm
% for matrix completion based on singular value thresholding and studied its behavior for both noiseless and noisy cases.
%  Moreover, Srebro and colleagues \cite{SreRJ05, SreS05, FoyS11} studied another rank-promoting function
% called the max-norm, particularly for collaborative filtering problems.
%
%

 The nuclear norm technique has been demonstrated to be a successful approach to encourage a low-rank solution for matrix completion.
 However, its efficiency may be challenged in some circumstances. For example, Salakhutdinov and Srebro \cite{SalS10} showed that when certain rows and/or columns are sampled with high probability, the nuclear norm minimization may fail in the sense that the number of observations required for recovery is much more than the setting of most matrix completion problems. It means that the efficiency of the nuclear norm techniques could be highly weakened under a general sampling scheme. Negahban and Wainwright \cite{NegW12} also pointed out the impact of such heavy sampling schemes on the recovery error bound. As a remedy for this, a weighted nuclear norm (trace norm), based on row- and column-marginals of the sampling distribution, was suggested in \cite{NegW12, SalS10, FoySSS11} if the prior information on sampling distribution is available. Moreover, the conditions characterized by Bach \cite{Bac08} for rank consistency of the nuclear norm penalized least squares estimator may not be satisfied, especially when certain constraints are involved.

%When the true matrix possesses a symmetric/Hermitian positive semidefinite structure, the impact of general sampling schemes on the recoverability of the nuclear norm technique is %more remarkable. In this situation, the nuclear norm reduces to the trace and thus only depends on diagonal entries rather than all entries as the rank function does. As a result, if %diagonal entries are heavily sampled, the low-rank promoting ability of minimizing the nuclear norm, as well as the recoverability, will be highly weakened. This phenomenon is fully %reflected in the widely-used correlation matrix completion problem, for which the nuclear norm becomes a constant and severely loses effectiveness.
A concrete example of interest is to recover a density matrix of a quantum system from Pauli measurements in quantum state tomography (see, e.g., \cite{GroLFBE10, FlaGLE12, Wan12}). A density matrix is a Hermitian positive semidefinite matrix of trace one. Clearly, if the constraints of positive semidefiniteness and trace one are simultaneously imposed on the nuclear norm minimization, the nuclear norm completely fails in promoting a low-rank solution. Thus, one of the two constraints has to be abandoned in the nuclear norm minimization and then be restored in the post-processing stage. In fact, this idea has been much explored in \cite{GroLFBE10,FlaGLE12} and the numerical results there indicated its relative efficiency though it still has much room for improvement.

All the above examples motivate us to ask whether it is possible to go beyond the nuclear norm approach for practical use to seek for better performance in low-rank matrix completion. In this paper, we provide a positive answer to this question with both theoretical and empirical supports. We first establish a unified low-rank matrix completion model, which allows for the imposition of fixed basis coefficients so that the correlation and the density matrix completion are included as special cases. It means that in our setting, for any given basis of the matrix space, a few basis coefficients of the true matrix are assumed to be fixed due to a certain structure or some prior information, and the rest are allowed to be observed with noises under a general sampling scheme. To pursue a low-rank solution with a high recovery accuracy, we propose a rank-correction step to generate a new estimator. The rank-correction step solves a penalized least squares problem with its penalization being the nuclear norm minus a linear rank-correction term constructed on a reasonable initial estimator. A satisfactory choice of the initial estimator could be the nuclear norm penalized least squares estimator or one of its analogies. The resulting convex matrix optimization problem can be solved by the efficient algorithms recently developed in
 \cite{FazPST12, JiaST12, JiaST12_1, JiaST12_2} even for large-scale cases.

 The idea of using a two-stage or even multi-stage procedure is not brand new for dealing with sparse recovery in the statistical and
 machine learning literature.  The $l_1$-norm penalized least squares method, also known as the Lasso \cite{Tib96}, is very attractive and
 popular for variable selection in statistics, thanks to the invention of the fast and efficient LARS algorithm \cite{EfrHJT04}.
 On the other hand, the $l_1$-norm penalty has long been known by statisticians to yield biased estimators and cannot achieve the best estimation performance \cite{FanL01, FanP04}. The issue of bias can be overcome by nonconvex penalization methods, see, e.g., \cite{LenLW06,Fan97,Zha10}. %The common penalties include the $l_q$-norm penalty ($0<q<1$) by Frank and Friedman \cite{LenLW06}, the smoothly clipped absolute deviation (SCAD) penalty by Fan \cite{Fan97}, and the minimax concave penalty (MCP) by Zhang \cite{Zha10}.
 A multi-stage procedure naturally occurs if the nonconvex problem obtained is solved by an iterative algorithm \cite{ZouL08, LaiXY13}. In particular, once a good initial estimator is used, a two-stage estimator is enough to achieve the desired asymptotic efficiency, e.g., the adaptive Lasso proposed by Zou \cite{Zou06}. There are also a number of important works along this line on variable selection, including \cite{LenLW06,MeiB06,ZhaY07,HuaMZ10,ZhoVB09, Mei07,FanL08}, to name only a few. For a broad overview, the interested readers are referred to the recent survey papers \cite{FanL10, FanLQ11}. It is natural to extend the ideas from the vector case to the matrix case. Fazel, Hindi and Boyd \cite{FazHB03} first proposed the reweighted trace minimization for minimizing the rank of a positive semidefinite matrix.  In \cite{Bac08}, Bach made an important step in extending the adaptive Lasso of Zou \cite{Zou06} to the matrix case for rank consistency. However, it is not clear how to apply Bach's idea to our matrix completion model with fixed basis coefficients since the required rate of convergence of the initial estimator for achieving asymptotic properties is no longer valid, as far as we can see. More critically, there are numerical difficulties in efficiently solving the resulting optimization problems. Numerical difficulties also occur in the reweighted nuclear norm approach proposed by Mohan and Fazel \cite{MohF10} as an extension of \cite{FazHB03} for rectangular matrices. Iterative reweighted least squares minimization is an alternative extension of \cite{FazHB03} independently proposed by Mohan and Fazel \cite{MohF13} and Fornasier, Rauhut and Ward \cite{ForRW11}, taking advantage of the property that the rank of a matrix is equal to the rank of the product of this matrix and its transpose. However, the resulting smoothness of inner-iteration subproblems is weak in encouraging a low-rank solution so much more iterations are needed in general and thus the computational cost is high especially when hard constraints such as fixed basis coefficients are involved.

 The rank-correction step to be proposed in this paper is for overcoming the above difficulties. This approach is inspired by the majorized penalty method   proposed by Gao and Sun \cite{GaoS10} for solving structured matrix optimization problems with a low-rank constraint.
 %Unlike \cite{GaoS10}, we here mainly consider the case where the rank of the true matrix is unknown.
 For our proposed rank-correction step, we establish a non-asymptotic recovery error bound in Frobenius norm, following a similar argument adopted by Klopp in \cite{Klo12}. We also discuss the impact of adding the rank-correction term on recovery error. More importantly, we provide an affirmative  guarantee that under mild condition the rank-correction step highly improves the recoverability, compared with the nuclear norm penalized least squares estimator.  As the estimator is expected to be of low-rank, we also study the asymptotic property --- rank consistency in the sense of Bach \cite{Bac08}, under the setting that the matrix size is assumed to be fixed. This setting may not be ideal for analyzing asymptotic properties for matrix completion, but it does allow us to take the crucial first step to gain insights into the limitation of the nuclear norm penalization. Among others, the concept of constraint nondegeneracy for conic optimization problem plays a key role in our analysis. Interestingly, our results of recovery error bound and rank consistency suggest a consistent criterion for constructing a suitable rank-correction function.
 %When the rank of the true matrix is known, this criterion can be satisfied by a natural choice, which coincidentally reduces our rank-correction step to an iteration in the majorized penalty method \cite{GaoS10}. For more common cases where the rank of the true matrix is unknown, since the natural choice is
 %no longer available, we propose two rank-correction functions to imitate its behavior.
 In particular, for the correlation and the density matrix completion problems, we prove that rank consistency automatically holds for a broad selection of rank-correction functions. For most cases, a single rank-correction step is sufficient for a substantial improvement, unless the sample ratio is rather low so that the rank-correction step may be iteratively used for two or three times to achieve the limit of improvement. Owing to this property, the advantage of our proposed method is more apparent in practical computations especially when fixed basis coefficients are involved. Finally, we remark that our results can also be used to provide a theoretical foundation in the statistical setting for the majorized penalty method of Gao and Sun
 \cite{GaoS10} and Gao \cite{Gao10} for  structured low-rank matrix optimization problems.
 %a theoretical foundation for the reweighted trace norm minimization advocated by Fazel et al. \cite{FazHB03, MohF10} for matrix rank minimization problems.

 %in the semidefinite case, the for the positive semidefinite case, the proposed rank-correction step is related to the reweighted trace norm minimization for matrix rank minimization problems, e.g., Fazel et.al. \cite{FazHB03, MohF10}. The key difference lies in the construction of the rank-correction function, which allows us to. However, our consideration is based on an entirely different perspective. The results presented in this paper may also be helpful to illustrate why the reweighted variations of nuclear norm minimization could achieve the goal in many situations.
%To the best of our knowledge, even when reduced to the vector case, our proposed rank-correction step still differs from existing methods in sparse recovery.

 This paper is organized as follows. In Section \ref{section2}, we introduce the observation model of matrix completion with fixed basis coefficients and formulate the rank-correction step. In Section \ref{section3}, we establish a non-asymptotic recovery error bound for the estimator generated from the rank-correction step and provide a quantification of the  improvement in recoverability. Section \ref{section4} provides necessary and sufficient conditions for rank consistency. Section \ref{section5} is devoted to the construction of the rank-correction function. In Section \ref{section6}, we report numerical results to validate the efficiency of our proposed rank-corrected procedure. We conclude this paper in Section \ref{section7}. All relevant material and all proofs of theorems are left in the appendices.

 \noindent {\bf Notation.} Here we provide a brief summary of the notation used in this paper.
 \begin{itemize}
   \item[$\bullet$\ ] Let $\mathbb{R}^{n_1\times n_2}$ and $\mathbb{C}^{n_1\times n_2}$ denote the space of all $n_1\times n_2$
                     real and complex matrices, respectively. Let $\mathcal{S}^n(\mathcal{S}_{+}^n,\,\mathcal{S}_{++}^n)$ denote the set of all $n\times n$ real symmetric (positive semidefinite, positive definite) matrices and $\mathcal{H}^n(\mathcal{H}_{+}^n,\,\mathcal{H}_{++}^n)$ denote the set of all $n\times n$ Hermitian (positive semidefinite, positive definite) matrices.

   \item[$\bullet$\ ] Let $\mathbb{V}^{n_1\times n_2}$ represent $\mathbb{R}^{n_1\times n_2}$, $\mathbb{C}^{n_1\times n_2}$, $\mathcal{S}^n$ or $\mathcal{H}^n$. We define  $n:=\min(n_1,n_2)$ for the previous two cases  and  stipulate $n_1=n_2=n$ for  the latter two cases. Let $\mathbb{V}^{n_1\times n_2}$ be endowed with the trace inner product $\langle \cdot, \cdot \rangle$ and its induced norm $\|\cdot\|_{F}$, i.e., $\langle X, Y \rangle:= \text{Re}\big(\text{Tr}(X^{\mathbb{T}}Y)\big)$ for $X, Y\in\mathbb{V}^{n_1\times n_2}$, where $``\text{Tr}"$ stands for the trace of a matrix and $``\text{Re}"$ means the real part of a complex number.

   \item[$\bullet$\ ] For the real case, i.e., $\mathbb{V}^{n_1\times n_2} = \mathbb{R}^{n_1\times n_2}$ or $\mathbb{V}^{n_1\times n_2} = \mathcal{S}^n$, let $\mathbb{S}^n\,(\mathbb{S}^n_+,\,\mathbb{S}^n_{++})$ represent $\mathcal{S}^n\,(\mathcal{S}^n_+,\,\mathcal{S}^n_{++})$; and for the complex case, i.e., $\mathbb{V}^{n_1\times n_2} = \mathbb{C}^{n_1\times n_2}$ or $\mathbb{V}^{n_1\times n_2} = \mathcal{H}^n$, let $\mathbb{S}^n\,(\mathbb{S}^n_+,\,\mathbb{S}^n_{++})$ represent $\mathcal{H}^n\,(\mathcal{H}^n_+,\,\mathcal{H}^n_{++})$.

   \item[$\bullet$\ ] For the real case, $\mathbb{O}^{n\times k}$ denotes the set of all $n\times k$ real matrices with orthonormal columns, and for the complex case, $\mathbb{O}^{n\times k}$ denotes the set of all $n\times k$ complex matrices with orthonormal columns. When $k=n$, we write $\mathbb{O}^{n\times k}$ as $\mathbb{O}^n$ for short.

   \item[$\bullet$\ ] The notation $^{\mathbb{T}}$ denotes the transpose for the real case and
                     the conjugate transpose for the complex case. The notation $^\ast$ means the adjoint of a linear operator.

   \item[$\bullet$\ ] For any index set $\pi$, let $|\pi|$ denote the cardinality of $\pi$, i.e., the number of elements in $\pi$. For any $x \in \mathbb{R}^n$, let $|x|$ denote the vector in $\mathbb{R}^n_+$ whose $i$-th component is $|x_i|$, let $x_+$ denote the vector in $\mathbb{R}^n_+$ whose $i$-th component is $\max(x_i, 0)$ and let $x_-$ denote the vector in $\mathbb{R}^n_+$ whose $i$-th component is $\min(-x_i, 0)$.

   \item[$\bullet$\ ] For any given vector $x$, $\text{Diag}(x)$ denotes a rectangular diagonal matrix of suitable size with the $i$-th diagonal entry being $x_i$.
   \item[$\bullet$\ ] For any $x \in \mathbb{R}^n$, let $\|x\|_2$ and $\|x\|_\infty$ denote
          the Euclidean norm and the maximum norm, respectively. For any $X\in \mathbb{V}^{n_1\times n_2}$, let $\|X\|$ and $\|X\|_*$ denote
          the spectral norm and the nuclear norm, respectively.
   \item[$\bullet$\ ] The notations $\stackrel{a.s.}{\rightarrow}$, $\stackrel{p}{\rightarrow}$ and $\stackrel{d}{\rightarrow}$ mean almost sure convergence, convergence in probability and convergence in distribution, respectively. We write $x_m = O_p(1)$ if $x_m$ is
       bounded in probability.
   \item[$\bullet$\ ] For any set $K$, let $\delta_{K}(x)$ denote the indicator function of $K$,
                       i.e., $\delta_{K}(x)=0$ if $x\in K$, and $\delta_{K}(x)=+\infty$ otherwise.
                       Let $I_n$  denote the $n\times n$ identity matrix.
 \end{itemize}

\section{Problem formulation}\label{section2}

 In this section, we formulate the model of the matrix completion problem with fixed basis coefficients,
 and then propose an adaptive nuclear semi-norm penalized least squares estimator for solving this class of problems.

 \subsection{The observation model} \label{subsecobs}

 Let $\{\Theta_1,\ldots,\Theta_d\}$ be a given orthonormal basis of the given real inner product space $\mathbb{V}^{n_1\times n_2}$.
 Then, any matrix $X\in\mathbb{V}^{n_1\times n_2}$ can be uniquely expressed in the form of $X=\sum_{k=1}^d \langle \Theta_k, X \rangle\Theta_k$,
 where $\langle \Theta_k, X\rangle$ is called the basis coefficient of $X$ relative to $\Theta_k$. Throughout this paper, let $\overline{X} \in \mathbb{V}^{n_1\times n_2}$
 be the unknown low-rank matrix to be recovered and let $\text{rank}(\overline{X}) = r$. In some practical applications, for example, the correlation and density matrix completion,
 a few basis coefficients of the unknown matrix $\overline{X}$ are fixed (or assumed to be fixed) due to a certain structure or
 reliable prior information. We let $\alpha\subseteq\{1,2,\ldots,d\}$ denote the set of the indices relative
 to which the basis coefficients are fixed, and $\beta$ denote the complement of $\alpha$ in $\{1,2,\ldots,d\}$,
 i.e., $\alpha\cap\beta=\emptyset$ and $\alpha\cup\beta = \{1,\ldots, d\}$. We define $d_1: = |\alpha|$ and $d_2 := |\beta|$.

When a few basis coefficients are fixed, one only needs to observe the rest for recovering the unknown matrix $\overline{X}$.
 Assume that we are given a collection of $m$ noisy observations of the basis coefficients relative to
 $\{\Theta_{k}:k\in\beta\}$ in the following form
 \begin{equation}\label{eqnobsori}
    y_i = \left\langle \Theta_{\omega_i}, \overline{X}\right\rangle + \nu \xi_i, \quad i = 1, \ldots, m,
 \end{equation}
 where $\omega_i$ are the indices randomly sampled from the index set $\beta$,
 $\xi_i$ are the independent and identically distributed (i.i.d.) noises with
 $\mathbb{E}(\xi_i)=0$ and $\mathbb{E}(\xi^2_i)=1$, and $\nu>0$ controls the magnitude of noise. Unless otherwise stated,
 we assume a general weighted sampling (with replacement) scheme with the sampling distributions of $\omega_i$ as follows.
 \begin{assumption}\label{asmpprob}
  The indices $\omega_1,\ldots, \omega_m$  are i.i.d. copies of a random variable $\omega$ that has a probability
  distribution $\Pi$ over $\{1,\ldots, d\}$ defined by
  $${\rm Pr}(\omega = k) =\left\{\begin{array}{ll}
                                0 & {\rm if}\ k\in \alpha,\\
                                p_k>0 & {\rm if}\ k\in \beta.
                          \end{array}\right.$$
 \end{assumption}
  Note that each $\Theta_k, k\in\beta$ is assumed to be sampled with a positive probability in this sampling scheme.
  In particular, when the sampling probability of all $k\in \beta$ are equal, i.e., $p_k =1/d_2 \ \forall\, k\in \beta$, we say that the observations are
  sampled uniformly at random.

For notational simplicity, let $\Omega$ be the multiset of all the sampled indices from the index set $\beta$, i.e., $\Omega =\{\omega_1, \ldots, \omega_m\}$.
 With a slight abuse on notation, we define the sampling operator $\mathcal{R}_\Omega$: $\mathbb{V}^{n_1\times n_2} \rightarrow \mathbb{R}^m$ associated with $\Omega$ by
  $$\mathcal{R}_\Omega(X) := \big(\langle \Theta_{\omega_1}, X \rangle, \ldots, \langle \Theta_{\omega_m}, X \rangle \big)^{\mathbb{T}}, \quad X \in \mathbb{V}^{n_1\times n_2}.$$
 Then, the observation model (\ref{eqnobsori}) can be expressed in the following vector form
  \begin{equation}\label{eqnobs}
    y = \mathcal{R}_\Omega(\overline{X})+\nu\xi,
  \end{equation}
  where $y = (y_1,\ldots,y_m)^{\mathbb{T}} \in \mathbb{R}^m$ and $\xi =\!(\xi_1,\ldots, \xi_m)^{\mathbb{T}} \in \mathbb{R}^m$ denote the observation vector
  and the noise vector, respectively.

  Next, we present some examples of low-rank matrix completion problems in the above settings.
\begin{description}
\item[(1)] {\bf Correlation matrix completion.} A correlation matrix is an $n\times n$ real symmetric or Hermitian positive semidefinite matrix with all diagonal entries being ones. Let $e_i$ be the vector with the $i$-th entry being one and the others being zeros. Then,
    $\langle e_ie_i^{\mathbb{T}}, \overline{X}\rangle =\overline{X}_{ii}=1 \ \forall\, 1\leq i\leq n$.
    The recovery of a correlation matrix is based on the observations of entries.
    For the real case, $\mathbb{V}^{n_1\times n_2} = \mathcal{S}^n$, $d = n(n+1)/2$, $d_1=n$,
    $$\Theta_{\alpha} = \big\{e_ie_i^{\mathbb{T}} \ | \ 1\leq i \leq n\big\}  \quad \text{and} \quad \Theta_\beta = \left\{\frac{1}{\sqrt{2}}(e_ie_j^{\mathbb{T}}+e_je_i^{\mathbb{T}})  \ \Big| \ 1\leq i < j \leq n\right\};$$
    and for the complex case, $\mathbb{V}^{n_1\times n_2} = \mathcal{H}^n$, $d=n^2$, $d_1 =n$,
    $$\Theta_\alpha =\!\big\{e_ie_i^{\mathbb{T}} \ | \ 1\leq i \leq n\big\}  \ \ \text{and} \ \
    \Theta_\beta\!=\!\left\{\!\frac{1}{\sqrt{2}}(e_ie_j^{\mathbb{T}}\!+e_je_i^{\mathbb{T}}),\frac{\sqrt{-1}}{\sqrt{2}}(e_ie_j^{\mathbb{T}}\!-e_je_i^{\mathbb{T}})
    \ \Big| \ i<j\!\right\}.$$
    Here, $\sqrt{-1}$ represents the imaginary unit. Of course, one may fix some off-diagonal entries in specific applications.
\item[(2)] {\bf Density matrix completion.} A density matrix of dimension $n=2^l$ for some positive integer $l$ is
    an $n\times n$ Hermitian positive semidefinite matrix with trace one.  In quantum state tomography, one aims to recover a density matrix from Pauli measurements (observations of the coefficients relative to the Pauli basis) \cite{GroLFBE10, FlaGLE12}, given by
    $$\Theta_\alpha = \left\{\frac{1}{\sqrt{n}}I_n\right\} \  \text{and} \  \Theta_\beta = \left\{\frac{1}{\sqrt{n}} (\sigma_{s_1}\otimes \cdots \otimes \sigma_{s_l}) \ \Big | \  (s_1,\ldots,s_l) \in \{0,1,2,3\}^l\right\}\Big \backslash \Theta_\alpha,$$
    where ``$\otimes$'' means the Kronecker product of two matrices and $$\sigma_0 = \begin{pmatrix} 1 & 0 \\ 0 & 1 \end{pmatrix},\ \sigma_1 = \begin{pmatrix} 0 & 1 \\ 1 & 0 \end{pmatrix}, \ \sigma_2 = \begin{pmatrix} 0 & -\sqrt{-1} \\ \sqrt{-1} & 0 \end{pmatrix}, \ \sigma_3 = \begin{pmatrix} 1 & 0 \\ 0 & -1 \end{pmatrix}$$
    are the Pauli matrices.  In this setting, $\mathbb{V}^{n_1\times n_2} = \mathcal{H}^{n}$, $\text{Tr}(\overline{X}) = \langle I_n, \overline{X} \rangle = 1$, $d = n^2$, and  $d_1 = 1$.

\item[(3)] {\bf Rectangular matrix completion.} Assume that a few entries of a rectangular matrix are known and let $\mathcal{I}$ be the index set of these entries. One aims to recover this rectangular matrix from the observations of the rest entries. For the real case, $\mathbb{V}^{n_1\times n_2} = \mathbb{R}^{n_1\times n_2}$, $d= n_1n_2$, $d_1 = |\mathcal{I}|$,
    $$\Theta_\alpha = \big\{e_ie_j^{\mathbb{T}} \ | \ (i,j) \in \mathcal{I}\big\} \quad \text{and} \quad \Theta_\beta = \big\{e_ie_j^{\mathbb{T}} \ | \ (i,j) \notin \mathcal{I}\big\};$$
    and for the complex case, $\mathbb{V}^{n_1\times n_2} = \mathbb{C}^{n_1\times n_2}$, $d= 2n_1n_2$,  $d_1 = 2|\mathcal{I}|$,
    $$\Theta_\alpha = \big\{e_ie_j^{\mathbb{T}}, \sqrt{-1}e_ie_j^{\mathbb{T}} \ | \ (i,j) \in \mathcal{I}\big\} \quad \text{and} \quad \Theta_\beta = \big\{e_ie_j^{\mathbb{T}}, \sqrt{-1}e_ie_j^{\mathbb{T}} \ | \ (i,j) \notin \mathcal{I}\big\}.$$
 \end{description}

Now we introduce some linear operators that are frequently used in the subsequent sections.
For any given index set $\pi \subseteq \{1,\ldots,d\}$, say $\alpha$ or $\beta$,
we define the linear operators $\mathcal{R}_\pi$: $\mathbb{V}^{n_1\times n_2} \rightarrow \mathbb{R}^{|\pi|}$, $\mathcal{P}_\pi$: $\mathbb{V}^{n_1\times n_2} \rightarrow \mathbb{V}^{n_1\times n_2}$ and $\mathcal{Q}_\pi$: $\mathbb{V}^{n_1\times n_2} \rightarrow \mathbb{V}^{n_1\times n_2}$ respectively, by
\begin{equation*}\label{operator-R-P-Q}
\mathcal{R}_\pi(X):= \big(\langle \Theta_k, X \rangle \big)^{\mathbb{T}}_{k\in\pi}, \quad
\mathcal{P}_\pi(X) := \sum_{k\in \pi} \langle \Theta_k, X \rangle \Theta_k \quad \text{and} \quad \mathcal{Q}_\pi(X): = \sum_{k\in \pi} p_k \langle \Theta_k, X \rangle \Theta_k.
\end{equation*}

For convenience of discussions, in the rest of this paper, for any given $X\in\mathbb{V}^{n_1\times n_2}$,
 we denote by $\sigma(X)=\big(\sigma_1(X), \ldots, \sigma_n(X)\big)^{\mathbb{T}}$ the singular value vector of $X$ arranged in the nonincreasing order
 and define
 $$\mathbb{O}^{n_1,n_2}(X) :=\big\{(U,V) \in \mathbb{O}^{n_1}\times \mathbb{O}^{n_2}\mid X = U\text{Diag}\big(\sigma(X)\big)V^\mathbb{T}\big\}.$$ In particular, when $\mathbb{V}^{n_1\times n_2} = \mathbb{S}^n$, we denote by $\lambda(X)=\big(\lambda_1(X), \ldots, \lambda_n(X)\big)^{\mathbb{T}} $ the eigenvalue vector of $X$ with $|\lambda_1(X)|\geq \ldots \geq |\lambda_n(X)|$ and define  $$\mathbb{O}^n (X) :=\big\{P \in \mathbb{O}^n \mid X = P\text{Diag}(\lambda(X))P^\mathbb{T}\big\}.$$
For any $X\in\mathbb{V}^{n_1\times n_2}$ and any $(U,V)\in \mathbb{O}^{n_1,n_2}(X)$, we write $U = [U_1 \  U_2 ]$
and $V =[V_1 \  V_2 ]$ with
$U_1 \in \mathbb{O}^{n_1\times r}$, $U_2 \in \mathbb{O}^{n_1\times (n_1-r)}$, $V_1 \in \mathbb{O}^{n_2\times r}$ and $V_2 \in \mathbb{O}^{n_2\times (n_2-r)}$. In particular, for any $X\in\mathbb{S}_+^n$ and any $P \in \mathbb{O}^n(X)$, we write $P = [P_1 \ P_2]$ with $P_1 \in \mathbb{O}^{n \times r}$ and $P_2 \in \mathbb{O}^{n\times (n-r)}$.

\subsection{The rank-correction step}
In many situations, the nuclear norm penalization performs well for matrix recovery,
but its efficiency may be challenged if the observations are sampled at random obeying a general distribution
 such as the one considered in \cite{SalS10}. The setting of fixed basis coefficients in our matrix completion model can also be regarded to be under an extreme sampling scheme. In particular, for the correlation and density matrix completion,
 the nuclear norm completely loses its efficiency since it reduces to a constant  in these two cases. In order to overcome
 the shortcomings of the nuclear norm penalization, we propose a rank-correction step to generate an estimator in pursuit of a better recovery performance.

 Recall that $\overline{X}$ is the unknown true matrix of rank $r$. Given an initial estimator $\widetilde{X}_m$ of $\overline{X}$, say, the nuclear norm penalized least squares estimator or one of its analogies, our proposed rank-correction step is to solve the convex optimization problem
\begin{equation}\label{eqnrcs}
\begin{aligned}
\widehat{X}_m \ \in \ \mathop{\arg\min}_{X \in \mathbb{V}^{n_1\times n_2}} & \ \ \frac{1}{2m} \left\|y - \mathcal{R}_\Omega(X)\right\|_2^2 + \rho_m \big(\|X\|_* - \langle F(\widetilde{X}_m),X\rangle \big)\\
\text{s.t.} & \ \ \mathcal{R}_\alpha(X) = \mathcal{R}_\alpha(\overline{X}), \ \ \|\mathcal{R}_\beta(X)\|_\infty \leq b,\ \ X \in \mathcal{C},
\end{aligned}
\end{equation}
where $\rho_m > 0$ is the penalty parameter (depending on the number of observations), $b$ is an upper bound of the magnitudes of basis coefficients of $\overline{X}$, $\mathcal{C} \subseteq \mathbb{V}^{n_1\times n_2}$ is a closed convex set that contains $\overline{X}$, and $F:\mathbb{V}^{n_1\times n_2} \rightarrow \mathbb{V}^{n_1\times n_2}$ is a spectral operator associated with a symmetric function $f:\mathbb{R}^n\rightarrow \mathbb{R}^n$. One may refer to \ref{AppendixSpeOpe} for more information on the concept of spectral operators. (Indeed, based on the subsequent analysis for better recovery performance, the choice $f:\mathbb{R}^n \rightarrow [0,1]^n$ is much preferred, for which the penalization $\|X\|_* - \langle F(\widetilde{X}_m), X\rangle$ is indeed a nuclear semi-norm. But this choice criterion is not compulsory). The bound restriction is very mild since such a bound is often available in applications, for example, the correlation and the density matrix completion. This boundedness setting can also be found in previous works done by Negahban and Wainwright \cite{NegW12} and Klopp \cite{Klo12}.

Hereafter, we call $F$ the rank-correction function and $\langle F(\widetilde{X}_m), X \rangle$ the rank-correction term.
Note that, when $F \equiv 0$, the rank-correction step (\ref{eqnrcs}) reduces to the nuclear norm penalized least squares estimator, which equally penalizes singular values to promote a low-rank solution for matrix completion. Certainly, for this purpose, penalizing more on small singular values or even directly penalizing the rank function could serve better, but only theoretically rather than practically, due to the lack of convexity. Also note that an initial estimation, if deviates not too much from the true matrix, could contain some information of the singular values and/or the rank of the true matrix to a certain extent. Therefore, provided such an initial estimator is available, it is achievable to construct a rank-correction term with a suitable $F$ to substantially offset the penalization of large singular values from the nuclear norm penalty. Consequently, we can expect the rank-correction step (\ref{eqnrcs}) to have a better low-rank promoting ability and outperform the nuclear norm penalized least squares estimator.

The key issue is then how to construct a favored rank-correction function $F$. In the next two sections, we provide theoretical supports to our proposed rank-correction step, from which some important guidelines on the construction of $F$ can be captured. In particular, if one chooses the nuclear norm penalized least squares estimator to be the initial estimator $\widetilde{X}_m$, and also suitably chooses the spectral operator $F$ so that $\|X\|_* - \langle F(\widetilde{X}_m), X \rangle$ is a semi-norm, called nuclear semi-norm, then the estimator $\widehat{X}_m$ generated from this two-stage procedure is called the adaptive nuclear semi-norm penalized least squares estimator associated with $F$.

\subsection{Relation with the majorized penalty approach}

The rank-correction step above is inspired by the majorized penalty approach   proposed by Gao and Sun \cite{GaoS10}
for solving the rank constrained matrix optimization problem:
\begin{equation}\label{major1}
\min_{X \in \mathcal{C}} \big\{h(X):\ {\rm rank}(X) \leq r\big\},
\end{equation}
where $r \geq 1$, $h:\mathbb{V}^{n_1\times n_2} \rightarrow \mathbb{R}$ is a given continuous function
and $\mathcal{C} \in \mathbb{V}^{n_1\times n_2}$ is a closed convex set.
Note that for any $X \in \mathbb{V}^{n_1\times n_2}$, the constraint $\text{rank}(X)\leq r$ is equivalent to
$$ 0=\sigma_{r+1}(X)+\cdots + \sigma_n(X) = \|X\|_* - \|X\|_{(r)},$$
where $\|X\|_{(r)}:=\sigma_1(X)+\cdots + \sigma_r(X)$ denotes the Ky Fan $r$-norm.
The central idea of the majorized penalty approach is to solve the following penalized version of (\ref{major1}):
 $$ \min_{X \in \mathcal{C}}\ h(X) + \rho\big(\|X\|_*-\|X\|_{(r)}\big),$$
 where $\rho>0$ is the penalty parameter. With the current iterate $X^k$, the majorized penalty approach
 yields the next iterate $X^{k+1}$ by solving the convex optimization problem
 \begin{equation}\label{eqnmpa}
  \min_{X\in\mathcal{C}}\ \widehat{h}^k(X) + \rho\big(\|X\|_* - \langle G^k, X \rangle\big),
 \end{equation}
 where $G^k$ is a subgradient of the convex function $\|X\|_{(r)}$ at $X^k$, and $\widehat{h}^k$ is
 a convex majorization function of $h$ at $X^k$. By comparing with (\ref{eqnrcs}), one may notice that our proposed rank-correction step is close to a single step of the majorized penalty approach.

Note that the rank constrained least squares problem is of great consideration in matrix completion especially when the rank information is known. However, different from the noiseless case, for matrix completion with noise, the solution to the rank constrained least squares problem (assuming the uniqueness) is in general not the true matrix though quite close to it. Indeed, there may exist many candidate matrices surrounding the true matrix and having its rank. The rank constrained least squares solution is only one of them. It deviates the least from the noisy observations rather than the true matrix. Naturally, it is conceivable that some candidate matrices may deviate a bit more from the noisy observations but less from the true matrix. So, for the purpose of matrix completion, there is no need to aim precisely at the rank constrained least squares solution and find this solution accurately. An approach roughly towards it such as our proposed rank-correction step (\ref{eqnrcs}) is good enough to bring similar good recovery performance.

 \section{Error bounds}\label{section3}

In this section, we aim to derive a recovery error bound in Frobenius norm for the estimator generated from the rank-correction step (\ref{eqnrcs}) and discuss the impact of the rank-correction term on the resulting bound. The analysis mainly  follows Klopp's arguments in \cite{Klo12}, which is also
in line with those used by  Negahban and Wainwright \cite{NegW12}.
%The results are also applicable if more constraints are added to the estimator (\ref{eqnrcs}), e.g., for the positive
%semidefinite case, since adding more prior information does not have any harm in recoverability but only improvement.

We start the analysis by defining a quantity, which plays a key role in the subsequent analysis, as
\begin{equation}\label{eqndelalpbet}
a_m := \frac{1}{\sqrt{r}} \|F(\widetilde{X}_m) - \overline{U}_1 \overline{V}_1^{\mathbb{T}}\|_F.
\end{equation}
A basic relation between the true matrix $\overline{X}$ and its estimate $\widehat{X}_m$ can be obtained by using the optimality of $\widehat{X}_m$ to the problem (\ref{eqnrcs}) as follows.

\begin{theorem}\label{thmopbd}
For any $\kappa > 1$, if
$\rho_m \geq \kappa \nu \Big\| \frac{1}{m} \mathcal{R}_\Omega^*(\xi) \Big\|,$
then the following inequality holds:
\begin{equation}\label{eqnopbd}
\frac{1}{2m} \big\|\mathcal{R}_\Omega(\widehat{X}_m -\overline{X})\big\|_2^2 \leq\bigg(\! \frac{\sqrt{2}}{\kappa} + a_m\!\bigg)\rho_m \sqrt{r}\|\widehat{X}_m\!-\overline{X}\|_F.
\end{equation}
\end{theorem}

We emphasize  that $\kappa$ is not restricted to be a constant in Theorem \ref{thmopbd} but could be set to depend on the size of matrix. This realization is important as can be seen in the sequel. According to Theorem \ref{thmopbd}, the choice of the penalty parameter $\rho_m$ depends on the observation noises $\xi_i$ and the sampling operator $\mathcal{R}_\Omega$. Therefore, we make the following  assumption on the noises $\xi_i$ as follows:
\begin{assumption}\label{asmpnoi}
The i.i.d. noise variables $\xi_i$ are sub-exponential, i.e., there exist positive constants $c_1$, $c_2$ and $c_3$ such that for all $t>0$, ${\rm Pr}(|\xi_i| \geq t)\leq c_1\exp(-c_2t^{c_3}).$
\end{assumption}
Moreover, based on Assumption \ref{asmpprob}, we further define quantities $\mu_1$ and $\mu_2$ that control the sampling probability for observations as
\begin{equation}\label{eqndefL}
\mu_1\geq \frac{1}{d_2} \cdot \max_{k\in\beta} \left\{\frac{1}{p_k}\right\} \quad \text{and}  \quad \mu_2 \geq \sqrt{d_2}\cdot \max\Bigg\{\Bigg\|\sum_{k\in\beta} p_k\Theta_k\Theta_k^\mathbb{T}\Bigg\|,\ \Bigg\|\sum_{k\in\beta} p_k\Theta_k^\mathbb{T}\Theta_k\Bigg\|\Bigg\}.
\end{equation}
It is easy to obtain that  $\mu_1 \geq 1$ and $\mu_2 \geq 1$, according to the facts $\sum_{k\in\beta} p_k = 1$ and $\text{Tr}\big(\sum_{k\in\beta} p_k\Theta_k\Theta_k^\mathbb{T}\big) = \text{Tr}\big(\sum_{k\in\beta} p_k\Theta_k^\mathbb{T}\Theta_k\big)=1$, respectively. In general, the values of $\mu_1$ and $\mu_2$ depend on the sampling distribution. The more extreme the sampling distribution is, the larger these two values have to be. Assume that there exist some positive constants $\gamma_1$ and $\gamma_2$ such that $\gamma_1 / d_2 \leq  p_k \leq \gamma_2 / d_2$, $\forall\, k \in \beta$. Then we can easily set $\mu_1:=1/\gamma_1$. The setting of $\mu_2$ is not universal for different cases. For example, consider the cases described in Section \ref{section2}.  For correlation matrix completion, we can set $\mu_2:=\gamma_2/\sqrt{2}$ for the real case and $\mu_2:=\gamma_2$ for the complex case. For density matrix completion, we can set $\mu_2:=1$ for any sampling distribution. For rectangular matrix completion, we can set $\mu_2 : =\gamma_2$ for the real case and $\mu_2:=\sqrt{2}\gamma_2$ for the complex case. Note that $\gamma_1=\gamma_2=1$ for uniform sampling.

Theorem \ref{thmopbd} reveals the key to deriving a recovery error bound in Frobenius norm, that is, to establish the relation between $\frac{1}{m} \|\mathcal{R}_\Omega(\widehat{X}_m-\overline{X})\|_{2}^2$ and $\|\widehat{X}_m-\overline{X}\|_F^2$. This can be achieved by looking into some RIP-like property of the sampling operator $\mathcal{R}_\Omega$, as done previously in \cite{NegW12, KolLT11, Klo12, Liu11}. Following this idea, we obtain an explicit recovery error bound as follows:

\begin{theorem}\label{thmstobd}
Under Assumptions \ref{asmpprob} and \ref{asmpnoi},
there exist  some positive absolute constants $c_0, c_1, c_2, c_3$ and some positive constants $C_0, C_1$ (only depending on the $\psi_1$ Orlicz norm of $\xi_k$) such that when $m \geq c_3 \sqrt{d_2}\log^3(n_1+n_2)/\mu_2$, for any $\kappa>1$, if $\rho_m$ is chosen as
\begin{equation}\label{bestrho}
\rho_m = C_1\kappa \nu\sqrt{\frac{\mu_2\log(n_1+n_2)}{\sqrt{d_2}m}},
\end{equation}
then with probability at least $1-c_1(n_1+n_2)^{-c_2}$,
\begin{align}\label{eqnstobd}
 \frac{\|\widehat{X}_m\!-\!\overline{X}\|_F^2}{d_2}\!\leq \!C_0 \bigg(\!{c_0}^2\big(\!\sqrt{2}\!+\!\kappa a_m\big)^2\!\nu^2\!+\!\Big(\!\frac{\kappa}{\kappa\!-\!1}\!\Big)^{\!2}\! \big(\!\sqrt{2}\!+\!a_m\big)^2 b^2\!\bigg) \mu_1^2 \mu_2 \!\frac{\sqrt{d_2}r\log(n_1\!+\!n_2)}{m}.
\end{align}
\end{theorem}

Theorem \ref{thmstobd} shows that for any rank-correction function $F$, controlling the recovery error only needs the samples size $m$ to be of   roughly the degree of freedom of a rank $r$ matrix up to a logarithmic factor in the matrix size. Besides the information on the order of magnitude, Theorem \ref{thmstobd} also provides us more details on the constant part in the recovery error bound, which also plays an important role in practice. The impact of different choices of rank-correction functions on recovery error is fully embodied with the value of $a_m$. Note that the smaller $a_m$ is, the smaller the error bound (\ref{eqnstobd}) is for a fixed $\kappa$, and thus the smaller value this error bound can achieve for the best $\kappa$ (as well as the best $\rho_m$). Therefore, we aim to establish an explicit relationship between $a_m$ and $F$ in the next theorem.

\begin{theorem}\label{thmambound}
For any given $\widetilde{X}_m \in \mathbb{V}^{n_1\times n_2}$ such that $\|\widetilde{X}_m - \overline{X}\|_F/\sigma_r(\overline{X}) < 1/2$, we have
$$a_m \leq -\frac{1}{\sqrt{2r}} \log\bigg(1- \sqrt{2} \, \frac{\|\widetilde{X}_m - \overline{X}\|_F}{\sigma_r(\overline{X})}\bigg) + \varepsilon_F(\widetilde{X}_m),$$
where $\varepsilon_{F}(\widetilde{X}_m):= \frac{1}{\sqrt{r}}\|F(\widetilde{X}_m) - \widetilde{U}_{m,1}\widetilde{V}_{m,1}^{\mathbb{T}}\|_F$.
\end{theorem}

It is immediate from Theorem \ref{thmambound} that
\begin{equation}\label{eqnamsmallcond}
\frac{\|\widetilde{X}_m - \overline{X}\|_F}{\sigma_r(\overline{X})} < \frac{1}{\sqrt{2}}\Big(1-e^{-\sqrt{2r}(1-\varepsilon_F(\widetilde{X}_m))}\Big)  \qquad \Longrightarrow \qquad a_m < 1.
\end{equation}
Recall that the nuclear norm penalized least squares estimator corresponds to the rank-correction step with $F \equiv 0$ so that $a_m = 1$. Therefore, Theorem \ref{thmambound}  guarantees that if the initial estimator $\widetilde{X}_m$ does not deviate too much from $\overline{X}$,  the rank-correction step outperforms the nuclear norm penalized least squares estimator in the sense of recovery error, provided that $F(\widetilde{X}_m)$ is close to $\widetilde{U}_{m,1}\widetilde{V}_{m,1}^{\mathbb{T}}$. For example, consider the case when the rank of the true matrix is known. One may simply choose $F(X) = U_1V_1^{\mathbb{T}}$ to take advantage of the rank information. In this case, the requirement in (\ref{eqnamsmallcond}) ensuring  $a_m < 1$ simply reduces to $\frac{\|\widetilde{X}_m - \overline{X}\|_F}{\sigma_r(\overline{X})} < 0.535 <  \frac{1}{\sqrt{2}}(1-e^{-\sqrt{2r}})$. Moreover, further suppose that $\widetilde{X}_m$ is the nuclear norm penalized least squares estimator. Then, according to Theorems \ref{thmstobd} and \ref{thmambound},  one only needs samples with size
$$ m = O\bigg(\sqrt{d_2}r^2\log^{1+2\tau}(n_1+n_2) \cdot \frac{d_2}{\sigma_r^2(\overline{X})}\bigg) \quad \Longrightarrow \quad a_m = O(\log^{-\tau}(n_1+n_2)),$$
where $\tau > 0$. As can be seen, the larger the matrix size $n$ is, the easier $a_m$ becomes less than $1$ or even close to $0$. If the rank of the true matrix is unknown, one could construct the rank-correction function $F$ on account of the tradeoff between optimality and robustness, to be discussed in Section \ref{section5}. An experimental example of the relationship between $a_m$ and $F$ can be found in Table \ref{tab1}.

Next, we demonstrate the power of the rank-correction term with more details. It is interesting to notice that the value of $\kappa$ (as well as $\rho_m$) has a substantial impact on the recovery error bound (\ref{eqnstobd}). The part related to the magnitude of noise $\nu$ increases as $\kappa$ increases, while the part related to the upper bound $b$ of entries slightly decreases to its limit as $\kappa$ increases. Therefore, our first target is to find the smallest error bound in terms of (\ref{eqnstobd}) among all possible $\kappa >1$.   It is possible to work on the error bound (\ref{eqnstobd}) directly for its minimum in $\kappa$ but the subsequent analysis is much more tedious.  For simplicity of illustration, instead, we perform our analysis on a slightly relaxed version instead as
\begin{align*}
\frac{\|\widehat{X}_m- \overline{X}\|_F^2}{d_2}  \leq
C_0 \, \eta_m^2 \, \mu_1^2\,\mu_2 \, \frac{\sqrt{d_2}r\log(n_1 + n_2)}{mn},
\end{align*}
where
$$\eta_m := c_0\big(\sqrt{2} + \kappa a_m \big) \nu + \bigg(\frac{\kappa}{\kappa - 1}\bigg) \big(\sqrt{2} + a_m\big)b.$$
Direct calculation shows that over $\kappa > 1$, $\eta_m$ attains its minimum
$$\overline{\eta}_m = \big(\sqrt{2}+a_m\big)(c_0\nu + b) + 2\sqrt{a_m\big(\sqrt{2}+a_m\big)c_0\nu b}\quad \text{at} \quad \overline{\kappa} = 1+\sqrt{\bigg(1+\frac{\sqrt{2}}{a_m}\bigg)\frac{b}{c_0\nu}}.$$
It is worthwhile to note that $\overline{\kappa} = O\big(1/\sqrt{a_m}\,\big)$ when $a_m \ll 1$,
meaning that the optimal choice of $\kappa$ is inversely proportional to  $\sqrt{a_m}$ rather than a simple constant. (This observation is important for achieving the rank consistency in Section \ref{section4}.)
In other words, for achieving the best possible recovery error, the penalty parameter $\rho_m$ chosen for the rank-correction step (\ref{eqnrcs}) with $a_m <1$ should be larger than that for the nuclear norm penalized least squares estimator.
In addition, consider two extreme cases with $a_m = 1$ and $a_m = 0$ respectively:
$$\overline{\eta}_m = \begin{cases}
\overline{\eta}^0 := \sqrt{2} (c_0 \nu + b) & \quad \text{if} \ \ {\displaystyle a_m = 0}, \\
\overline{\eta}^1 := \big(\sqrt{2}+1\big)(c_0\nu + b) + 2\sqrt{\big(\sqrt{2}+1\big)c_0\nu b} & \quad \text{if} \ \ {\displaystyle a_m = 1}. \end{cases}$$
By direct calculations, we obtain $\overline{\eta}^0 / \overline{\eta}^1 \in (0.356, 0.586)$, where the lower bound is attained when $c_0\nu = b$ and the upper bound is approached when $c_0\nu/b \rightarrow 0$ or $c_0\nu/b \rightarrow \infty$. This finding motivates us to wonder whether the recovery error can be reduced by around half in practice. This inference is further validated by numerical experiments in Section \ref{section6}.

\section{Rank consistency}\label{section4}

In this section we consider the asymptotic behavior of the estimator generated from the rank-correction step (\ref{eqnrcs}) in term of its rank. We expect that the resulting $\widehat{X}_m$ has the same rank as the true matrix $\overline{X}$. Theorem \ref{asmpnoi} only reveals a flavored parameter $\rho_m$ in terms of the optimal order but rather its exact value. In practice, for a chosen parameter $\rho_m$, there is hardly any clue to know the recovery performance of the resulting solution since the true matrix is unknown. However, if the rank property holds as expected, the observable rank information may be used to infer the recovery quality of the resulting solution of a parameter and thus help in parameter searching. Numerical experiments in Section \ref{section6} demonstrate the practicability of this idea.

For the purpose above, we study the rank consistency in the sense of Bach \cite{Bac08} under the setting that the matrix size is fixed. An estimator $X_m$ of the true matrix $\overline{X}$ is said to be rank consistent if $$\lim\limits_{m\rightarrow \infty}{\rm Pr}\big({\rm rank}(X_m)= {\rm rank}(\overline{X})\big)=1.$$
Throughout this section, we make the following assumptions:
\begin{assumption}\label{asmpfun}
The spectral operator $F$ is continuous at $\overline{X}$.
\end{assumption}
\begin{assumption}\label{asmpini}
The initial estimator $\widetilde{X}_m$ satisfies $\widetilde{X}_m \stackrel{p}{\rightarrow} \overline{X}$ as $m\rightarrow \infty$.
\end{assumption}

Epi-convergence in distribution gives us an elegant way in analyzing the asymptotic behavior of optimal solutions of a sequence of constrained optimization problems. Based on this technique, we obtain the following result.
\begin{theorem}\label{thmcons}
If $\rho_m \rightarrow 0$, then $\widehat{X}_m \stackrel{p}{\rightarrow} \overline{X}$ as $m\rightarrow \infty$.
\end{theorem}

We first focus on the characterization of necessary and sufficient conditions for rank consistency of $\widehat{X}_m$.
Unlike in the analysis of recovery error bound, additional information represented by the set $\mathcal{C}$  could affect the path along which $\widehat{X}_m$ converges to  $\overline{X}$ and thus may break the rank consistency. In the sequel, we only discuss two most common cases: the rectangular case $\mathcal{C} = \mathbb{V}^{n_1\times n_2}$ (recovering a rectangular matrix
or a symmetric/Hermitian matrix) and the positive semidefinite case $\mathcal{C} = \mathbb{S}_+^n$ (recovering a
symmetric/Hermitian positive semidefinite matrix).

For notational simplicity, we divide the index set $\beta$ into three subsets as
\begin{equation}\label{defbetadivision}
\beta^+ := \{k \in \beta \mid \langle \Theta_k, \overline{X} \rangle = b\}, \ \ \beta^- := \{k \in \beta \mid \langle \Theta_k, \overline{X} \rangle = -b\},  \ \ \beta^\circ := \beta \backslash (\beta^+ \cup \beta^-).
\end{equation}
Then, we define a linear operator $\mathcal{Q}_\beta^\dag : \mathbb{V}^{n_1\times n_2} \rightarrow \mathbb{V}^{n_1\times n_2}$ as
$$\mathcal{Q}_\beta^\dag(X) := \sum_{k\in\beta^\circ} \frac{1}{p_k} \langle \Theta_k, X\rangle \Theta_k + \sum_{k \in \beta^+} \frac{1}{p_k} (\langle \Theta_k, X \rangle)_- \Theta_k + \sum_{k \in \beta^-} \frac{1}{p_k} (\langle \Theta_k, X \rangle)_+ \Theta_k.$$
Here, we use the superscript ``$\dag$'' because of its inverse-like property in terms of
\begin{equation*}
\mathcal{Q}_\beta(\mathcal{Q}_\beta^\dag(Z)) = \mathcal{Q}_\beta^\dag(\mathcal{Q}_\beta(Z)) = \mathcal{P}_\beta(Z) \quad \forall \, Z \in \{Z \in \!\mathbb{V}^{n_1\times n_2} \mid  \mathcal{R}_{\beta^+}(Z) \leq 0, \mathcal{R}_{\beta^-}(Z) \geq 0\}.
\end{equation*}

By extending the arguments of Bach \cite{Bac08} for the nuclear norm penalized least squares estimator from the unconstrained case to the constrained case, we obtain the following results.

\begin{theorem}\label{thmnessuf}
For the rectangular case $\mathcal{C} = \mathbb{V}^{n_1\times n_2}$, consider the linear system
\begin{equation}\label{eqndelcond}
\overline{U}_2^\mathbb{T}\mathcal{Q}_\beta^\dag(\overline{U}_2 \Gamma\overline{V}_2^\mathbb{T})\overline{V}_2 = \overline{U}_2^\mathbb{T}\mathcal{Q}_\beta^\dag\big(\overline{U}_1 \overline{V}_1^\mathbb{T}-F(\overline{X})\big)\overline{V}_2.
\end{equation}
If $\rho_m \rightarrow 0$ and $\sqrt{m}\rho_m\rightarrow \infty$, then for the rank consistency of $\widehat{X}_m$,
\begin{description}
\setlength{\itemsep}{0pt}
\item[(i)] a necessary condition: (\ref{eqndelcond}) has a solution $\widehat{\Gamma} \in \!\mathbb{V}^{(n_1-r)\times (n_2-r)}$ with $\|\widehat{\Gamma}\|\le 1$;
\item[(ii)] a sufficient condition: (\ref{eqndelcond}) has a unique solution $\widehat{\Gamma} \in \mathbb{V}^{(n_1-r)\times (n_2-r)}$ with $\|\widehat{\Gamma}\| < 1$.
\end{description}
\end{theorem}

For the positive semidefinite case, the nuclear norm $\|X\|_*$ in (\ref{eqnrcs}) simply reduces to the trace $\langle I_n, X\rangle$. We assume that the Slater condition holds.

\begin{assumption}\label{asmpslater}
For the positive semidefinite case $\mathcal{C} = \mathbb{S}_+^n$, the Slater condition holds, i.e., there exists some $X^0 \in \mathbb{S}^n_{++}$ such that $\mathcal{R}_\alpha(X^0) = \mathcal{R}_\alpha(\overline{X})$ and $\|\mathcal{R}_\beta(\!X^0)\|_\infty < b$.
\end{assumption}
\begin{theorem}\label{thmnessufpos}
For the positive semidefinite case $\mathcal{C} = \mathbb{S}_+^n$, consider the linear system
\begin{equation}\label{eqndelcondpos}
 \overline{P}_2^\mathbb{T}\mathcal{Q}_\beta^\dag(\overline{P}_2 \Lambda \overline{P}_2^\mathbb{T})\overline{P}_2 = \overline{P}_2^\mathbb{T}\mathcal{Q}_\beta^\dag\big(I_n- F(\overline{X})\big)\overline{P}_2.
\end{equation}
Under Assumption \ref{asmpslater}, if $\rho_m \rightarrow 0$ and $\sqrt{m}\rho_m\rightarrow \infty$, then for the rank consistency of $\widehat{X}_m$,
\begin{description}
\setlength{\itemsep}{0pt}
\item[(i)] a necessary condition:  (\ref{eqndelcondpos}) has a solution $\widehat{\Lambda} \in \mathbb{S}^{n-r}_+$;
\item[(ii)] a sufficient condition:  (\ref{eqndelcondpos}) has a unique solution $\widehat{\Lambda} \in \mathbb{S}^{n-r}_{++}$.
\end{description}
\end{theorem}

Next, we provide a theoretical guarantee on the uniqueness of the solution to the linear systems (\ref{eqndelcond})
and (\ref{eqndelcondpos}) with the help of constraint nondegeneracy. The concept of constraint nondegeneracy was pioneered by Robinson \cite{Rob84} and later extensively developed
by Bonnans and Shapiro \cite{BonS00}.
We say that the constraint nondegeneracy holds at $\overline{X}$ to (\ref{eqnrcs}) with $\mathcal{C} = \mathbb{V}^{n_1\times n_2}$ if
\begin{equation}\label{eqncndc}
\mathcal{R}_{\alpha \cup \beta^+ \cup \beta^-} \big(\mathcal{T}(\overline{X})\big) = \mathbb{R}^{|\alpha \cup \beta^+ \cup \beta^-|},
\end{equation}
where $\mathcal{T}(\overline{X}) = \big\{H \in \mathbb{V}^{n_1\times n_2} \mid  \overline{U}_2^\mathbb{T} H \overline{V}_2 = 0 \big\}$. Meanwhile, we say that the constraint nondegeneracy holds at $\overline{X}$ to (\ref{eqnrcs}) with $\mathcal{C} = \mathbb{S}_+^n$ if
\begin{equation}\label{eqncndcpos}
\mathcal{R}_{\alpha \cup \beta^+ \cup \beta^-} \big(\text{lin}(\mathcal{T}_{\mathbb{S}^n_+}(\overline{X}))\big) = \mathbb{R}^{|\alpha \cup \beta^+ \cup \beta^-|},
\end{equation}
where $\text{lin}(\mathcal{T}_{\mathbb{S}^n_+}(\overline{X})) = \big\{H \in \mathbb{S}^{n} \mid  \overline{P}_2^\mathbb{T} H \overline{P}_2 = 0 \big\}$. One may refer to \ref{AppendixConNon} for more details of constraint nondegeneracy.

To take a closer look at the linear systems (\ref{eqndelcond})
and (\ref{eqndelcondpos}), we define  linear operators $\mathcal{B}_1:\mathbb{V}^{n_1\times n_2} \to \mathbb{V}^{(n_1-r)\times (n_2-r)}$
and $\mathcal{B}_2:\mathbb{V}^{(n_1-r)\times(n_2-r)} \to \mathbb{V}^{(n_1-r)\times(n_2-r)}$ associated with $\overline{X}$, respectively, by
\begin{equation}\label{linear-operator}
\mathcal{B}_1(Y): = \overline{U}_2^\mathbb{T}\mathcal{Q}_{\beta}^{\dag}(Y)\overline{V}_2 \ \ {\rm and}\ \
\mathcal{B}_2(Z): = \overline{U}_2^\mathbb{T}\mathcal{Q}_{\beta}^{\dag}(\overline{U}_2 Z \overline{V}_2^T)\overline{V}_2,
\end{equation}
where $Y \in \mathbb{V}^{n_1\times n_2}$ and $Z \in \mathbb{V}^{(n_1-r)\times (n_2-r)}$. From
the definition of $\mathcal{Q}_{\beta}^{\dag}$, we know that the operator $\mathcal{B}_2$ is self-adjoint and positive semidefinite. %Moreover, in the light of the definition of the spectral operator $F$,  define $\widehat{g}(\overline{X})$ to be a vector in
%$\mathbb{R}^r$ as
%\begin{equation}\label{eqnvecgr}
%\widehat{g}(\overline{X}) := \big(1-f_1(\sigma(\overline{X})),\ldots,1-f_r(\sigma(\overline{X}))\big)^{\mathbb{T}}.
%\end{equation}
Then, for the rectangular case $\mathcal{C} = \mathbb{V}^{n_1\times n_2}$, the linear system (\ref{eqndelcond}) can be rewritten as
\begin{equation}\label{eqndelcond1eqiv}
 \mathcal{B}_2(\Gamma) = \mathcal{B}_1(\overline{U}_1\overline{V}_1^{\mathbb{T}}-F(\overline{X})), \quad \Gamma \in \mathbb{V}^{(n_1-r)\times (n_1-r)},
\end{equation}
and for the positive semidefinite case $\mathcal{C} = \mathbb{S}^n_+$, the linear system (\ref{eqndelcondpos}) can be rewritten as
\begin{equation}\label{eqndelcond2eqiv}
 \mathcal{B}_2(\Lambda) =\mathcal{B}_2(I_{n-r})+\mathcal{B}_1(\overline{P}_1\overline{P}_1^{\mathbb{T}} - F(\overline{X})), \quad \Lambda \in \mathbb{S}^{n-r},
\end{equation}
since both $\overline{U}_i$ and $\overline{V}_i$ reduce to $\overline{P}_i$ for $i=1,2$ for $\overline{X} \in \mathbb{S}_+^n.$

Clearly, the invertibility of $\mathcal{B}_2$ is equivalent to the uniqueness of the solution to the linear systems (\ref{eqndelcond}) and (\ref{eqndelcondpos}). The following result provides a link between the constraint nondegeneracy and the positive definiteness of $\mathcal{B}_2$.

\begin{theorem}\label{thmsoluni}
For either the rectangular case $\mathcal{C} = \mathbb{V}^{n_1\times n_2}$ or the positive semidefinite case $\mathcal{C} = \mathbb{S}_+^n$, if the constraint nondegeneracy holds at $\overline{X}$ to the problem (\ref{eqnrcs}), then the self-adjoint linear operator $\mathcal{B}_2$ defined by (\ref{linear-operator}) is positive definite.
\end{theorem}

Combining Theorems \ref{thmnessuf}, \ref{thmnessufpos} and \ref{thmsoluni} together with (\ref{eqndelcond1eqiv}) and (\ref{eqndelcond2eqiv}), we immediately have the following result of rank consistency.

\begin{theorem}\label{thmgenconsis}
Suppose that $\rho_m \rightarrow 0$ and $\sqrt{m}\rho_m\rightarrow \infty$. If
\begin{description}
\item[(i)] for the rectangular case $\mathcal{C} = \mathbb{V}^{n_1\times n_2}$, the constraint nondegeneracy (\ref{eqncndc}) holds at $\overline{X}$ to the problem (\ref{eqnrcs}) and
\begin{equation}\label{eqnsufso}
 \big\|\mathcal{B}_2^{-1}\mathcal{B}_1 (\overline{U}_1\overline{V}_1^{\mathbb{T}}-F(\overline{X}))\big)\big\| <1;
\end{equation}
\item[(ii)] for the positive semidefinite case $\mathcal{C} = \mathbb{S}_+^n$, the constraint nondegeneracy (\ref{eqncndcpos}) holds at $\overline{X}$ to the problem (\ref{eqnrcs}) and
\begin{equation}\label{eqnsufsopos}
I_{n-r} + \mathcal{B}_2^{-1}\mathcal{B}_1 (\overline{P}_1\overline{P}_1^{\mathbb{T}}-F(\overline{X})) \in \mathbb{S}_{++}^{n-r},
\end{equation}
\end{description}
then the estimator $\widehat{X}_m$ generated from the rank-correction step (\ref{eqnrcs}) is rank consistent.
\end{theorem}

From Theorem \ref{thmgenconsis}, it is not difficult to see that when $F(\overline{X})$ is sufficiently close to $\overline{U}_1\overline{V}_1^\mathbb{T}$, the conditions (\ref{eqnsufso}) and (\ref{eqnsufsopos}) hold automatically and so does the rank consistency. Thus, Theorem \ref{thmgenconsis} provides us a guideline to construct a suitable rank-correction function $F$ to achieve the rank consistency.
In particular, for the positive semidefinite matrix completion, we further consider two important classes as follows.

\begin{description}
\item[Class I:] The covariance matrix completion with partial positive diagonal entries fixed.

Due to the positive semidefinite structure, the magnitudes of off-diagonal entries are fully controlled by the magnitudes of diagonal entries. Therefore, we remove all the bounded constraints corresponding to off-diagonal entries from the rank-correction step (\ref{eqnrcs}) as they are redundant. Thus, the constraints are reduced to
$$X_{ii} = \overline{X}_{ii} \ \ \forall\, i \in \pi, \quad X_{ii} \leq b \ \ \forall\, i \in \pi^c, \quad X \in \mathbb{S}_+^n,$$ where $(\pi, \pi^c)$ is a partition of the index set $\{1,\ldots, n\}$.
This class of problems includes the correlation matrix completion as a special case, in which  all diagonal entries are fixed to be ones.
\item[Class II:] The density matrix completion with its trace fixed to be one.

Due to the positive semidefinite structure, all the coefficients of Pauli basis are controlled because of the trace one constraint. Therefore, we remove all the bounded constraints from the rank-correction step (\ref{eqnrcs}) as they are redundant. Thus, in this case the constraints are reduced to
$$\frac{1}{\sqrt{n}}\text{Tr}(X) = \frac{1}{\sqrt{n}}, \quad X \in \mathbb{S}_+^n.$$
\end{description}

Interestingly, for the matrix completion problems of Classes I and II, the constraint nondegeneracy automatically holds at $\overline{X}$. More importantly, if observations are sampled uniformly at random, the rank consistency can be guaranteed for a broad class of rank-correction functions $F$.

\begin{theorem}\label{thmrccordencons}
For the matrix completion problems of Classes I and II under uniform sampling, if $\rho_m \rightarrow 0, \ \sqrt{m}\rho_m\rightarrow \infty$ and $F$ is a spectral operator associated with a symmetric function $f:\mathbb{R}^n \rightarrow \mathbb{R}^n$ such that for $i=1,\ldots,n$,
\begin{equation}\label{eqnfcorden}
\left\{\begin{array}{ll}  f_i(x) > 0 \quad & \text{if} \ x_i >0, \\  f_i(x)  =0 \quad &  \text{if} \ x_i =0,\end{array}\right. \qquad \forall\,  x\in \mathbb{R}_+^n \ \ \text{and}  \ \ \forall\, i=1,\ldots,n,
\end{equation}
then the estimator $\widehat{X}_m$ generated from the rank-correction step (\ref{eqnrcs}) is rank consistent.
\end{theorem}

\section{Construction of the rank-correction function}\label{section5}

 In this section, we focus on the construction of  a suitable rank-correction function $F$ based on the results in Sections \ref{section3} and \ref{section4}. For achieving a smaller recovery error, according to Theorem \ref{thmstobd}, we desire a construction such that $F(\widetilde{X}_m)$ is close to $\overline{U}_1\overline{V}_1^\mathbb{T}$. Meanwhile, for achieving the rank consistency, according to
 Theorem  \ref{thmgenconsis}, we desire a construction such that $F(\overline{X})$ is close to $\overline{U}_1\overline{V}_1^\mathbb{T}$.  Therefore, these two guidelines consistently suggest a natural idea, i.e., if possible, choosing $$F(X) \approx U_1V_1^{\mathbb{T}} \quad \text{near} \ \overline{X}.$$
 Next, we proceed with the construction of the rank-correction function $F$ for the rectangular case. For the positive semidefinite case, one only needs to replace  the singular value decomposition with the eigenvalue decomposition and conduct exactly the same analysis.

\subsection{The  rank is known}\label{subsecrankknown}

If the rank of the true matrix $\overline{X}$ is known, it is clear that the best choice of $F$ is
\begin{equation}\label{eqnchofun1}
F(X) :=  U_1 V_1^\mathbb{T},
\end{equation}
where $(U,V)\in \mathbb{O}^{n_1,n_2}(X)$ and $X \in \mathbb{V}^{n_1\times n_2}$.  Note that $F$ defined by (\ref{eqnchofun1}) is not a spectral operator over the whole space of $\mathbb{V}^{n_1\times n_2}$, but in a neighborhood of $\overline{X}$ it is indeed a spectral operator and is actually twice continuously differentiable (see, e.g., \cite[Proposition 8]{DinST10}). With this rank-correction function,
the rank-correction step is essentially the same as a single  step of the majorized penalty method developed in \cite{GaoS10}.
%By Theorems \ref{thmnessuf} and \ref{thmsoluni}, we immediately obtain the following result.
%
%\begin{corollary}
%For either the rectangular case $\mathcal{C} = \mathbb{V}^{n_1\times n_2}$ or the positive semidefinite case $\mathcal{C} = \mathbb{S}_+^n$, suppose that the rank of the true matrix $\overline{X}$ is known and the constraint nondegeneracy holds at $\overline{X}$ to the problem (\ref{eqnrcs}).
%If $\rho_m\rightarrow 0$, $\sqrt{m}\rho_m\rightarrow \infty$ and  $F$ is chosen by (\ref{eqnchofun1}), then the estimator $\widehat{X}_m$ generated from the rank-correction step (\ref{eqnrcs}) is rank consistent.
%\end{corollary}

\subsection{The rank is unknown}

If the rank of the true matrix $\overline{X}$ is unknown, we intend to construct a spectral operator $F$ to imitate the case when the rank is known. Here, we propose $F$ to be a spectral operator
 \begin{equation}\label{eqnchofun}
  F(X) := U \text{Diag}\big(f(\sigma(X))\big)V^\mathbb{T}
 \end{equation}
 associated with the symmetric function $f:\mathbb{R}^n \rightarrow \mathbb{R}^n$ defined by
 \begin{equation}\label{eqnchofun2}
  f_i(x) = \begin{cases} {\displaystyle  \phi\left(\frac{x_i}{\|x\|_\infty}\right)} \quad  & \text{if}\ x \in \mathbb{R}^n \backslash \{0\},\\
   0  \quad &  \text{if}\ x =0,\end{cases}
 \end{equation}
 where $(U,V)\in \mathbb{O}^{n_1,n_2}(X)$, $X \in \mathbb{V}^{n_1\times n_2}$, and the scalar function $\phi:\mathbb{R} \rightarrow \mathbb{R}$ takes the form
 \begin{equation}\label{eqnchofun3}
  \phi(t): = \text{sgn}(t) (1+\varepsilon^\tau)\frac{|t|^\tau}{|t|^\tau+\varepsilon^\tau},  \quad   t \in \mathbb{R},
 \end{equation}
 for some $\tau >0$ and $\varepsilon >0$.
 %By noting that for each $t$, $\phi(t) \rightarrow \text{sgn}(t)$ as $\varepsilon %\downarrow 0$, we directly obtain the following result.

\begin{corollary}\label{cororankcons}
Let $F$ be a spectral operator defined by (\ref{eqnchofun}), (\ref{eqnchofun2}) and (\ref{eqnchofun3}).
\
\begin{description}
\item[(i)] If $\frac{\|\widetilde{X}_m - \overline{X}\|_F}{\sigma_r(\overline{X})} <  \frac{1}{\sqrt{2}}\big(1-e^{-\sqrt{2r}}\big)$,
    %$ \|\widetilde{X}_m - \overline{X}\|_F / \sigma_r(\overline{X}) <  0.535$,
    then for any $\varepsilon $ satisfying
    $ \frac{\sigma_{r+1}(\widetilde{X}_m)}{\sigma_1(\widetilde{X}_m)} < \varepsilon < \frac{\sigma_r(\widetilde{X}_m)}{\sigma_1(\widetilde{X}_m)}$,
    %$\sigma_{r+1}(\widetilde{X}_m)/\sigma_1(\widetilde{X}_m) < \varepsilon < %\sigma_r(\widetilde{X}_m)/\sigma_1(\widetilde{X}_m)$,
    there exists some  $\overline{\tau}_1 >0$ such that $a_m  < 1$ for any $F$ with $\tau \geq \overline{\tau}_1$.
\item[(ii)] Suppose that the constraint nondegeneracy holds at $\overline{X}$ to the problem (\ref{eqnrcs}). If $\rho_m \rightarrow 0$ and $\sqrt{m}\rho_m\rightarrow \infty$, then for any $\varepsilon$ satisfying
    $0 < \varepsilon < \frac{\sigma_r(\overline{X})}{\sigma_1(\overline{X})}$,
    % $\varepsilon < \sigma_r(\overline{X})/\sigma_1(\overline{X})$,
    there exists some $\overline{\tau}_2 > 0$ such that the rank consistency of $\widehat{X}_m$ holds for any $F$ with $\tau \geq \overline{\tau}_2$.
\end{description}
\end{corollary}

The proof of Corollary \ref{cororankcons} is straightforward so we omit it. Corollary \ref{cororankcons} suggests an ideal choice of $\varepsilon$ for the recovery error reduction,  i.e.,  $\varepsilon \in \left(\frac{\sigma_{r+1}(\widetilde{X}_m)}{\sigma_1(\widetilde{X}_m)}, \frac{\sigma_{r}(\widetilde{X}_m)}{\sigma_1(\widetilde{X}_m)} \right)$, provided that $\widetilde{X}_m$ does not deviate too much from $\overline{X}_m$, and also an ideal choice of $\varepsilon$ for rank consistency, i.e., $\varepsilon \in \left(0, \frac{\sigma_{r}(\overline{X}_m)}{\sigma_1(\overline{X}_m)} \right)$. Note that these two intervals may not overlap each other, implying the theoretical possibility that the recovery error reduction and the rank consistency may not be achieved simultaneously if the initial estimator $\widetilde{X}_m$ is not close to $\overline{X}_m$.

\begin{figure}
  \begin{center}
  \subfigure[$\varepsilon=0.1$ with different $\tau>0$]{\includegraphics[height=5cm]{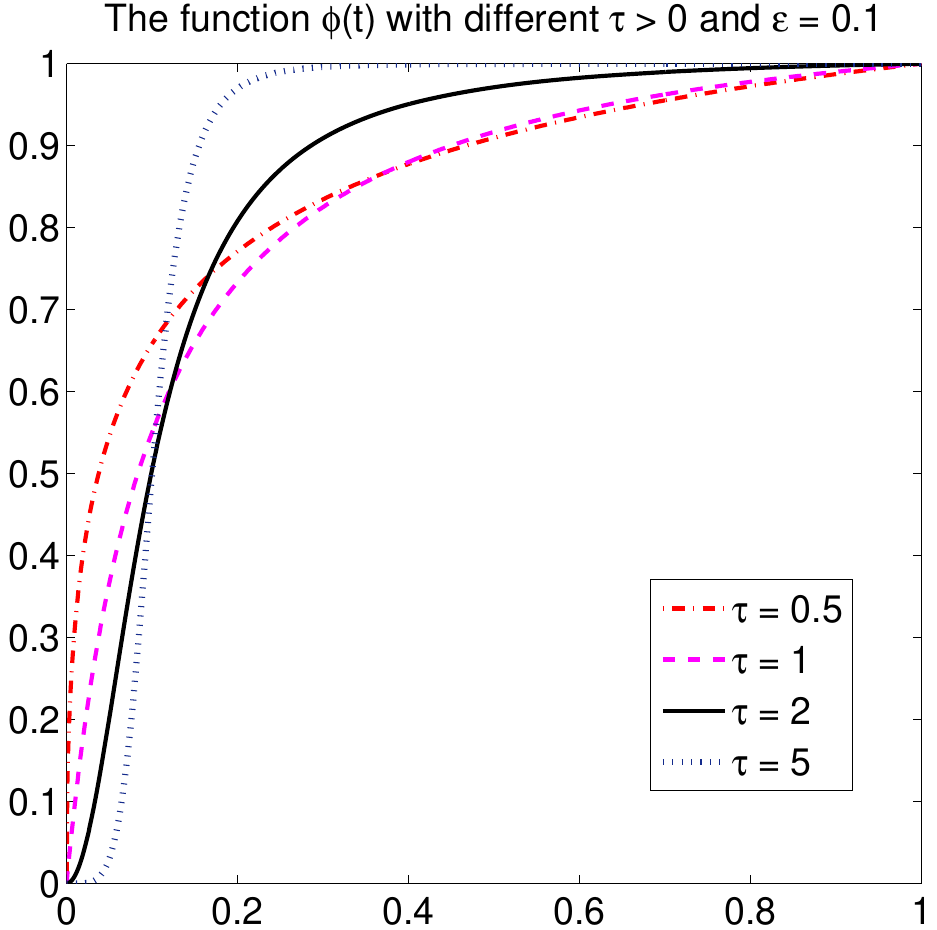}} \hspace{2cm}
  \subfigure[$\tau=2$ with different $\varepsilon>0$]{\includegraphics[height=5cm]{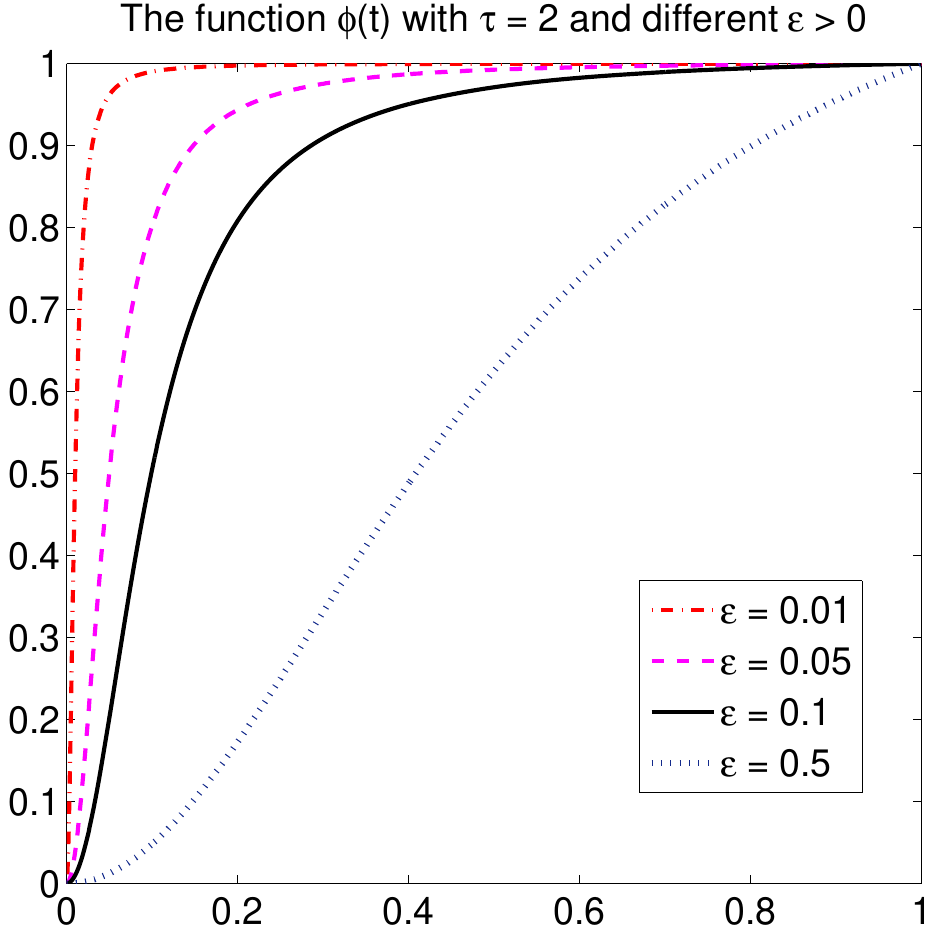}}
  \caption{Shapes of the function $\phi$ with different $\varepsilon>0$ and $\tau>0$}\label{dad}
  \end{center}
\end{figure}

The interval of $\varepsilon$ for the recovery error reduction is disclosed if the true rank is accessible. Therefore, this ideal interval is an important insight that can be used to guide the choice of $\varepsilon$ in practice since the initial $\widetilde{X}_m$ should contain some information of the true rank in general.  Indeed, the value of $\varepsilon$ can be regarded as a divide of confidence on whether $\sigma_i(\widetilde{X}_m)$ is believed to come from a nonzero singular values of $\overline{X}$ with perturbation --- positive confidence if $\sigma_i(\widetilde{X}_m) > \varepsilon \sigma_1(\widetilde{X}_m)$ and negative confidence if $\sigma_i(\widetilde{X}_m) < \varepsilon \sigma_1(\widetilde{X}_m)$. Next we look for a suitable $\tau$. It is observed from Figure \ref{dad} that the parameter $\tau>0$ mainly controls the shape of $\phi$ over $t\in [0,1]$. The function $\phi$ is concave if $0<\tau\leq 1$ and $S$-shaped with a single inflection point at $\varepsilon\big(\frac{\tau-1}{\tau+1}\big)^{1/\tau}$ if $\tau>1$. It should be good to choose an $S$-shaped function $\phi$. But one also needs to take account of the steepness of $\phi$, which increases when $\tau$ increases. In particular for any $\varepsilon$ satisfying $0<\varepsilon<1$, $\phi$ approaches to the step function taking the value $0$ if $0\leq t<\varepsilon$ and the value $1$ if $ \varepsilon < t \leq 1$ as $\tau \rightarrow \infty$. Since the rank of $\overline{X}$ is unknown and the singular values of $\widetilde{X}_m$ are unpredictable, choosing a large $\tau$ could be risky. Therefore, one needs to choose $\tau$ with certain conservation, sacrificing certain recovery quality in exchange for robustness strategically. Here, we provide a recommendation of the choices $\varepsilon \approx  0.05$ (or within $0.01\sim 0.1$) and $\tau = 2$ (or within $1 \sim 3$) for most cases, particularly when the initial estimator is generated from the nuclear norm penalized least squares problem. These choices have performed very  stably  for plenty of problems, as validated in Section \ref{section6}.

 We also remark that for the positive semidefinite case, the rank-correction function defined by (\ref{eqnchofun}), (\ref{eqnchofun2}) and
 (\ref{eqnchofun3}) is related to the reweighted trace norm for the matrix rank minimization proposed by Fazel et al. \cite{FazHB03, MohF10}. The reweighted trace norm in \cite{FazHB03, MohF10} for the positive semidefinite case is $\langle (X^k+\varepsilon I_n)^{-1}, X\rangle$, which arises from the derivative of the surrogate function $\log\det(X+\varepsilon I_n)$ of the rank   at an iterate $X^k$, where $\varepsilon$ is a small positive constant. Meanwhile, in our proposed rank-correction step, if we choose $\tau = 1$,  then $I_n-\frac{1}{1+\varepsilon}F(\widetilde{X}_m) = \varepsilon'(\widetilde{X}_m+\varepsilon' I_n)^{-1}$ with $\varepsilon' =\varepsilon \|\widetilde{X}_m\|$. Superficially, similarity occurs; however, it is notable that $\varepsilon'$ depends on $\widetilde{X}_m$, which is different from the constant $\varepsilon$ in \cite{FazHB03, MohF10}. More broadly speaking, the rank-correction function $F$ defined by (\ref{eqnchofun}), (\ref{eqnchofun2}) and (\ref{eqnchofun3}) is not a gradient of any real-valued function. This distinguishes our proposed rank-correction step from the reweighted trace norm minimization in \cite{FazHB03, MohF10} even for the positive semidefinite case.

\section{Numerical experiments}\label{section6}

 In this section, we validate the power of our proposed rank-correction step on the recovery
 by applying it to different matrix completion problems. We adopted the proximal alternating direction method
of multipliers (proximal ADMM) to solve the optimization problem (\ref{eqnrcs}). For more details of the proximal ADMM, the readers may
refer to Appendix B of \cite{FazPST12}.  For convenience, in the sequel, the {\rm NNPLS} estimator and the {\rm RCS} estimator, respectively, stand for
the estimators from the nuclear norm penalized least squares problem (i.e., $F \equiv 0$) and the rank-correction step (\ref{eqnrcs}) with $F$ specified in Section \ref{section5}.
Given an estimator $X_m$ of $\overline{X}_m$, the {\bf relative error} ({\bf relerr} for short) is defined by
$${\rm relerr}=\frac{\|X_m - \overline{X}\|_{F}}{\max(10^{-8},\|\overline{X}\|_{F})}.$$

%---------------------------------------------------------------------------------------------------Subsection6.1
 \subsection{Influence of fixed basis coefficients on the recovery}\label{subsec6.1}

 In this subsection, we test
 the performance of the {\rm NNPLS} estimator and the {\rm RCS} estimator
 for different patterns of fixed basis coefficients.
 We randomly generated a correlation matrix by the following command:
 \vspace{-0.3cm}
 \begin{verbatim}
   M = randn(n,r)/sqrt(sqrt(n));  ML = weight*M(:,1:k);  M(:,1:k) = ML;
     Xtemp = M*M'; D = diag(1./sqrt(diag(Xtemp))); X_bar = D*Xtemp*D.
 \end{verbatim}
 \vspace{-0.8cm}
 \noindent
 We took the true matrix $\overline{X}=$ {\ttfamily X\_bar} with dimension {\ttfamily n} $=500$, rank {\ttfamily r} $=5$, {\ttfamily weight} $=5$ and {\ttfamily k} $=1$. Here, the parameter {\ttfamily weight} is used to control the relative magnitude difference between the first $k$ largest eigenvalues and the left $r-k$ nonzero eigenvalues. We randomly fixed partial diagonal and off-diagonal entries of $\overline{X}$ and then uniformly sampled the rest entries
 with i.i.d. Gaussian noise. The noise level, defined by $\|\nu \xi\|_2/\|y\|_2$ in (\ref{eqnobs}) hereafter, was set to be $10\%$ and the upper bound of the non-fixed diagonal entries was set to be $1$. We further assumed that the rank of the true matrix was known so that for RCS estimator we chose the rank-correction function (\ref{eqnchofun1}).

\begin{figure}[htbp]
  \begin{center}
     \subfigure[Nuclear norm: diag=0, off-diag=0]{\includegraphics[width=7cm]{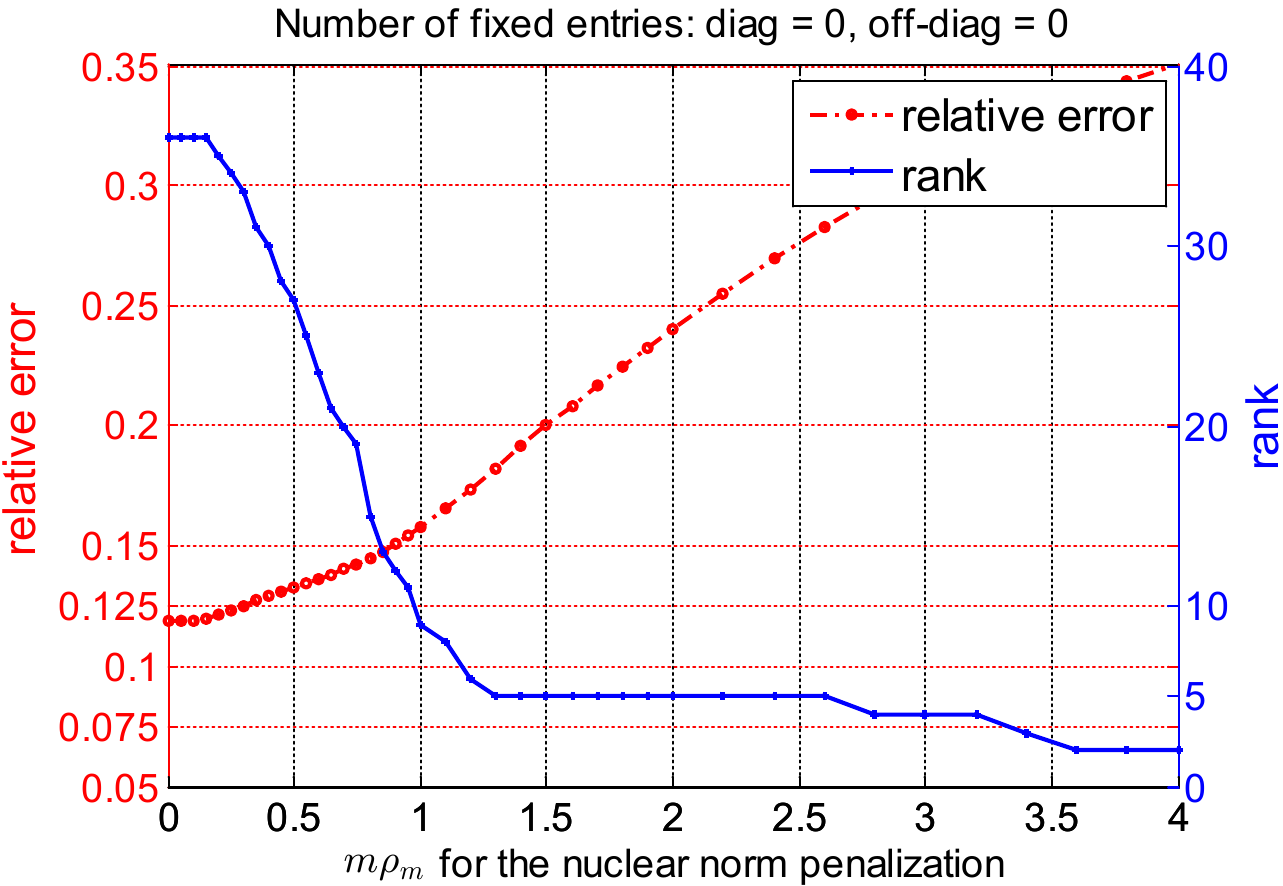}}
     \subfigure[Rank-correction step: diag=0, off-diag=0]{\includegraphics[width=7cm]{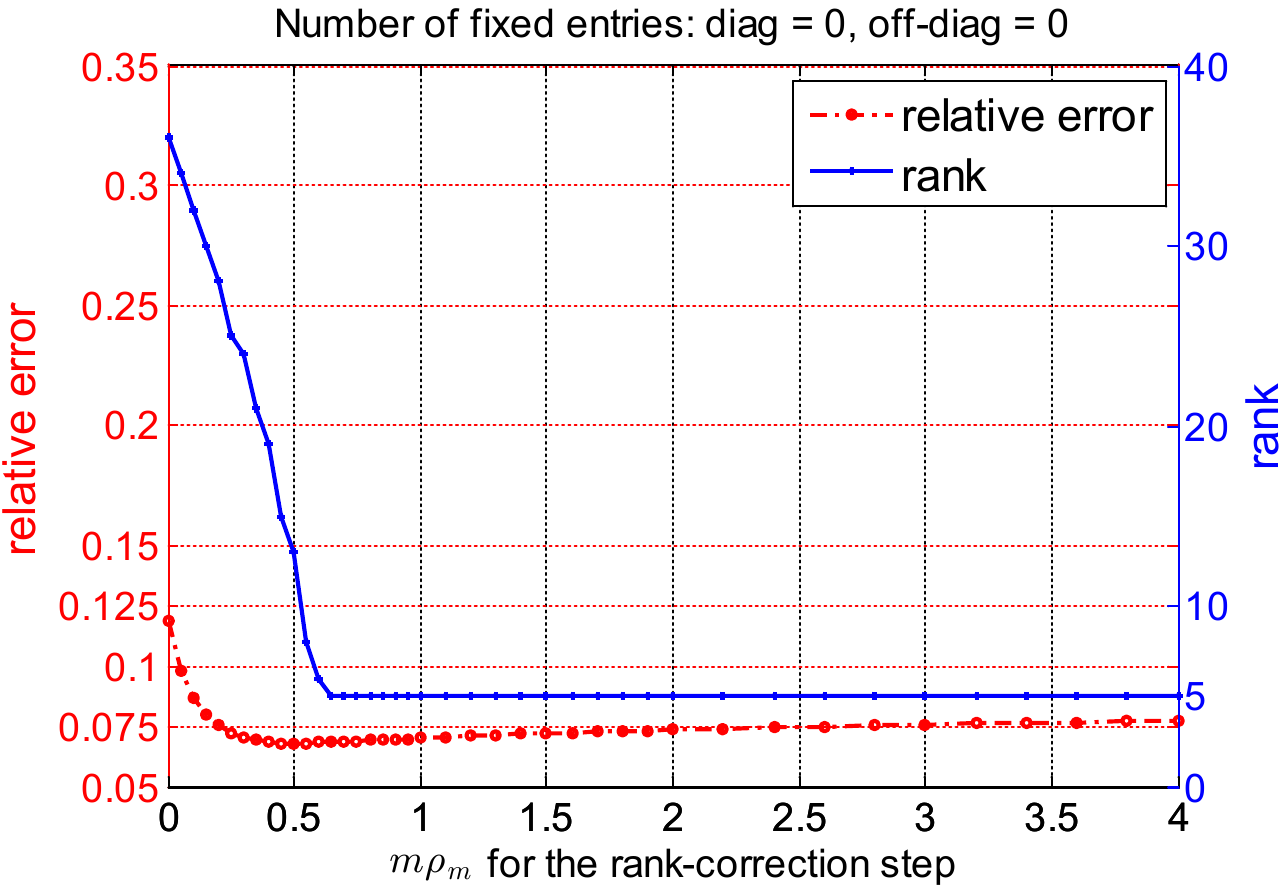}}\\
     \subfigure[Nuclear norm: diag=n/2, off-diag=0]{\includegraphics[width=7cm]{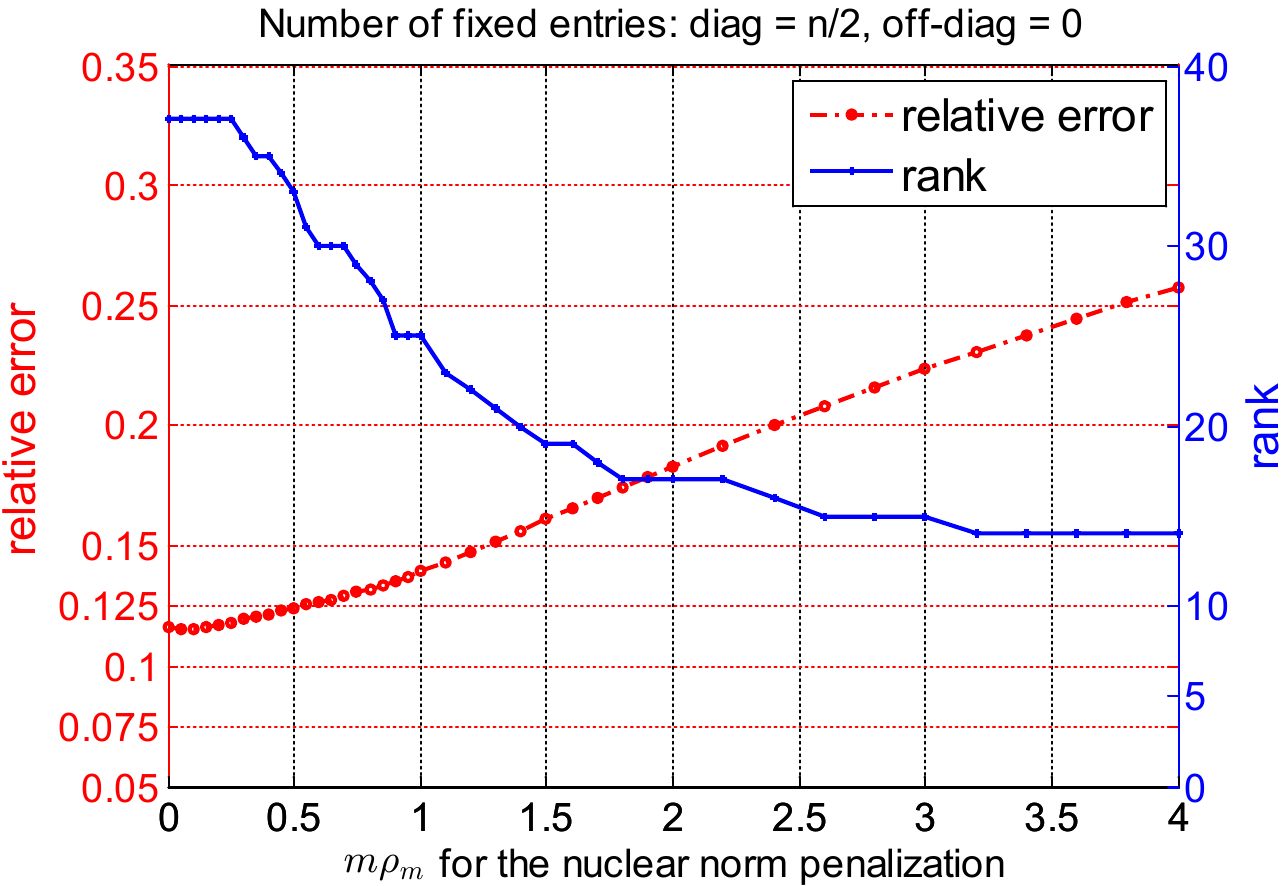}}
     \subfigure[Rank-correction step: diag=n/2, off-diag=0]{\includegraphics[width=7cm]{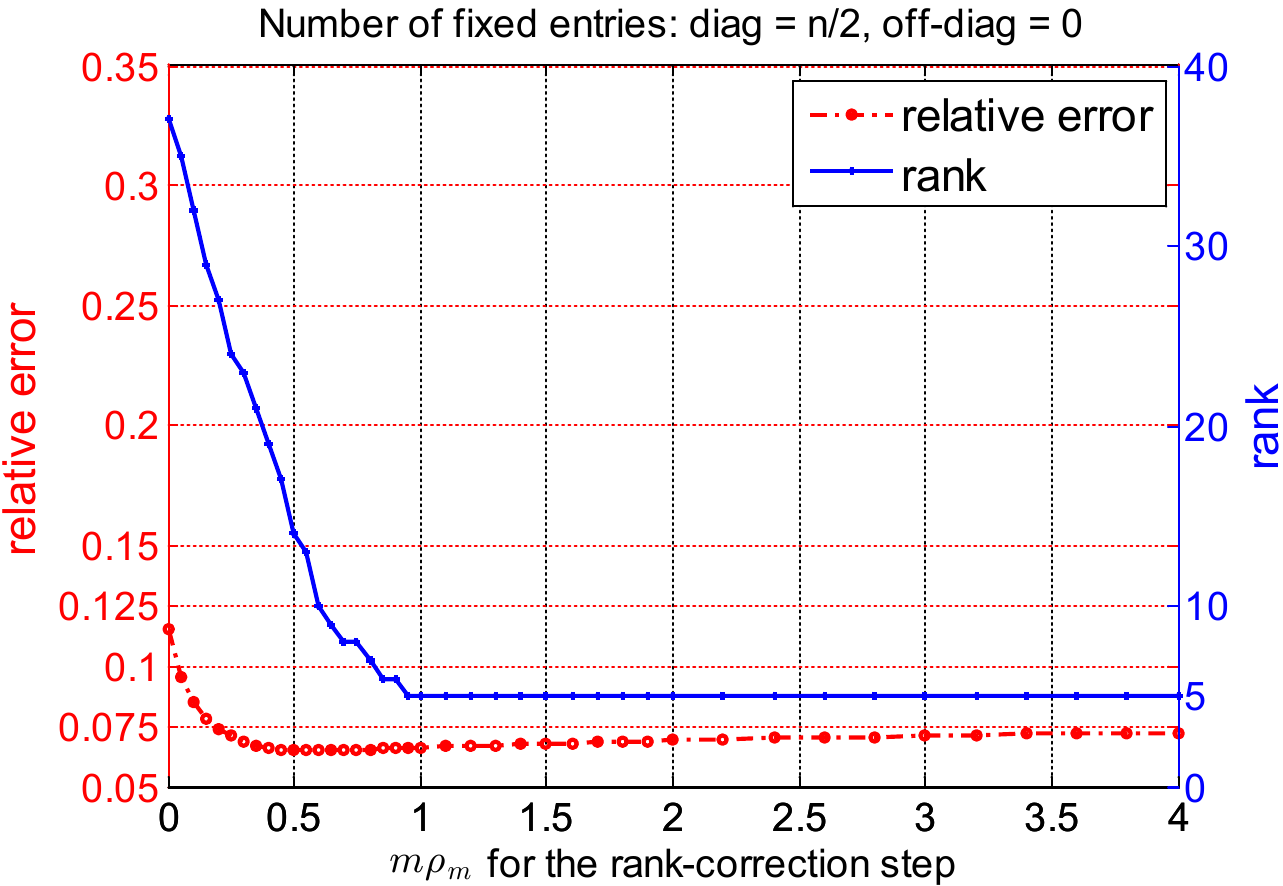}}\\
     \subfigure[Nuclear norm: diag=n, off-diag=0]{\includegraphics[width=7cm]{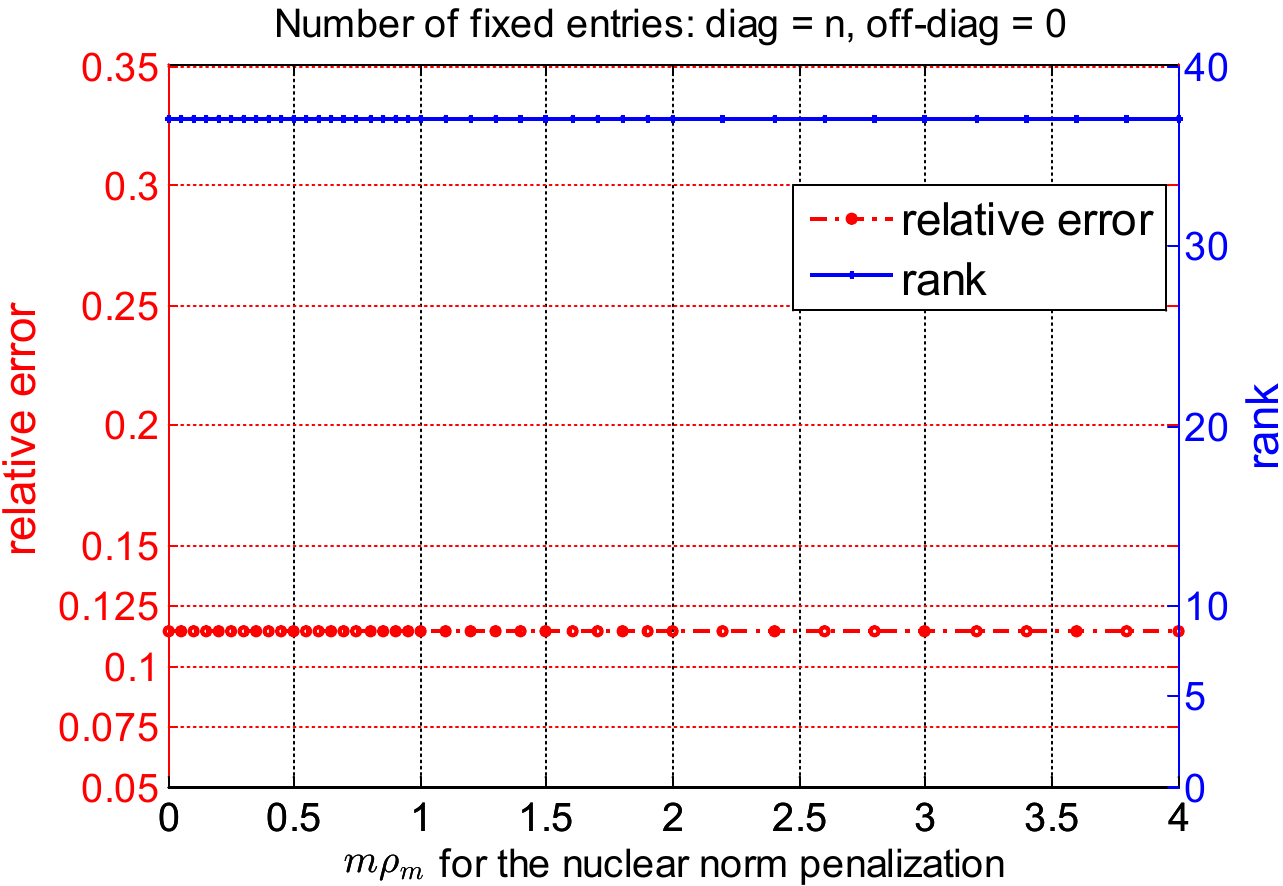}}
     \subfigure[Rank-correction step: diag=n, off-diag=0]{\includegraphics[width=7cm]{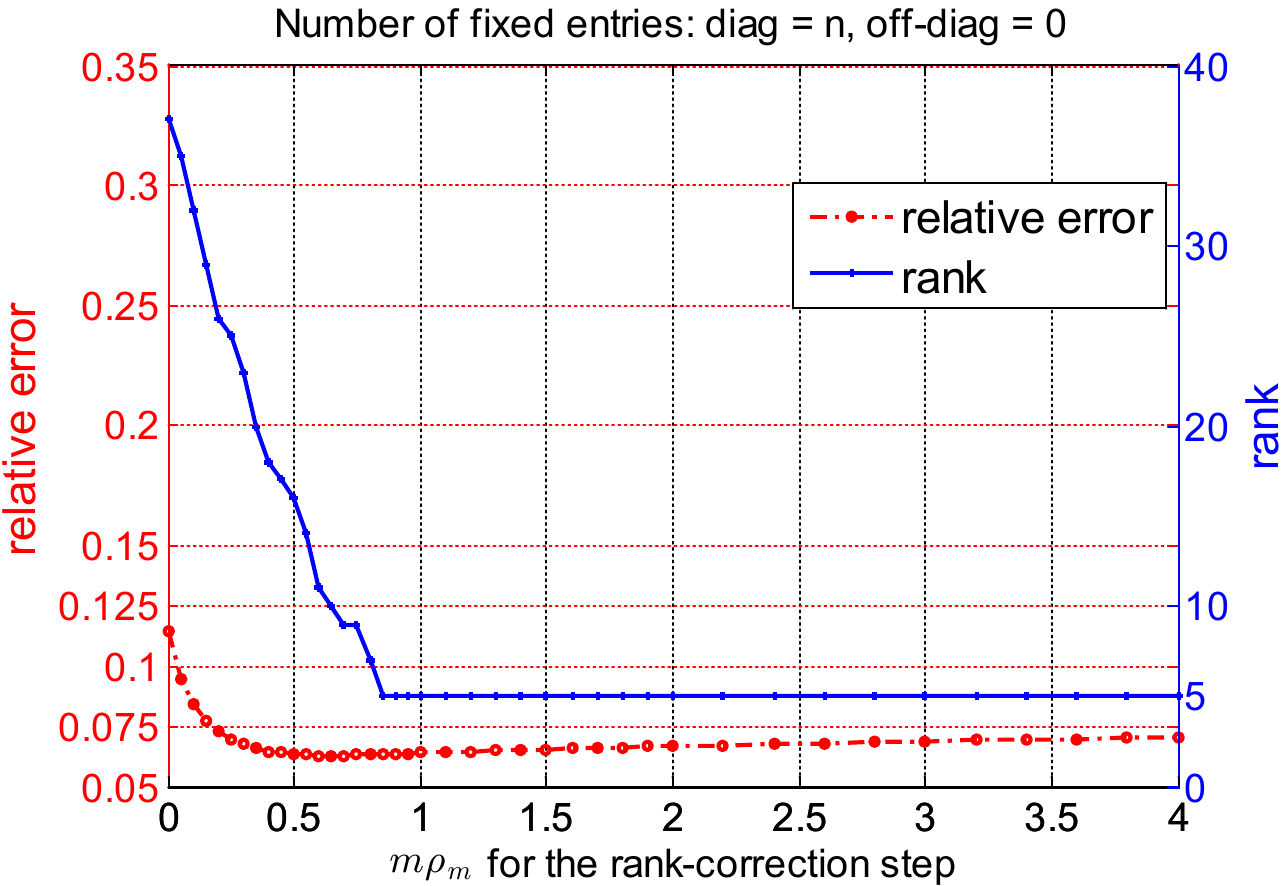}}\\
     \subfigure[Nuclear norm: diag=n, off-diag=n/2]{\includegraphics[width=7cm]{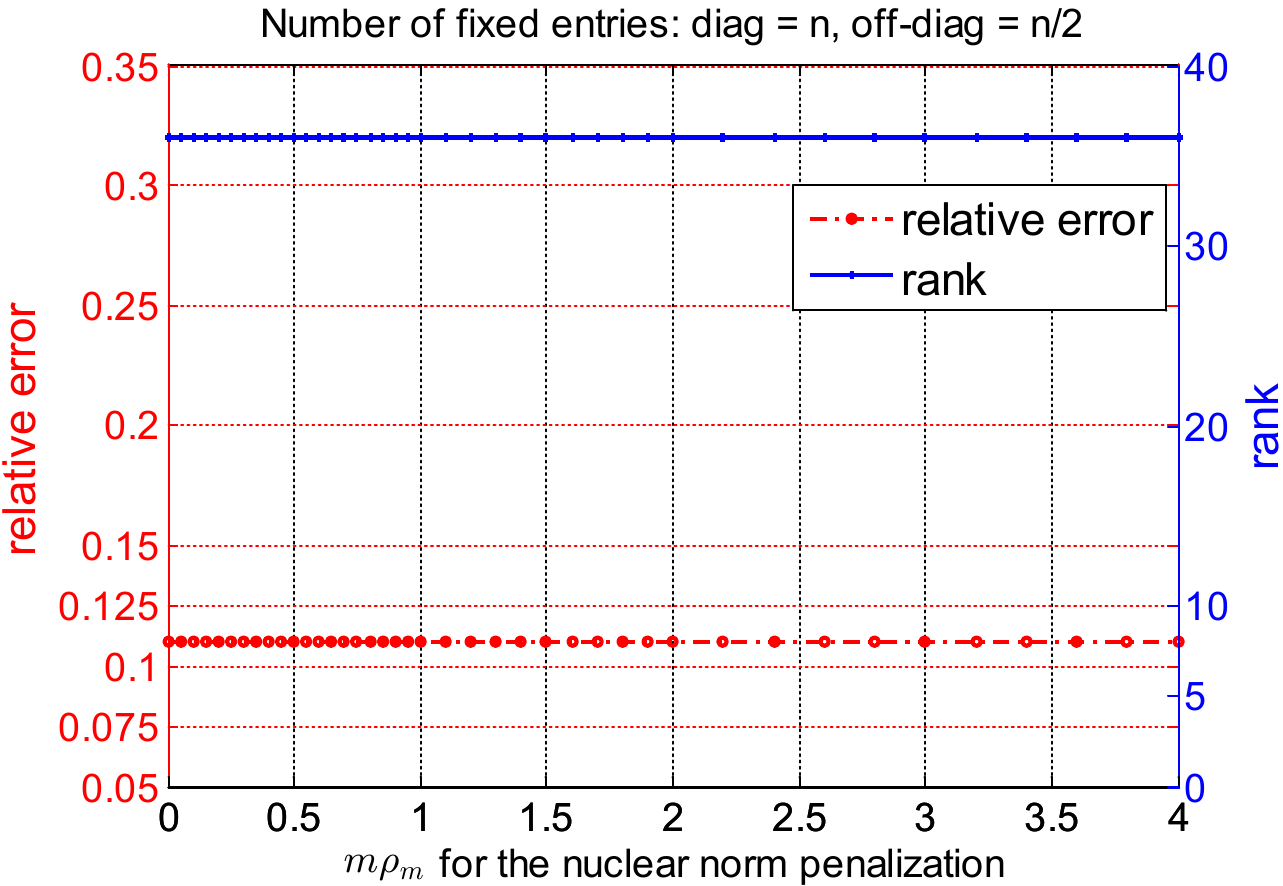}} \hspace{0cm}
     \subfigure[Rank-correction step: diag=n, off-diag=n/2]{\includegraphics[width=7cm]{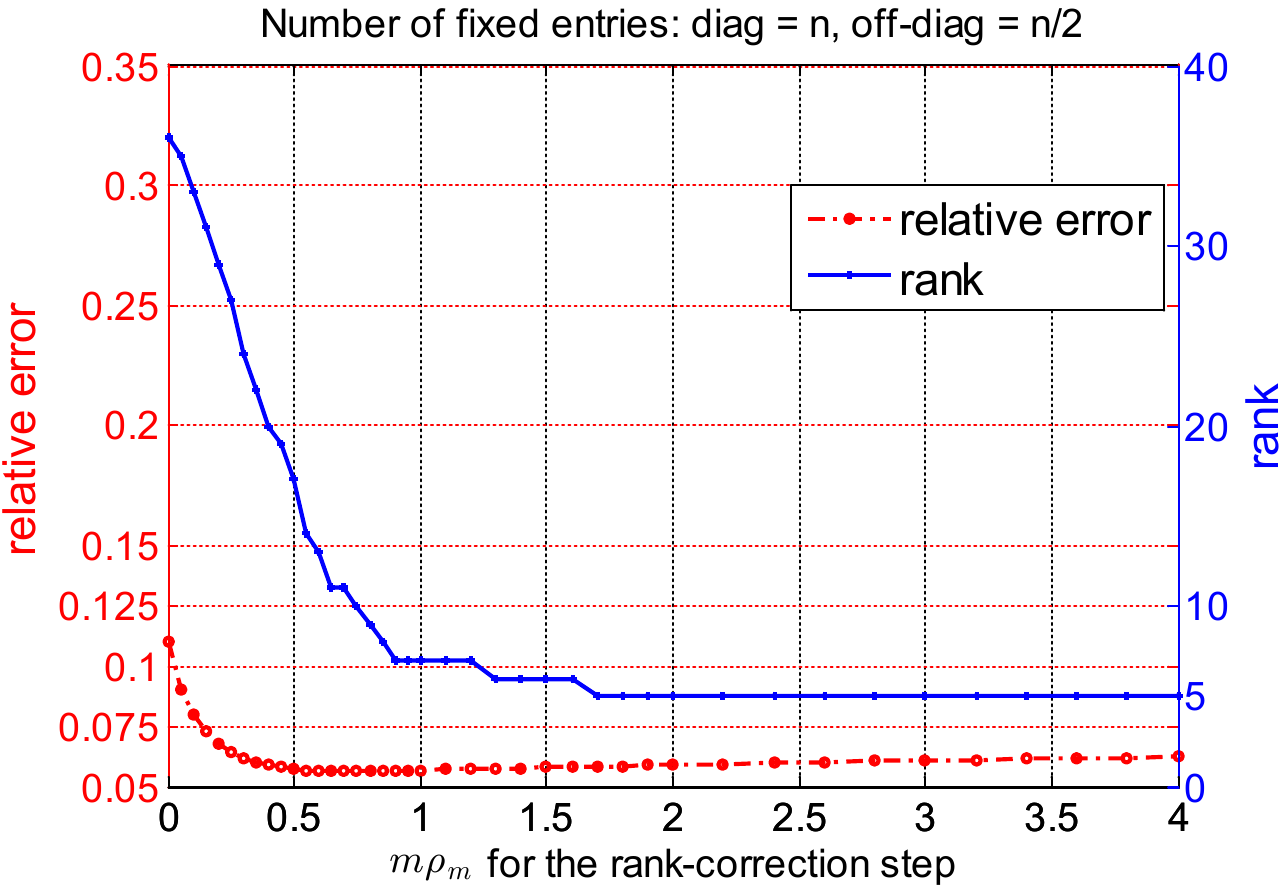}}\\
  \caption{Influence of fixed basis coefficients on recovery (sample ratio $=6.4\%$)\  \label{figure2}}
  \end{center}
\end{figure}

 In Figure \ref{figure2}, we plot the curves of the relative recovery error  and the rank of both the NNPLS estimator (the subfigures on the left) and the RCS estimator (the subfigures on the rigth) for different patterns of fixed entries. Note that both $m$ and $\rho_m$ in the rank-correction step (\ref{eqnrcs}) depend on the problem of consideration. Thus, we report $m\rho_m$ as a whole in the $x$-axis. (Note that for a specific problem, only $\rho_m$ is adjustable.) In the captions of  subfigures, {\bf diag} means the number of fixed diagonal entries, and {\bf off-diag} means the number of  fixed off-diagonal entries. For each subfigure on the right side, the initial $\widetilde{X}_m$ for the RCS estimator is the point with the smallest recovery error from the corresponding subfigure on the left side.

Figure \ref{figure2} fully manifests the advantage of the RCS estimator over the NNPLS estimator. It is shown that compared with the NNPLS estimator, the RCS estimator substantially reduces the recovery error and significantly improves the rank consistency. Moreover, the RCS estimator possesses a wide rage of the parameter $\rho_m$ to achieve a desired small recovery error and the rank of the true matrix simultaneously. It indicates that whether the resulting solution of a parameter $\rho_m$ achieves the true rank can be used to infer the recovery quality. Even if the true rank is unknown in advance, it is still possible to pick out a satisfied solution via monitoring the change of rank in parameter searching. Such advantages are far beyond the reach of the NNPLS estimator.

%--------------------------------------------------------------------------------------------Subsection6.2
\subsection{Performance of different rank-correction functions for recovery}

\begin{figure}[htbp]
\begin{center}
  \subfigure{\includegraphics[width=\textwidth]{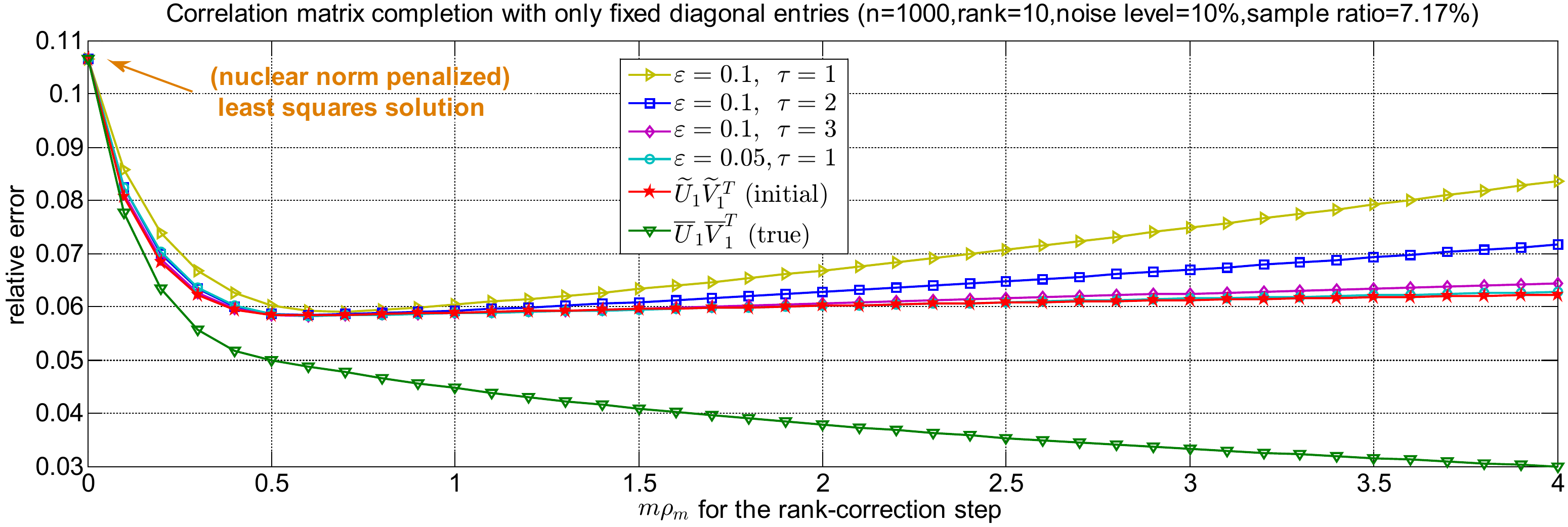}}\\
  \subfigure{\includegraphics[width=\textwidth]{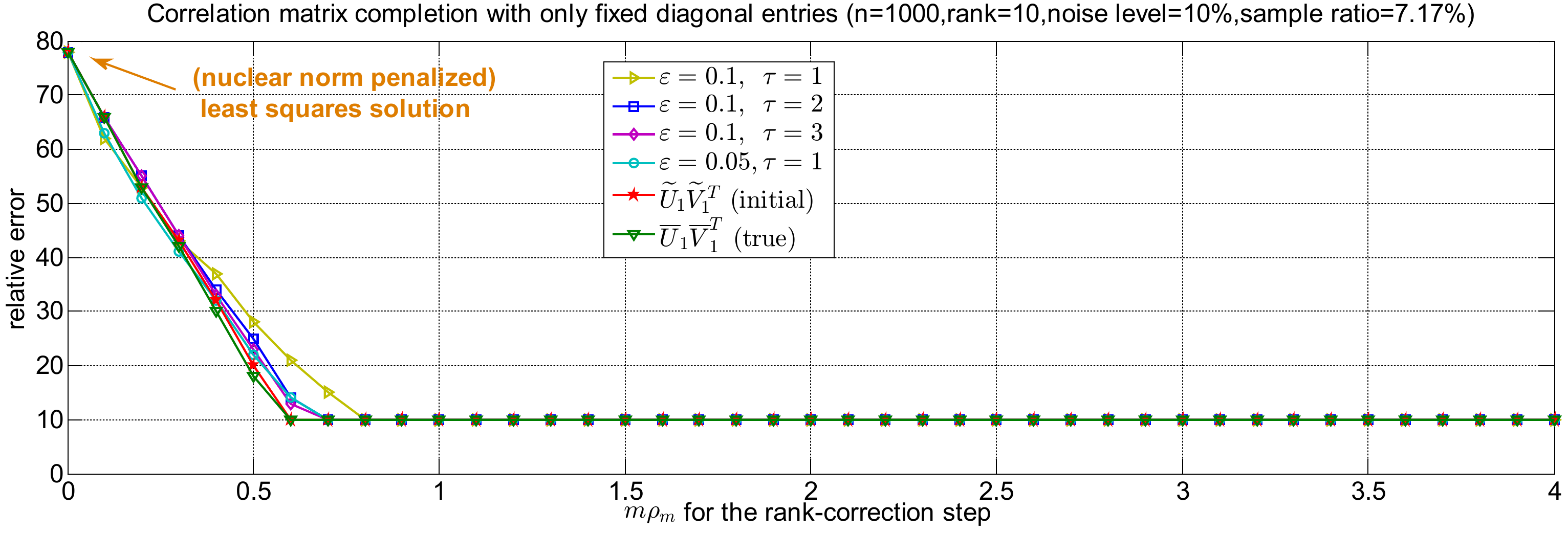}}
  \caption{Influence of the rank-correction term on the recovery\  \label{figure3}}
\end{center}
\end{figure}

 In this subsection, we test the performance of different rank-correction functions for recovering a correlation matrix. We randomly generated the true matrix $\overline{X}$ by the command in Subsection \ref{subsec6.1} with {\ttfamily n} $=1000$, {\ttfamily r} $=10$, {\ttfamily weight} $=2$ and {\ttfamily k} $ = 5$.
 We fixed all the diagonal entries of $\overline{X}$ and then sampled partial off-diagonal entries uniformly at random with i.i.d. Gaussian noise. The noise level was set to be $10\%$. We chose the (nuclear norm penalized) least squares estimator to be the initial estimator  $\widetilde{X}_m$. In Figure \ref{figure3}, we plot four curves corresponding to the rank-correction functions $F$ defined by (\ref{eqnchofun}), (\ref{eqnchofun2}) and (\ref{eqnchofun3}) with different $\varepsilon$ and $\tau $, and additional two curves corresponding to the rank-correction functions $F$ defined by (\ref{eqnchofun1}) at $\widetilde{X}_m$ (i.e., $\widetilde{U}_1\widetilde{V}_1^\mathbb{T}$) and $\overline{X}$ (i.e., $\overline{U}_1\overline{V}_1^\mathbb{T}$), respectively. The values of $a_m$ and the best recovery error are listed in Table \ref{tab1}.

  For all the rank-correction functions plotted in Figure \ref{figure3}, when $\rho_m$ increases, the recovery error first decreases together with the rank and then increases after the rank of the true matrix is attained. The only exception is $\overline{U}_1\overline{V}_1^\mathbb{T}$. This exactly validates our discussion about the recovery error at the end of Section \ref{section3}.
  %Moreover, for a smaller $\varepsilon$, the curve of recovery error changes more gently, though a certain optimality in the sense of %recovery error is sacrificed. This means that the choice of a relatively small $\varepsilon$,  say $0.01$ or $0.02$, is more robust for %those ill-conditioned problems.
  %From Table \ref{tab1}, we see that a smaller $a_m/b_m$ corresponds to a better optimal recovery error.
  It is worthwhile to point out that, according to our observations of many tests, in practice, if $a_m$ is larger than $1$ but not too much, the recovery performance of the RCS estimator still has a high chance to be much better than that of the NNPLS estimator.

\begin{table}[htbp]
 \begin{center}
 \caption{Influence of the rank-correction term on the recovery error\ \ \label{tab1}}
 \vspace{0.1cm}
 \begin{tabular}{|c|c|c|c|c|c|c|c|}
  \hline
  & zero & $\varepsilon = 0.1$ & $\varepsilon = 0.1$ & $\varepsilon = 0.1$ & $\varepsilon = 0.05$ & & \\
 \raisebox{1.5ex}[0pt]{$F$} & function &  $\tau = 1$ & $\tau = 2$ & $\tau = 3$ & $\tau = 2$ & \raisebox{1.5ex}[0pt]{$\widetilde{U}_1\widetilde{V}_1^{\mathbb{T}}$} & \raisebox{1.5ex}[0pt]{$\overline{U}_1\overline{V}_1^{\mathbb{T}}$}\\
  \hline
 $a_m$       & $1$       & $0.3126$ & $0.1652$ & $0.1402$ & $0.1849$ & $0.1355$  & 0\\
  \hline
 optimal relerr & $10.66\%$ & $5.92\%$ & $5.84\%$ & $5.83\%$ & $5.83\%$ & $5.84\%$  & $3.00\%$   \\
  \hline
\end{tabular}
\end{center}
\end{table}

\subsection{Performance of  different initial NNPLS estimators for recovery}
In this subsection, we take the covariance matrix completion for example to test the performance of the RCS estimator with different initial NNPLS estimators $\widetilde{X}_m$. We generated the true matrix $\overline{X}$ by the command in Subsection \ref{subsec6.1} with {\ttfamily n} $=500$, {\ttfamily r} $=5$, {\ttfamily weight} $=3$ and {\ttfamily k} $=1$ except that {\ttfamily D = eye(n)}. The upper bound of the non-fixed diagonal entries was set to be double of the largest absolute value among all the noisy observations of entries together with the fixed entries. We assumed that the rank of the true matrix was known so that we chose the rank-correction function  (\ref{eqnchofun1}).

\begin{figure}[htbp]
 \begin{center}
 \includegraphics[width=\textwidth]{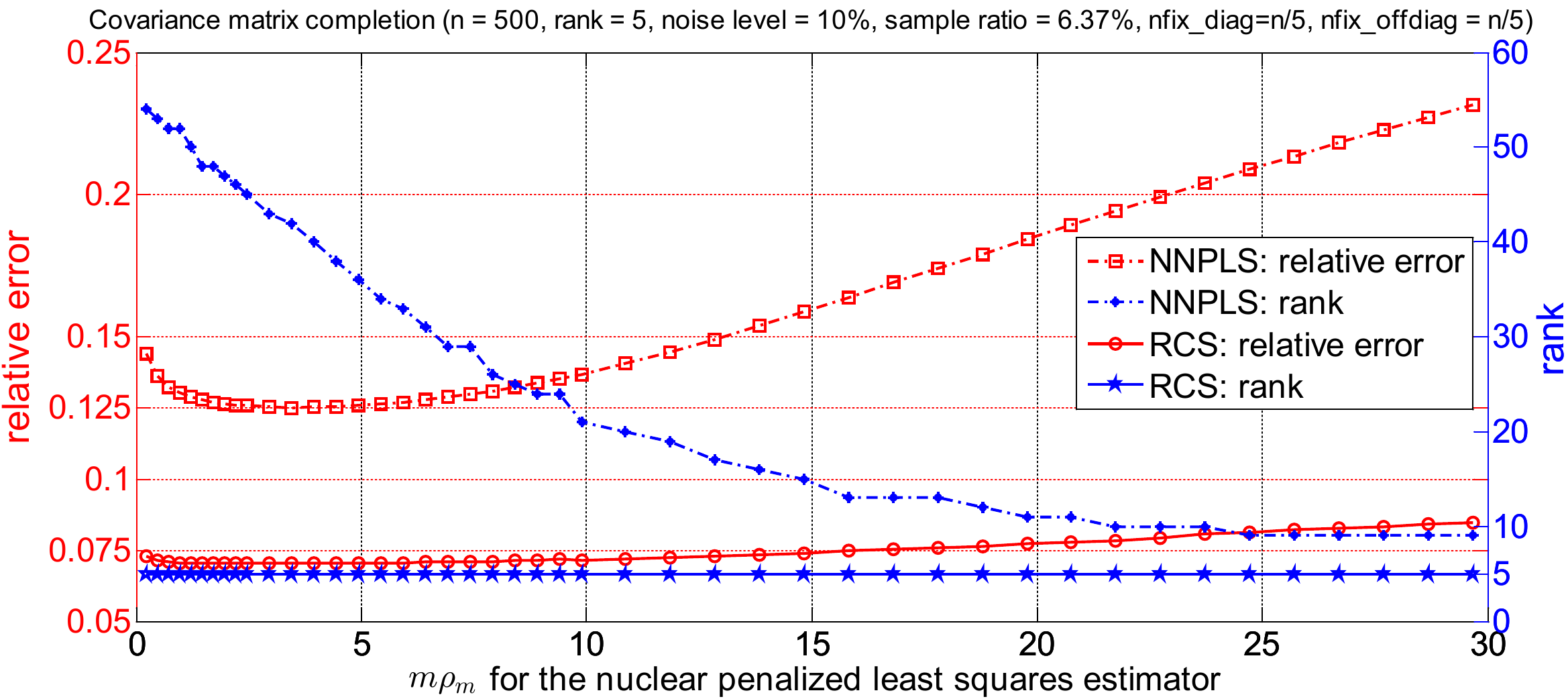}
 \caption{Performance of the RCS estimator with different initial NNPLS estimators \label{figure4}}
 \end{center}
 \end{figure}

For each $\rho_m$, we first produced the NNPLS estimator, and then use it as the initial point to produce a sequence of RCS estimators with different penalty parameters. Next we choose the RCS estimators that attains the correct rank with the smallest penalty parameter. As can be seen from Figure \ref{figure2}, this choice of the RCS estimator results in the desired small recovery error. The test results are plotted in Figure \ref{figure4}, where the dash curves represent for the NNPLS estimator and the solid curves represent for the chosen RCS estimator. We clearly observe from Figure \ref{figure4} that, no matter which NNPLS estimator is given to be the initial estimator, the RCS estimator can always substantially improve the recovery quality in terms of both the error and the rank.

%-------------------------------------------------------------------------------------------------Subsection 6
\subsection{Performance for different matrix completion problems}

   In this subsection, we test the performance of the RCS estimator for different matrix completion problems. Figure \ref{figure2} has revealed that a good choice of the parameter $\rho_m$ for the RCS estimator could be the smallest value that attains a stable rank. Therefore, the bisection search method can be used to find such a parameter $\rho_m$. This is actually what we benefit from rank consistency. In the following experiments, we apply this strategy to find a suitable $\rho_m$ for the RCS estimator.

  A natural question then arises: Will multiple rank-correction steps further improve the recovery quality? The answer can be found in Tables \ref{tab2}, \ref{tab3} and \ref{tab4} below, which report the experimental results for covariance matrix completion, rectangular matrix completion and density matrix completion, respectively. The reported NNPLS estimator is the one with the smallest recovery error among all different $\rho_m$ presuming the true matrix is known. The initial estimator of the first RCS estimator is the NNPLS estimator with a single preset $\rho_m = 0.4  \frac{\eta\|y\|_2}{\sqrt{m}} \sqrt{\frac{\log(n_1+n_2)}{mn}}$, where $\eta$ is the noise level. This choice of $\rho_m$ follows (\ref{bestrho}) with $C = 0.4$, $\kappa=1$, $\mu=1$ and $\nu$ taken its expected value based on observations. The second (third) RCS estimator takes the first (second) RCS estimator to be the initial estimator. The rank-correction function $F$ is defined by (\ref{eqnchofun}), (\ref{eqnchofun2}) and (\ref{eqnchofun3}) with $\varepsilon = 0.05$ and $\tau=2$.

\begin{table}[htbp]
  \begin{center}
  \renewcommand\arraystretch{1.2}
   {\caption{\label{tab2} Performance for covariance matrix completion problems with $n=1000$}}
   \vspace{0.1cm}
   \begin{tabular}{|c|c|c|c|c|c|c|}
  \hline
   &  &   &  {\rm NNPLS}& {\rm 1st RCS}&  {\rm 2st RCS}& {\rm 3rd RCS}\\
   \cline{4-7}
  \raisebox{1.5ex}[0pt]{\!$r$\!} & \raisebox{1.5ex}[0pt]{$\renewcommand{\arraystretch}{0.85} \begin{array}{c} {\rm  \!\!\!diag/\!\!\!} \\  {\rm \!\!\!off\!\!-\!\!diag\!\!\!} \end{array}$}   & \raisebox{1.5ex}[0pt]{$\renewcommand{\arraystretch}{0.85} \begin{array}{c} {\rm  \!\!sample\!\!} \\  {\rm \!\!ratio\!\!} \end{array}$} &{\rm \!relerr (rank)\!}&{\rm \!relerr(rank)\!}&{\rm \!relerr (rank)\!}&{\rm \!relerr (rank)\!}\\
  \hline \hline
     & \!1000/0\! & 2.40\% & 1.94e-1 (47) & 8.84e-2 (5) & 8.03e-2 (5) & 7.85e-2 (5)\\
     & \!1000/0\! & 7.99\% & 6.08e-2 (50) & 3.39e-2 (5) & 3.38e-2 (5) & 3.38e-2 (5)\\
  \cline{2-7}
   \raisebox{2.25ex}[0pt]{\!5\!}  & \!500/500\! & 2.39\% & 2.28e-1 (56) & 1.07e-1 (5) & 8.99e-2 (5) & 8.48e-2 (5)\\
     & \!500/500\! & 7.98\% & 1.16e-1 (56) & 5.62e-2 (5) & 5.42e-2 (5) & 5.40e-2 (5)\\
  \hline
     & \!1000/0\! & 5.38\%  & 1.59e-1 (77) & 7.42e-2 (10) & 7.23e-2 (10) & 7.22e-2 (10)\\
     & \!1000/0\! & 8.96\% & 9.15e-2 (81) & 5.06e-2 (10) & 5.05e-2 (10) & 5.05e-2 (10)\\
  \cline{2-7}
   \raisebox{2.25ex}[0pt]{\!10\!}  & \!500/500\! & 5.38\%  & 1.65e-1 (82) & 7.70e-2 (10) & 7.29e-2 (10) & 7.28e-2 (10)\\
     & \!500/500\! & 8.96\%  & 9.54e-2 (85) & 5.16e-2 (10) & 5.11e-2 (10) & 5.11e-2 (10)\\
  \hline
  \end{tabular}
  \end{center}
  \end{table}

 For the covariance matrix completion problems, we generated the true matrix $\overline{X}$ by the command in Subsection \ref{subsec6.1} with {\ttfamily n} $=1000$, {\ttfamily weight} $=2$ and {\ttfamily k} $=1$ except that {\ttfamily D = eye(n)}. The rank of $\overline{X}$ and
 the number of fixed diagonal and non-diagonal entries of $\overline{X}$ are reported in the first and the second columns of Table \ref{tab2}, respectively. We sampled partial off-diagonal entries uniformly at random with i.i.d. Gaussian noise at the noise level $10\%$.  The upper bound of the non-fixed diagonal entries was set to be double of the largest absolution value among all the noisy observations of entries together with the fixed entries. From Table \ref{tab2}, we see that when the sample ratio is reasonable, a single rank-correction step is fully capable to yield a desired result. However, when the sample ratio is very low, especially if some off-diagonal entries are fixed, one or two further rank-correction steps could still bring some improvement in recovery quality.

For the density matrix completion problems, we generated the true density matrix $\overline{X}$ by the following command:
 \vspace{-0.3cm}
 \begin{verbatim}
   M = randn(n,r)+i*randn(n,r);  ML = weight*M(:,1:k);  M(:,1:k) = ML;
              Xtemp = M*M';  X_bar = Xtemp/sum(diag((Xtemp))).
 \end{verbatim}
 \vspace{-0.8cm}
 During the testing, we set {\ttfamily n} $=1024$, {\ttfamily weight} $=2$ and {\ttfamily k} $=1$, and sampled
 partial Pauli measurements except the trace of $\overline{X}$ uniformly at random with $10\%$ i.i.d. Gaussian noise.
 Besides this statistical noise, we further added the depolarizing noise, which frequently appears in quantum systems. The strength of the depolarizing noise was set to be $0.01$. This case is labeled as the mixed noise in the last four rows of Table \ref{tab3}. We remark here that the depolarizing noise differs from our assumption on noise since it does not have randomness. One may refer to \cite{GroLFBE10, FlaGLE12} for details of the quantum depolarizing channel. In \cite{FlaGLE12}, Flammia et al. proposed a two-step method for seeking a feasible solution of low-rank --- (1) evaluating an NNPLS estimator by dropping the trace one constraint; (2) normalizing the resulting solution to be of trace one. We tested this method in our experiments, with the NNPLS estimator without trace one constraint chosen to be the one with the smallest recovery error among all that attain the true rank, presuming that the true matrix is known. The two-step results are reported as NNPLS1 and NNPLS2, respectively, in Table \ref{tab3}. Besides the relative recovery error ({\bf relerr}), we also report the (squared) {\bf fidelity}, which is a measure of the closeness of two quantum states defined by $\big\|\widehat{X}_m^{1/2} \overline{X}^{1/2}\big\|_{*}^2$.
 From Table \ref{tab3}, we can see that the RCS estimator is superior to the NNPLS2 estimator in terms of both
 the fidelity and the relative error.

  \begin{table}[htbp]
  \begin{center}
  \renewcommand\arraystretch{1.2 }
   {\caption{\label{tab3} Performance for density matrix completion problems with $n=1024$}}
   \vspace{0.1cm}
   \begin{tabular}{|c|c|c|c|c|c|c|}
   \hline
  \multirow{2}{*}{\rotatebox{90}{noise\ \ }}& \multirow{2}{*}{\!$r$\!}& & &{\rm NNPLS1}&{\rm NNPLS2}&{\rm RCS}\\
   \cline{5-7}
   & & \raisebox{1.5ex}[0pt]{$\renewcommand{\arraystretch}{0.85} \begin{array}{c} {\rm \!\!\!\! noise \!\!\!\!} \\  {\rm \!\!\!\!level\!\!\!\!} \end{array}$} & \raisebox{1.5ex}[0pt]{$\renewcommand{\arraystretch}{0.85} \begin{array}{c} {\rm  \!\!\!\!sample\!\!\!\!} \\  {\rm \!\! ratio\!\!} \end{array}$} & {\rm \!\!\!\! fidelity \!\! \!}\ \  {\rm \!\!relerr\!\!}\ \ {\rm  rank \!\!\!\!}& {\rm \!\!\!\! fidelity \!\! \!}\ \  {\rm \!\!relerr\!\!}\ \ {\rm  rank \!\!\!\!} & {\rm \!\!\!\! fidelity \!\! \!}\ \  {\rm \!\!relerr\!\!}\ \ {\rm  rank \!\!\!\!}\\
   \hline \hline
   \multirow{4}{*}{\rotatebox{90}{statistical\ \ }}
   &                         & \!10.0\%\! & 1.5\%  & 0.716\ \ \ \!\!2.49e-1\ \ 3\ \ & 0.962\ \ \ \!\!2.34e-1\ \  3\ & 0.992\ \ \ \!\!8.47e-2\ \  3\ \\
   &\raisebox{2.25ex}[0pt]{\!3\!} & \!10.0\%\! & 4.0\% & 0.915\ \ \ \!\!8.14e-2\ \ 3\ \ & 0.997\ \ \ \!\!6.88e-2\ \  3\ & 0.998\ \ \ \!\!4.13e-2\ \  3\ \\
   \cline{2-7}
   &                         & \!10.0\%\! & 2.5\%  & 0.696\ \ \ \!\!2.56e-1\ \ 5\ \ & 0.959\ \ \ \!\!2.71e-1\ \  5\ & 0.992\ \ \ \!\!8.28e-2\ \  5\ \\
   &\raisebox{2.25ex}[0pt]{\!5\!} & \!10.0\%\! & 5.0\% & 0.886\ \ \ \!\!1.04e-1\ \ 5\ \ & 0.994\ \ \ \!\!9.61e-2\ \  5\ & 0.997\ \ \ \!\!4.81e-2\ \  5\ \\
   \hline \hline
   \multirow{4}{*}{\rotatebox{90}{mixed\ \ }}
   &                         & \!12.5\%\! & 1.5\%  & 0.657\ \ \ \!\!2.95e-1\ \ 3\ \ & 0.959\ \ \ \!\!2.41e-1\ \  3\ & 0.990\ \ \ \!\!9.89e-2\ \  3\ \\
   &\raisebox{2.25ex}[0pt]{\!3\!} & \!12.4\%\!  & 4.0\% & 0.842\ \ \ \!\!1.42e-1\ \  3\ \ & 0.996\ \ \ \!\!7.48e-2\ \  3\ & 0.997\ \ \ \!\!6.20e-2\ \  3\ \\
   \cline{2-7}
   &                         & \!12.4\%\! & 2.5\%  & 0.631\ \ \ \!\!3.05e-1\ \  5\ \ & 0.954\ \ \ \!\!2.87e-1\ \ 5\ & 0.990\ \ \ \!\!9.81e-2\ \  5\ \\
   &\raisebox{2.25ex}[0pt]{\!5\!} & \!12.5\%\! & 5.0\% & 0.814\ \ \ \!\!1.62e-1\ \  5\ \ & 0.994\ \ \ \!\!1.03e-1\ \  5\ & 0.996\ \ \ \!\!6.94e-2\ \  5\ \\
    \hline
  \end{tabular}
  \end{center}
  \end{table}

For the rectangular matrix completion problems, we generated the true matrix $\overline{X}$ by the following command:
\vspace{-0.3cm}
\begin{verbatim}
   ML = randn(nr,r);   MR = randn(nc,r);   MW = weight*ML(:,1:k);
               ML(:,1:k) = MW;   X_bar = ML*MR'.
 \end{verbatim}
 \vspace{-0.8cm}
 \noindent
 We set {\ttfamily weight} $=2$, {\ttfamily k} $=1$ and took $\overline{X}=$ {\ttfamily X\_bar} with different dimensions and ranks. Both the uniform sampling scheme and the non-uniform sampling scheme were tested for comparison. For the non-uniform sampling scheme, the probability to sample the first $1/4$ rows and the first $1/4$ columns were  $3$ times as much as that of other rows and columns respectively. In other words, the density of sampled entries in the top-left part was $3$ times as much as that in the bottom-left part and the top-right part respectively and $9$ times as much as that in the bottom-right part. We added $10\%$ i.i.d. Gaussian noise to the sampled entries. We also fixed partial entries of $\overline{X}$ uniformly from the rest un-sampled entries. The upper bound of the non-fixed entries was set to be double of the largest absolution value among all the noisy observations of entries together with the fixed entries. What we observe from Table \ref{tab4} for the rectangular matrix completion is similar to that for the covariance matric completion. Moreover, we can see that the non-uniform sampling scheme greatly weakens the recoverability of the NNPLS estimator in terms of both the recovery error and the rank, especially when the sample ratio is low. Meanwhile, the advantage of the RCS estimators in such  cases becomes more remarkable.

\medskip

\begin{table}[htbp]
  \begin{center}
  \renewcommand\arraystretch{1.2}
   {\caption{\label{tab4} Performance for rectangular matrix completion problems}}
   \vspace{0.1cm}
   \begin{tabular}{|c|c|c|c|c|c|c|c|}
  \hline
   \multirow{2}{*}{\rotatebox{90}{{\rm setting}}} & \multirow{2}{*}{\rotatebox{90}{{\rm sample}}} &  &   &  {\rm NNPLS}& {\rm 1st RCS}&  {\rm 2st RCS}& {\rm 3rd RCS}\\
   \cline{5-8}
  \raisebox{1.75ex}[0pt]{ } & \raisebox{1.5ex}[0pt]{} & \raisebox{1.5ex}[0pt]{\!{\rm fixed}\!} & \raisebox{1.5ex}[0pt]{$\renewcommand{\arraystretch}{0.85} \begin{array}{c} {\rm  \!\!sample\!\!} \\  {\rm \!\!ratio\!\!} \end{array}$} &{\rm \!relerr (rank)\!}&{\rm \!relerr (rank)\!}&{\rm \!relerr (rank)\!}&{\rm \!relerr (rank)\!}\\
  \hline \hline
   \multirow{8}{*}{\rotatebox{90}{\ \ \ \ \ \ $\text{dim}\!=\!1000\!\times\!1000, \text{rank}\!=\!10\!\!$ }}
     & \multirow{4}{*}{\rotatebox{90}{\rm uniform\ \ }}
        & 0    &  5.97\%  & 1.98e-1 \!(119)\!\! & 7.69e-2 (10)  & 7.31e-2 (10)  & 7.30e-2 (10) \\
     &  & 0    &  11.9\%  & 8.34e-2 \!(114)\!\! & 4.49e-2 (10)  & 4.48e-2 (10)  & 4.48e-2 (10) \\
   \cline{3-8}
     &  & 1000 &  5.98\%  & 1.93e-1 \!(120)\!\! & 7.45e-2 (10)  & 7.01e-2 (10)  & 7.00e-2 (10) \\
     &  & 1000 &  12.0\%  & 8.20e-2 \!(108)\!\! & 4.35e-2 (10)  & 4.34e-2 (10)  & 4.34e-2 (10) \\
  \cline{2-8}
     & \multirow{4}{*}{\rotatebox{90}{\rm non-uniform\  }}
        & 0    &  5.97\%  & 3.20e-1 \!(144)\!\! & 1.22e-1 (10)  & 9.31e-2 (10)  & 8.77e-2 (10) \\
     &  & 0    &  11.9\%  & 1.27e-1 \!(171)\!\! & 5.32e-2 (10)  & 5.12e-2 (10)  & 5.11e-2 (10) \\
  \cline{3-8}
     &  & 1000 &  5.98\%  & 3.07e-1 \!(146)\!\! & 1.16e-1 (10)  & 8.78e-2 (10)  & 8.30e-2 (10) \\
     &  & 1000 &  12.0\%  & 1.24e-1 \!(173)\!\! & 5.14e-2 (10)  & 4.93e-2 (10)  & 4.92e-2 (10) \\
  \hline
     \multirow{8}{*}{\rotatebox{90}{ \ \ \ \ $\text{dim}\!=\!500\!\times\!1500, \text{rank}\!=\!5$  }}
     & \multirow{4}{*}{\rotatebox{90}{\rm uniform\ \ }}
        & 0    &  3.99\%  & 2.31e-1 (73) & 9.15e-2 (5)  & 8.10e-2 (5)  & 7.95e-2 (5) \\
     &  & 0    &  7.98\%  & 9.01e-2 (78) & 4.60e-2 (5)  & 4.58e-2 (5)  & 4.58e-2 (5) \\
   \cline{3-8}
     &  & 1000 &  4.00\%  & 2.15e-1 (74) & 8.77e-2 (5)  & 7.58e-2 (5)  & 7.36e-2 (5) \\
     &  & 1000 &  7.99\%  & 8.72e-2 (69) & 4.34e-2 (5)  & 4.31e-2 (5)  & 4.31e-2 (5) \\
  \cline{2-8}
     & \multirow{4}{*}{\rotatebox{90}{\rm non-uniform\ }}
        & 0    &  3.99\%  & 3.37e-1 (91) & 1.53e-1 (5)  & 1.18e-1 (5)  & 1.07e-1 (5) \\
     &  & 0    &  7.98\%  & 1.37e-1 \!(128)\!\! & 5.62e-2 (5)  & 5.33e-2 (5)  & 5.31e-2 (5) \\
  \cline{3-8}
     &  & 1000 &  4.00\%  & 3.11e-1 (93) & 1.39e-1 (6)  & 1.06e-1 (5)  & 9.55e-2 (5) \\
     &  & 1000 &  7.99\%  & 1.29e-1 \!(104)\!\! & 5.21e-2 (5)  & 4.91e-2 (5)  & 4.89e-2 (5) \\
  \hline
  \end{tabular}
  \end{center}
  \end{table}

%&  3.99\%  & 2.24e-1 (43) & 8.75e-2 (5)  & 7.93e-2 (5)  & 7.76e-2 (5) \\
%&  4.00\%  & 2.23e-1 (49) & 8.77e-2 (6)  & 7.67e-2 (5)  & 7.44e-2 (5) \\
%&  3.99\%  & 3.29e-1 (57) & 2.64e-1 (5)  & 1.39e-1 (5)  & 1.23e-1 (5) \\

\section{Conclusions}\label{section7}

 In this paper, we proposed a rank-corrected procedure for low-rank matrix completion problems with fixed basis coefficients. This approach can substantially overcome the limitation of the nuclear norm technique for recovering a low-rank matrix. We confirmed the improvement of the rank-correction step in both the reduction of recovery error and the achievement of rank consistency (in the sense of Bach \cite{Bac08}). Due to the presence of fixed basis coefficients, constraint nondegeneracy plays an important role in our analysis. Extensive numerical experiments show that our approach can significantly improve the recovery performance compared with the nuclear norm penalized least square estimator. As a byproduct, our results also provide a theoretical foundation for the majorized penalty method of Gao and Sun
 \cite{GaoS10} and Gao \cite{Gao10} for structured low-rank matrix optimization problems.

 Our proposed rank-correction step also allows additional constraints according to other possible prior information. In order to better fit the under-sampling setting of matrix completion, in the future work, it would be of great interest to extend the asymptotic rank consistency results to the case where the   matrix size is allowed to grow.  It would also be interesting to extend this approach to deal with other low-rank matrix problems.

\section*{Acknowledgements}

The authors would like to thank Professor  Wotao Yin for his valuable comments on  possibly choosing the optimal penalty parameter for recovery error bounds  and Dr. Kaifeng Jiang for helpful discussions  on efficiently solving  the density matrix completion problem.

\appendix
\gdef\thesection{Appendix \Alph{section}}

\section{Spectral operator}\label{AppendixSpeOpe}

The concept of spectral operator is associated with
 a symmetric vector-valued function.
 A function $f: \mathbb{R}^n \rightarrow \mathbb{R}^n$ is said to be symmetric if
 $$f(x) = Q^{\mathbb{T}} f(Qx) \quad \forall \, \text{signed permutation matrix} \ Q \ \text{and} \ x \in \mathbb{R}^n,$$
  where a signed permutation matrix
 is a real matrix that contains exactly one nonzero entry $1$ or $-1$ in each row and column and $0$ elsewhere.
 From this definition, we see that $$f_i(x) = 0 \quad \text{if} \  x_i=0.$$
 The spectral operator $F: \mathbb{V}^{n_1\times n_2}\rightarrow \mathbb{V}^{n_1\times n_2}$ associated with the function
 $f$ is defined by
 \begin{align}\label{Foperator}
  F(X): = U \text{Diag}\big(f(\sigma(X))\big) V^\mathbb{T},
 \end{align}
 where $(U,V)\in \mathbb{O}^{n_1,n_2}(X)$ and $X \in \mathbb{V}^{n_1\times n_2}$. From \cite[Theorems 3.1 \& 3.6]{Din12}, the symmetry of $f$ guarantees
 the well-definiteness of the spectral operator $F$, and the (continuous) differentiability of $f$ implies the (continuous)
 differentiability of $F$. When $\mathbb{V}^{n_1\times n_2}=\mathbb{S}^n$, we have that
 $$F(X) = P\text{Diag}\big(f(|\lambda(X)|)\big)\big(P\text{Diag}(s(X))\big)^\mathbb{T},$$ where $P \in \mathbb{O}^n(X)$
 and $s(X)\in\mathbb{R}^n$ with the $i$-th component $s_i(X)=-1$ if $\lambda_i(X)<0$ and $s_i(X)=1$ otherwise.
 In particular for the positive semidefinite case, both $U$ and $V$ in (\ref{Foperator}) reduce to $P$.
 For more details on spectral operators, the readers may refer to the PhD thesis \cite{Din12}.
% the spectral operator can also be equivalently written
% as $$F(X) = P\text{Diag}\big(f(|\lambda(X)|)\big)\text{Diag}\big(\text{sgn}(\lambda(X))\big)P^\mathbb{T},$$ with $P \in \mathbb{O}^n(X)$ and % $\text{sgn}(\lambda(X)) = \big(\text{sgn}(\lambda_1(X)),\ldots \text{sgn}(\lambda_n(X)))^\mathbb{T}\big)$.

\section{Constraint nondegeneracy}\label{AppendixConNon}
Consider the following constrained optimization problem
 \begin{equation}\label{conic-prob}
   \min_{X\in\mathbb{V}^{n_1\times n_2}}\Big\{\Phi(X)+\Psi(X):\ \mathcal{A}(X)-b\in K\Big\},
 \end{equation}
 where $\Phi:\mathbb{V}^{n_1\times n_2} \to\mathbb{R}$ is a continuously differentiable function,
 $\Psi:\mathbb{V}^{n_1\times n_2}\to\mathbb{R}$ is a convex function, $\mathcal{A}: \mathbb{V}^{n_1\times n_2}\to\mathbb{R}^l$ is
 a linear operator and $K \subseteq \mathbb{R}^l$ is a closed
 convex set. Let $\widehat{X}$ be a given feasible point
 of (\ref{conic-prob}) and $\widehat{z}:=\mathcal{A}(\widehat{X})-b$.
 When $\Psi$ is differentiable at $\widehat{X}$,  we say that
 the constraint nondegeneracy holds at $\widehat{X}$ if
 \begin{equation}\label{nondegeneracy1}
  \mathcal{A}\,\mathbb{V}^{n_1\times n_2} + {\rm lin}\big(\mathcal{T}_{K}(\widehat{z})\big) = \mathbb{R}^l,
 \end{equation}
 where $\mathcal{T}_{K}(\widehat{z})$ denotes the tangent cone of $K$ at $\widehat{z}$
 and ${\rm lin}(\mathcal{T}_{K}(\widehat{z}))$ denotes the largest linearity space
 contained in $\mathcal{T}_{K}(\widehat{z})$, i.e.,
 ${\rm lin}(\mathcal{T}_{K}(\widehat{z}))=\mathcal{T}_{K}(\widehat{z})\cap(-\mathcal{T}_{K}(\widehat{z}))$.
 When the function $\Psi$ is nondifferentiable,  we can rewrite the optimization problem (\ref{conic-prob}) equivalently as
$$ \min_{X\in\mathbb{V}^{n_1\times n_2},t\in\mathbb{R}}\Big\{\Phi(X)+t:\ \widetilde{\mathcal{A}}(X,t)\in  K \times {\rm epi}\Psi \Big\},$$
where ${\rm epi}\Psi:=\left\{(X,t)\in\mathbb{V}^{n_1\times n_2}\times\mathbb{R}\ |\ \Psi(X)\le t\right\}$ denotes the epigraph of $\Psi$ and $\widetilde{\mathcal{A}}:\mathbb{V}^{n_1\times n_2}\times\mathbb{R}\rightarrow \mathbb{R}^l\times \mathbb{V}^{n_1\times n_2}\times\mathbb{R}$
is a linear operator defined by
$$\widetilde{\mathcal{A}}(X,t) := \begin{pmatrix}
                                   \mathcal{A}(X) - b\\ X\\ t
                                   \end{pmatrix}, \quad\   (X,t)\in  \mathbb{V}^{n_1\times n_2} \times \mathbb{R}.$$
From (\ref{nondegeneracy1}) and \cite[Theorem 6.41]{RocW98}, the constraint nondegeneracy holds at $(\widehat{X}, \widehat{t})$ with $\widehat{t}=\Psi(\widehat{X})$ if
$$\widetilde{\mathcal{A}}\begin{pmatrix}
                      \mathbb{V}^{n_1\times n_2} \\ \mathbb{R}
                \end{pmatrix}
    +\begin{pmatrix}
        {\rm lin}\big(\mathcal{T}_{K}(\widehat{z})\big) \\  {\rm lin}\big(\mathcal{T}_{{\rm epi}\Psi}(\widehat{X},\widehat{t})\big)
                \end{pmatrix}
     =  \begin{pmatrix}
               \mathbb{R}^l\\ \mathbb{V}^{n_1\times n_2}\\ \mathbb{R}
           \end{pmatrix}.$$
By the definition of $\widetilde{\mathcal{A}}$, it is not difficult to verify that this condition is equivalent to
\begin{equation}\label{nondegeneracy2}
[\mathcal{A}\ \ 0]\big({\rm lin}(\mathcal{T}_{{\rm epi}\Psi}(\widehat{X},\widehat{t}))\big)
+ {\rm lin}\big(\mathcal{T}_{K}(\widehat{z})\big) = \mathbb{R}^l.
\end{equation}

One can see that the problem (\ref{eqnrcs}) with $\mathcal{C} = \mathbb{V}^{n_1\times n_2}$ can be cast into (\ref{conic-prob}) with $\Psi=\|\cdot\|_*, \mathcal{A}=[\mathcal{R}_{\alpha} \ \mathcal{R}_{\beta}]$, $K\!=\{0\}^{|\alpha|} \times [-b, b]^{|\beta|}$, and meanwhile the problem  (\ref{eqnrcs}) with $\mathcal{C} = \mathbb{S}_+^n$ can be cast into (\ref{conic-prob}) with $\Psi =\delta_{\mathbb{S}_{+}^n}, \mathcal{A}=\mathcal{R}_{\alpha}$, $K=\{0\}$. In the previous case, the condition
(\ref{nondegeneracy2}) reduces to (\ref{eqncndc}) according to the expression of $\mathcal{T}_{{\rm epi}\Psi}(\overline{X},\overline{t})$ with $\overline{t}=\|\overline{X}\|_*$
(e.g., see \cite{JiaST12}). In the latter case, the condition (\ref{nondegeneracy2}) reduces to (\ref{eqncndcpos}) according to Arnold's characterization of  the tangent cone $\mathcal{T}_{\mathbb{S}^n_+}(\overline{X})=\big\{H \in\mathbb{S}^n \mid \overline{P}_2^\mathbb{T} H \overline{P}_2 \in \mathbb{S}_+^{n-r} \big\}$ in \cite{Arn71}.

\section{Proofs of Theorems}\label{AppendixTheorem}
\gdef\thesection{\Alph{section}}

\subsection{Proof of Theorem \ref{thmopbd}}

Let $\Delta_m:=\widehat{X}_m - \overline{X}$. Using the optimality of $\widehat{X}_m$ to the problem (\ref{eqnrcs}), we obtain that
\begin{equation}\label{eqnopbd1}
 \frac{1}{2m} \|\mathcal{R}_\Omega(\Delta_m)\|_2^2 \leq  \Big\langle \frac{\nu}{m}\mathcal{R}_\Omega^*(\xi), \Delta_m \Big\rangle - \rho_m\big(\|\widehat{X}_m\|_*-\|\overline{X}\|_* - \langle F(\widetilde{X}_m), \Delta_m\rangle\big).
\end{equation}
Then, we introduce an orthogonal decomposition $\mathbb{V}^{n_1\times n_2}=T\oplus T^\perp$  with
 \begin{equation*}
 \left\{
 \begin{aligned}
    &T:=   \big\{X\in \mathbb{V}^{n_1\times n_2} \ |\ X = X_1 + X_2\ {\rm with}\ {\rm col}(X_1)\subseteq {\rm col}(\overline{X}) \ {\rm and}\
    {\rm row}(X_2)\subseteq {\rm row}(\overline{X})\big\},\\
    &T^{\bot}:=   \big\{X\in \mathbb{V}^{n_1\times n_2}\ |\ {\rm row}(X)\perp {\rm row}(\overline{X})\ {\rm and}\
    {\rm col}(X)\perp {\rm col}(\overline{X})\big\},
 \end{aligned}
 \right.
 \end{equation*}
 where ${\rm row}(X)$ and ${\rm col}(X)$ denote the row space and column space of $X$, respectively. Let $\mathcal{P}_{T}$ and $\mathcal{P}_{T^\perp}$ be orthogonal projections onto $T$ and $T^\perp$, respectively, given by
 \begin{equation}\label{Operator-PT}
  \mathcal{P}_T(X) = \overline{U}_1\overline{U}_1^\mathbb{T}X + X \overline{V}_1\overline{V}_1^\mathbb{T} - \overline{U}_1\overline{U}_1^\mathbb{T}X \overline{V}_1 \overline{V}_1^\mathbb{T} \ \ \text{and}\ \
  \mathcal{P}_{T^\perp}(X) = \overline{U}_2\overline{U}_2^\mathbb{T} X \overline{V}_2\overline{V}_2^\mathbb{T}
 \end{equation}
 for any $X\in\mathbb{V}^{n_1\times n_2}$ and $(\overline{U},\overline{V})\in\mathbb{O}^{n_1,n_2}(\overline{X})$.
Then, it follows from the choice of $\rho_m$ that
\begin{equation}\label{eqnopbd3}
 \Big\langle \frac{\nu}{m}\mathcal{R}_\Omega^*(\xi), \Delta_m \Big\rangle
 \leq  \Big\|\frac{\nu}{m}\mathcal{R}_\Omega^*(\xi)\Big\| \|\Delta_m\|_*
 \leq \frac{\rho_m}{\kappa} \big(\|\mathcal{P}_T(\Delta_m)\|_* +\|\mathcal{P}_{T^\perp}(\Delta_m)\|_*\big).
\end{equation}
Moreover, from the directional derivative of the nuclear norm at $\overline{X}$, (see \cite[Theorem 1]{Wat92}), we have
\begin{align}
\|\widehat{X}_m \|_* - \|\overline{X}\|_* -\langle F(\widetilde{X}_m), \Delta_m\rangle \geq & \ \langle \overline{U}_1\overline{V}_1^\mathbb{T}, \Delta_m\rangle +\|\overline{U}_2^\mathbb{T}\Delta_m\overline{V}_2\|_* -\langle F(\widetilde{X}_m), \Delta_m\rangle \nonumber\\
  \geq & \ \|\mathcal{P}_{T^\perp}(\Delta_m)\|_* - \|\overline{U}_1\overline{V}_1^\mathbb{T} - F(\widetilde{X}_m)\|_F \|\Delta_m\|_F \nonumber \\
  = & \ \|\mathcal{P}_{T^\perp}(\Delta_m)\|_* - a_m \sqrt{r} \|\Delta_m \|_F. \label{eqnopbd2}
\end{align}
Then, by substituting (\ref{eqnopbd3}) and (\ref{eqnopbd2}) into  (\ref{eqnopbd1}), we have
\begin{equation}\label{eqnopbdtig}
\frac{1}{2m} \|\mathcal{R}_\Omega(\Delta_m)\|_2^2 \leq   \rho_m\Big(a_m \sqrt{r} \|\Delta_m \|_F +  \frac{1}{\kappa}\|\mathcal{P}_T(\Delta_m)\|_* - \frac{\kappa-1}{\kappa} \|\mathcal{P}_{T^\perp}(\Delta_m)\|_*\Big).
\end{equation}
Note that $\text{rank}(\mathcal{P}_T(\Delta_m))\leq 2r$. Hence, $\|\mathcal{P}_T(\Delta_m)\|_* \leq \sqrt{2r} \|\mathcal{P}_T(\Delta_m)\|_F \leq\sqrt{2r}\|\Delta_m\|_F$ and then the desired result (\ref{eqnopbd}) follows.

\subsection{Proof of Theorem \ref{thmstobd}}

\noindent We first show that the sampling operator $\mathcal{R}_\Omega$ satisfies some RIP-like property for matrices specified in a certain set with high probability. Similar results can also be found in \cite{NegW12, KolLT11, Klo12, Liu11}.

For this purpose, define
\begin{equation}\label{defvaritheta}
 \vartheta_m := \mathbb{E}\, \Big\|\frac{1}{m}\mathcal{R}_\Omega^*(\epsilon)\Big\|\ \ {\rm with}\ \
 \epsilon= (\epsilon_1,\ldots,\epsilon_m)^{\mathbb{T}},
 \end{equation}
where $\{\epsilon_1,\ldots, \epsilon_m\}$ is an i.i.d. Rademacher sequence, i.e., an i.i.d. sequence of Bernoulli random variables taking the values $1$ and $-1$ with probability $1/2$.

\begin{lemma}\label{lemiso}
Given any $s>0$ and $t>0$, define
\begin{align*}
K(s,t) := & \big\{\Delta \in \mathbb{V}^{n_1\times n_2} \ \big| \ \mathcal{R}_\alpha(\Delta) = 0, \|\mathcal{R}_\beta(\Delta)\|_\infty=1,  \|\Delta\|_* \leq s\|\Delta\|_F,  \langle \mathcal{Q}_\beta(\Delta), \Delta \rangle \geq t\big\}.
\end{align*}
Then, for any given $\gamma > 1$, $\tau_1 \in (0,1)$ and $\tau_2 \in (0, \tau_1/\gamma)$,
with probability at least $1-\frac{\exp(-(\tau_1-\gamma\tau_2)^2 m t^2/2)}{1-\exp(-(\gamma^2-1) (\tau_1-\gamma\tau_2)^2 m t^2/2)}$,
\begin{equation}\label{eqnlemiso}
 \frac{1}{m}\|\mathcal{R}_\Omega(\Delta)\|_2^2 \geq (1-\tau_1) \langle \mathcal{Q}_\beta(\Delta),\Delta\rangle - \frac{16}{\tau_2} s^2 \mu_1  d_2   \vartheta_m^2 \quad  \forall \, \Delta \in K(s,t).
\end{equation}
\end{lemma}
\begin{proof}
The proof is similar to that of \cite[Lemma 12]{Klo12}. For any $s, t>0$, $\gamma > 1$, $\tau_1 \in (0,1)$ and $\tau_2 \in (0, \tau_1/\gamma)$, we need to show that the event
\begin{align*}
E = \Big\{\exists \, \Delta \! \in \!  K(s,t) \, \text{such that}\, \Big|\frac{1}{m}\|\mathcal{R}_\Omega(\Delta)\|_2^2 - \langle \! \mathcal{Q}_\beta(\Delta),\Delta\rangle \Big| \geq \tau_1 \langle \mathcal{Q}_\beta(\Delta), \Delta\rangle + \frac{16}{\tau_2} s^2 \mu_1  d_2   \vartheta_m^2 \! \Big\}
\end{align*}
occurs with probability less than $\frac{\exp(-(\tau_1-\gamma\tau_2)^2 m t^2/2)}{1-\exp( -  (\gamma^2-1) (\tau_1-\gamma\tau_2)^2 m t^2/2)}$. We decompose $K(s,t)$ as
$$K(s,t) = \bigcup_{k=1}^\infty \left\{\Delta \in K(s,t) \ \big| \  \gamma^{k-1} t \leq \langle \mathcal{Q}_\beta(\Delta),\Delta\rangle \leq \gamma^{k} t \right\}.$$
For any $a\geq t$, we further define
$K(s,t,a):=\{\Delta \in K(s,t) \mid  \langle \mathcal{Q}_\beta(\Delta),\Delta\rangle \leq a\}.$
Then we get $E \subseteq \bigcup_{k=1}^\infty E_k$ with
\begin{align*}
E_k=\Big\{\exists \, \Delta \!\in \! K(s,t,\gamma^{k}t) \,\text{such that}\,  \Big|\frac{1}{m}\|\mathcal{R}_\Omega(\Delta)\|_2^2-\langle \mathcal{Q}_\beta(\Delta),\Delta\rangle \Big| \geq    \gamma^{k-1}\tau_1t + \frac{16}{\tau_2} s^2 \mu_1  d_2  \vartheta_m^2  \Big\}.
\end{align*}
Now we need to estimate the probability of each event $E_k$.  Define
$$Z_a:= \sup_{\Delta \in K(s,t,a)}  \Big| \frac{1}{m} \|\mathcal{R}_\Omega(\Delta)\|_2^2-\langle \mathcal{Q}_\beta(\Delta),\Delta \rangle  \Big|.$$ Notice that for any $\Delta \in \mathbb{V}^{n_1\times n_2}$,
$$\frac{1}{m} \|\mathcal{R}_\Omega(\Delta)\|_2^2 = \frac{1}{m}\sum_{i=1}^m \langle \Gamma_{\omega_i}, \Delta\rangle^2 \ \stackrel{a.s.}{\rightarrow} \ \mathbb{E}(\langle \Gamma_{\omega_i}, \Delta\rangle^2) = \langle \mathcal{Q}_\beta(\Delta),\Delta\rangle.$$
 Since $\|\mathcal{R}_\beta(\Delta)\|_\infty \leq 1$ for all $\Delta \in K(s,t)$, from Massart's Hoeffding type concentration inequality \cite[Theorem 1.4]{Mas98} for suprema of empirical processes, we have
  \begin{equation}\label{eqnmasineq}
 {\rm Pr}\big( Z_a \geq \mathbb{E}(Z_a) + \varepsilon \big) \leq \exp\bigg(\!\!-\!\frac{m\varepsilon^2}{2}\bigg) \quad \forall \, \varepsilon > 0.
 \end{equation}
 Next, we use the standard Rademacher symmetrization in the theory of empirical processes to further derive an upper bound of $\mathbb{E}(Z_a)$.  Let $\{\epsilon_1,\ldots, \epsilon_m\}$ be a Rademacher sequence. Then, we have
\begin{align}
\mathbb{E}(Z_a) = & \ \mathbb{E} \bigg(\sup_{\Delta \in K(s,t,a)}  \Big| \frac{1}{m} \sum_{i=1}^m  \langle \Gamma_{\omega_i}, \Delta \rangle^2 -\mathbb{E}\big(\langle \Gamma_{\omega_i}, \Delta \rangle^2\big)    \Big|\bigg) \nonumber\\
\leq & \ 2 \mathbb{E} \bigg(\sup_{\Delta \in K(s,t,a)} \Big|\frac{1}{m}\sum_{i=1}^m \epsilon_i \langle \Gamma_{\omega_i}, \Delta\rangle^2 \Big|\bigg)
\leq 8 \mathbb{E} \bigg(\sup_{\Delta \in K(s,t,a)} \Big|\frac{1}{m}\sum_{i=1}^m \epsilon_i \langle \Gamma_{\omega_i}, \Delta\rangle \Big|\bigg) \nonumber \\
= & \ 8 \mathbb{E} \bigg(\sup_{\Delta \in K(s,t,a)} \Big|\frac{1}{m}\sum_{i=1}^m  \langle \mathcal{R}_\Omega^*(\epsilon), \Delta\rangle  \Big|\bigg)
 \leq 8 \mathbb{E} \, \Big\|\frac{1}{m}\mathcal{R}_\Omega^*(\epsilon)\Big\|\bigg( \sup_{\Delta \in K(s,t,a)} \|\Delta\|_*\bigg), \label{eqnmasineq1}
\end{align}
 where the first inequality follows from the symmetrization theorem (e.g., see \cite[Lemma 2.3.1]{VanW96} and \cite[Theorem 14.3]{BuhV11})  and the second inequality follows from the contraction theorem (e.g., see \cite[Theorem 4.12]{LedT91} and \cite[Theorem 14.4]{BuhV11}).
Moreover, from (\ref{eqndefL}), we have
\begin{equation}\label{eqndefiot}
\langle \mathcal{Q}_\beta(\Delta), \Delta \rangle \geq (\mu_1 d_2)^{-1}\|\Delta\|_F^2 \quad \forall \, \Delta \in \{\Delta \! \in \! \mathbb{V}^{n_1\times n_2} \mid \mathcal{R}_\alpha(\Delta) = 0\}.
\end{equation}
This leads to
\begin{equation}\label{eqnmasineq2}
 \|\Delta\|_*\leq s \|\Delta\|_F \leq s \sqrt{\mu_1 d_2 \langle \mathcal{Q}_\beta(\Delta), \Delta\rangle} \leq s \sqrt{\mu_1  d_2 a} \quad \forall \, \Delta \in K(s,t,a).
\end{equation}
Combining (\ref{eqnmasineq1}) and (\ref{eqnmasineq2}) with the definition of $\vartheta_m$ in (\ref{defvaritheta}), we obtain that
 \begin{equation*}
 \begin{aligned}
  \mathbb{E}(Z_a) + \Big(\frac{\tau_1}{\gamma}-\tau_2\Big) a  \leq & \ 8 \vartheta_m s\sqrt{\mu_1  d_2 a}  + \Big(\frac{\tau_1}{\gamma}-\tau_2\Big) a
%=  \   s \vartheta_m \sqrt{\frac{32 \mu_1 d_2}{\tau_2}}\cdot
%\sqrt{2a\tau_2} + \Big(\frac{\tau_1}{\gamma}-\tau_2\Big) a \\
\leq  \frac{16}{\tau_2}  s^2 \mu_1 d_2  \vartheta^2_m + \frac{\tau_1}{\gamma}a,
 \end{aligned}
 \end{equation*}
 where the second inequality follows from the simple fact   $x_1 x_2 \leq (x_1^2+x_2^2)/2$ for any $x_1, x_2\geq 0$.
 Then,  it follows from (\ref{eqnmasineq}) that
\begin{align*}
{\rm Pr} \bigg(\! Z_a \geq \frac{\tau_1}{\gamma}a +  \frac{16}{\tau_2} s^2 \mu_1 d_2  \vartheta_m^2 \!\bigg) \leq  {\rm Pr} \bigg(\!Z_a \geq \mathbb{E}(Z_a)+ \Big(\frac{\tau_1}{\gamma}\!-\!\tau_2\Big)a \!\bigg)
\leq \exp\bigg(\!\!-\!{\Big(\frac{\tau_1}{\gamma}\!-\!\tau_2\Big)}^2\,\frac{ma^2}{2}\! \bigg).
\end{align*}
 This implies that
 ${\rm Pr}(E_k) \leq \exp\Big(\!\!-\!\frac{1}{2} \, \gamma^{2(k-1)}(\tau_1-\gamma\tau_2)^2 m t^2\Big).$
 Then, since $\gamma > 1$, by using $\gamma^k \geq 1 + k(\gamma -1)$ for any $k\geq1$, we have
\begin{align*}
{\rm Pr}(E) & \ \leq  \sum_{k=1}^\infty {\rm Pr}(E_k) \leq \sum_{k=1}^\infty \exp\Big(\!\!-\!\frac{1}{2} \, \gamma^{2(k-1)}(\tau_1-\gamma\tau_2)^2 m t^2\Big) \\
%& \ \leq \exp\Big(\!\!-\!\frac{1}{2}(\tau_1-\gamma\tau_2)^2 m t^2\Big) \sum_{k=1}^\infty \exp\Big(\!\!-\!\frac{1}{2} (\gamma^{2(k-1)}-1)(\tau_1-\gamma\tau_2)^2 m t^2\Big) \\
& \ \leq  \exp\Big(\!\!-\!\frac{1}{2}(\tau_1-\gamma\tau_2)^2
m t^2\Big) \sum_{k=1}^\infty \exp\Big(\!\!-\!\frac{1}{2} (k-1)(\gamma^2-1)(\tau_1-\gamma\tau_2)^2 m t^2\Big) \\
& \ \leq \frac{\exp(- (\tau_1-\gamma\tau_2)^2 m t^2/2)}{1-\exp(- (\gamma^2-1) (\tau_1-\gamma\tau_2)^2 m t^2/2)}.
\end{align*}
Thus, we complete the proof of Lemma \ref{lemiso}.
\end{proof}

\bigskip

Now we proceed with the proof of Theorem \ref{thmstobd}.
Let $\Delta_m :=\widehat{X}_m-\overline{X}$. Notice that the equality (\ref{eqnopbdtig}) implies that
$$\|\mathcal{P}_{T^\perp}(\Delta_m)\|_* \leq  \frac{1}{\kappa-1} \|\mathcal{P}_T(\Delta_m)\|_* + \frac{\kappa}{\kappa-1} a_m \sqrt{r} \|\Delta_m \|_F.$$
This, together with $\|\mathcal{P}_T(\Delta_m)\|_* \leq \sqrt{2r} \|\Delta_m\|_F$, leads to
\begin{equation} \label{eqnbdmid1}
\|\Delta_m\|_* \leq \|\mathcal{P}_T(\Delta_m)\|_*  + \|\mathcal{P}_{T^\perp}(\Delta_m)\|_* \leq \frac{\kappa}{\kappa-1}\big(\sqrt{2} + a_m\big)\sqrt{r} \|\Delta_m \|_F.
 \end{equation}
Let $b_m := \|\mathcal{R}_\beta(\Delta_m)\|_\infty \leq 2b$. For any fixed $c>0$, $\gamma>1$, $\tau_1 \in (0,1)$ and $\tau_2 \in (0, \tau_2/\gamma)$, define $t_m:= \sqrt{\frac{2c\log(n_1+n_2)}{(\tau_1-\gamma\tau_2)^2m}}$ so that direct calculation yields
\begin{align*}
 \frac{\exp (\!-(\tau_1 - \gamma\tau_2)^2 m t_m^2/2)}{1 - \exp (\!- (\gamma^2 - 1) (\tau_1 - \gamma\tau_2)^2 m t_m^2/2)} =  \frac{(n_1 + n_2)^{-c}} {1 - (n_1 + n_2)^{-(\gamma^2-1)c}} \leq  \frac{(n_1 + n_2)^{-c}}{1 - 2^{-(\gamma^2 - 1)c}}.
\end{align*}
Then we separate the discussion into two cases:

\noindent
{Case 1:} $\langle \mathcal{Q}_\beta(\Delta_m), \Delta_m \rangle \leq b_m^2 t_m$. It follows from (\ref{eqndefiot}) that $\|\Delta_m\|_F^2/d_2\leq 4 b^2 \mu_1  t_m$.

\noindent
{Case 2:} $\langle \mathcal{Q}_\beta(\Delta_m), \Delta_m \rangle> b_m^2 t_m$. It follows from (\ref{eqnbdmid1}) that $\Delta_m/b_m \in K(s_m,t_m)$ with $s_m := \frac{\kappa}{\kappa-1}\big(\sqrt{2}+a_m\big)\sqrt{r}$. Then for any given $\tau_3$ satisfying $0 < \tau_3 <1$, we obtain that with probability at least $1-\frac{(n_1 + n_2)^{-c}}{1 - 2^{-(\gamma^2 - 1)c}}$,
\begin{align*}
\frac{\|\Delta_m\|_F^2}{d_2} \leq & \  \mu_1 \langle \mathcal{Q}_\beta(\Delta_m), \Delta_m\rangle \leq \frac{\mu_1}{1-\tau_1}\bigg(\frac{1}{m}\|\mathcal{R}_\Omega(\Delta_m)\|_2^2+  \frac{16}{\tau_2} s_m^2 \mu_1  d_2 \vartheta_m^2 b_m^2\bigg)  \\
\leq & \ \frac{2}{1-\tau_1}\bigg(\frac{\sqrt{2}}{\kappa} + a_m\bigg)\mu_1\rho_m\sqrt{r}\|\Delta_m\|_F +\frac{16}{(1-\tau_1)\tau_2} s_m^2 \mu_1^2 d_2 \vartheta_m^2 b_m^2\\
\leq & \  \tau_3 \frac{\|\Delta_m\|_F^2}{d_2} + \frac{2}{(1-\tau_1)^2\tau_3}\bigg(\frac{\sqrt{2}}{\kappa} + a_m\bigg)^2\mu_1^2\rho_m^2r d_2 + \frac{16}{(1-\tau_1)\tau_2} s_m^2 \mu_1^2 d_2 \vartheta_m^2 b_m^2,
\end{align*}
where the first inequality follows from (\ref{eqndefiot}), the second inequality follows from Lemma \ref{lemiso} and the third inequality follows from Theorem \ref{thmopbd}.
Plugging in $s_m$ further leads to
$$\frac{\|\Delta_m\|_F^2}{d_2} \leq \frac{\mu_1^2 d_2 r}{1-\tau_3}\left(\frac{2}{(1\!-\!\tau_1)^2\tau_3}\bigg(
\frac{\sqrt{2}}{\kappa} + a_m\bigg)^2 \rho_m^2  + \frac{64}{(1-\tau_1)\tau_2} \bigg(\frac{\kappa}{\kappa\!-\!1}\bigg)^2\big(\sqrt{2}+a_m\big)^2 \vartheta_m^2 b^2\right).$$

\noindent
Combing the above two cases together, with $\gamma$, $\tau_1$, $\tau_2$ and $\tau_3$ chosen to be absolute constants, we arrive at an intermediate result that
there exist some positive absolute constants $c'_0,c'_1, c'_2$ and $C'_0$ such that for any $\kappa > 1$, if $\rho_m$ is chosen as in Theorem \ref{thmopbd}, then with probability at least $1-c'_1 (n_1+n_2)^{-c'_2}$,
\begin{align}
\frac{\|\widehat{X}_m - \overline{X}\|_F^2}{d_2} \leq & \ C'_0\max \Bigg\{ \mu_1^2 d_2 r \left({c'_0}^2\bigg(\frac{\sqrt{2}}{\kappa} + a_m \bigg)^2 \rho_m^2 +  \bigg( \frac{\kappa}{\kappa - 1} \bigg)^2\big( \sqrt{2} + a_m \big)^2 \vartheta_m ^2 b^2 \right) , \nonumber \\
& \ \hspace{7.0cm} b^2  \mu_1 \sqrt{ \frac{\log(n_1 + n_2)}{m} } \Bigg\}. \label{eqnthmbdmid}
\end{align}
%\begin{align*}
%\frac{\|\widehat{X}_m \!- \!\overline{X}\|_F^2}{d_2}\! \leq \! C_0 \!\max\!\Bigg\{\!\mu_1^2 d_2 r \!\left(\!c_0^2\Big(\frac{\sqrt{2}}{\kappa} \!+ \! a_m\Big)^2\!\rho_m^2\! +\! \Big(\!\frac{\kappa}{\kappa\!-\!1}\! \Big)^{\!2}\!\Big(\!\sqrt{2}\!+\!a_m\!\Big)^{\!2} \vartheta_m ^2 b^2 \!\right)\!,
%b^2 \mu_1\!\sqrt{\frac{\log(n_1\!+\!n_2)}{m}}\!\Bigg\}.
%\end{align*}

To further derive explicit estimations of $\rho_m$ and $\vartheta_m$, we introduce the noncommutative Bernstein inequality taken from \cite[Corollary 2.1]{Kol11}, which provides a probability control of  the deviation of  the sum of random matrices from its mean in the operator norm. The noncommutative Bernstein inequality introduced here is a recently-extended version, with the random matrices being controlled by the Orlicz norms (see \cite{Kol11, Kol12, KolLT11}) rather than the operator norm (see, e.g., \cite{Rec11,Tro11,Gro11}). The Orlicz norms are used to characterize the tail behavior of random variables. Given any $s\geq 1$, the $\psi_s$ Orlicz norm of a random variable $z$ is defined by
$\|z\|_{\psi_s}:=\inf\{t>0 \,\big| \, \mathbb{E}\exp(|z|^s/t^s) \leq 2\}$.

\begin{lemma}[Koltchinskii \cite{Kol11}]\label{lemnbi}
Let $Z_1,\ldots, Z_m\in\mathbb{V}^{n_1\times n_2}$ be independent random matrices with mean zero.
Suppose that
$\max\big\{\big\|\|Z_i\|\big\|_{\psi_s},2\mathbb{E}^{\frac{1}{2}}(\|Z_i\|^2)\big\} <\varpi_{s}$
for some constant $\varpi_{s}$.
Define $$\sigma_Z := \max\left\{\bigg\|\frac{1}{m} \sum_{i=1}^m \mathbb{E}(Z_iZ_i^\mathbb{T})\bigg\|^{1/2},\ \bigg\|\frac{1}{m} \sum_{i=1}^m \mathbb{E}(Z_i^\mathbb{T}Z_i)\bigg\|^{1/2}\right\}.$$
Then, there exists a constant $C$ such that for all $t>0$, with probability at least $1\!-\exp(-t)$,
$$\bigg\|\frac{1}{m} \sum_{i=1}^m Z_i\bigg\| \leq C \max\left\{ \sigma_Z\sqrt{\frac{t+\log(n_1+n_2)}{m}}, \varpi_{s}\left(\log\frac{\varpi_{s}}{\sigma_Z}\right)^{1/s}\frac{t+\log(n_1\!+n_2)}{m}\right\}.$$
\end{lemma}

With the help of Lemma \ref{lemnbi}, we obtain the following result, which is an extension of  \cite[Lemma 2]{KolLT11} and  \cite[Lemmas 5 \& 6]{Klo12} from the standard basis to an arbitrary orthonormal basis. A similar result can also be found in \cite[Lemma 6]{NegW12}.
\begin{lemma}\label{lemben}
Under Assumption \ref{asmpnoi}, there exists a positive constant $C'$ (only depending on the $\psi_1$ Orlicz norm of $\xi_k$) such that for all $t>0$, with probability at least $1-\exp(-t)$,
\begin{equation}\label{eqnnoibd}
\left\|\frac{1}{m}\mathcal{R}^*_\Omega(\xi)\right\|  \leq  C'\max\left\{\sqrt{\frac{\mu_2(t\!+ \log(n_1\!+ n_2))}{\sqrt{d_2}m}}, \frac{\log(d_2)(t\!+ \log(n_1\!+ n_2))}{2m}\right\}.
\end{equation}
In particular, when $m\geq \sqrt{d_2}\log^3(n_1+n_2)/\mu_2$, we also have
\begin{equation}\label{eqnnoiexpd}
\mathbb{E}\,\bigg\|\frac{1}{m}\mathcal{R}^*_\Omega(\xi)\bigg\| \leq C'\sqrt{\frac{2e\mu_2\log(n_1+n_2)}{\sqrt{d_2}m}},
\end{equation}
where $e$ is the exponential constant.
\end{lemma}

\begin{proof}
Recall that $\frac{1}{m} \mathcal{R}^*_\Omega(\xi) = \frac{1}{m}\sum_{i=1}^m \xi_i \Theta_{\omega_i}$. Let $Z_i := \xi_i \Theta_{\omega_i}$. Since $\mathbb{E}(\xi_i)=0$, the independence of $\xi_i$ and $\Theta_{\omega_i}$ implies that $\mathbb{E}(Z_i)=0$. Since $\|\Theta_{\omega_i}\|_F=1$, we have that $\|Z_i\|\leq \|Z_i\|_F = |\xi_i|\|\Theta_{\omega_i}\|_F = |\xi_i|$.
It follows that $\big\|\|Z_i\| \big\|_{\psi_1} \leq \|\xi_i\|_{\psi_1}$ and thus finite. (It is known that a random variable is sub-exponential if and only its $\psi_1$ Orlicz norm is finite \cite{VanW96}).  Meanwhile, $\mathbb{E}^\frac{1}{2}(\|Z_i\|^2)\leq \mathbb{E}^\frac{1}{2}(\|Z_i\|_F^2) = \mathbb{E}^\frac{1}{2}(\xi_i^2)=1$.
Then direct calculation yields
$$\mathbb{E}\big(Z_iZ_i^\mathbb{T}\big)=\mathbb{E}\big(\xi^2_i \Theta_{\omega_i}\Theta_{\omega_i}^\mathbb{T}\big)=\mathbb{E}\big( \Theta_{\omega_i}\Theta_{\omega_i}^\mathbb{T}\big) =\sum_{k\in\beta}p_k \Theta_k \Theta_k^\mathbb{T}.$$
The calculation for $\mathbb{E}\big(Z_i^\mathbb{T} Z_i\big)$ is similar. We obtain from (\ref{eqndefL}) that $1/\sqrt{d_2}\leq \sigma_Z^2 \leq \mu_2/\sqrt{d_2}$. Then, applying this to Lemma \ref{lemnbi} yields (\ref{eqnnoibd}). The remaining proof of (\ref{eqnnoiexpd}) follows the same as the proof of Lemma 6 in \cite{Klo12}. For simplicity, we omit it.
\end{proof}

\bigskip

A good estimation of $\rho_m$ can be achieved by choosing $t=c'_2\log(n_1 +n_2)$ in Lemma \ref{lemben} for an optimal order bound, where $c'_2$ is the same as that in (\ref{eqnthmbdmid}). With this choice, when $m \geq 4(1+c'_2)\sqrt{d_2}\log^2(d_2)\log(n_1+n_2)/\mu_2$,
the first term in the maximum of (\ref{eqnnoibd}) dominates the second one. Thus, with probability at least $1- (n_1+n_2)^{-c'_2}$, one can choose
$$\rho_m =  \kappa \nu \cdot C' \sqrt{\frac{(1+c'_2)\mu_2\log(n_1+n_2)}{\sqrt{d_2}m}}.$$
Moreover, since Bernoulli random variables are sub-exponential, Lemma \ref{lemben} also provides an upper bound of $\vartheta_m$ in (\ref{eqnnoiexpd}). It is worthwhile to note that after plugging the above estimations of $\rho_m$ and $\vartheta_m$, the second term in the maximum of (\ref{eqnthmbdmid}) is negligible compared with the first term. Therefore, the second term is further dropped for simplicity and thus we complete the proof.

\subsection{Proof of Theorem \ref{thmambound}}
For notational simplicity, we drop the subscript of $\widetilde{X}_m$ in this proof.
With $(\widetilde{U}, \widetilde{V}) \in \mathbb{O}^{n_1,n_2}(\widetilde{X})$, one immediately obtains from the definition of $a_m$ in (\ref{eqndelalpbet}) that
\begin{equation}\label{eqnambndinproof}
a_m \leq \frac{1}{\sqrt{r}} \big(\|F(\widetilde{X}) - \widetilde{U}_1\widetilde{V}_1^{\mathbb{T}}\|_F + \|\widetilde{U}_1\widetilde{V}_1^{\mathbb{T}} - \overline{U}_1\overline{V}_1^{\mathbb{T}}\|_F\big) \leq \varepsilon_F(\widetilde{X}) + \frac{1}{\sqrt{r}} \|\widetilde{U}_1\widetilde{V}_1^{\mathbb{T}} - \overline{U}_1\overline{V}_1^{\mathbb{T}}\|_F.
\end{equation}
The left proof is to find an upper bound of $\|\widetilde{U}_1\widetilde{V}_1^{\mathbb{T}} - \overline{U}_1 \overline{V}_1^{\mathbb{T}}\|_F $. Let $\delta := \|\widetilde{X} -\overline{X}\|_F$ and  $\mathcal{N}_\delta(\overline{X}) := \{X\in \mathbb{V}^{n_1\times n_2} \mid \|X - \overline{X}\|_F \leq \delta\}$.

Let $\widehat{F}:\mathbb{V}^{n_1\times n_2} \rightarrow \mathbb{V}^{n_1\times n_2}$ be a spectral operator associated with a symmetric function $\widehat{f}:\mathbb{R}^n \rightarrow \mathbb{R}^n$ given by $\widehat{f}_i(x) = \phi(x_i)$, $i=1,\ldots, n$, where $\phi:\mathbb{R} \rightarrow \mathbb{R}$ is an odd scalar function with $\phi(t) = -\phi(-t)$ for $t < 0$, and $\phi(t)$ for $t\geq 0$  is defined as
$$\phi(t) =  \begin{cases} 1 & \text{if} \ \ t \geq 2\sigma_r(\overline{X})/3 - \delta/3, \\
\frac{t-(\sigma_r(\overline{X})/3 + \delta/3)}{\sigma_r(\overline{X})/3 - 2\sigma/3} & \text{if}  \ \  \sigma_r(\overline{X})/3 + \delta/3 < t  < 2\sigma_r(\overline{X})/3 - \delta/3, \\ 0 & \text{if} \ \  0 \leq  t \leq \sigma_r(\overline{X})/3 + \delta/3. \end{cases}$$
Note that for any $X \in \mathcal{N}_\delta(\overline{X})$,
$$|\sigma_i(X) - \sigma_i(\overline{X})| \leq \sigma_1(X-\overline{X})\leq \|X - \overline{X}\|_F \leq \delta, \quad i = 1,\ldots, n.$$
Since $\delta/\sigma_r(\overline{X}) < 1/2$, we further have  $\sigma_r(X) \geq \sigma_r(\overline{X}) - \delta  > \delta \geq \sigma_{r+1}(X)$.
This means $$\widehat{F}(X) = U_1V_1^{\mathbb{T}} \qquad \forall\, X \in \mathcal{N}_\delta(\overline{X}).$$ Moreover, $\widehat{F}$ is continuously differentiable over $\mathcal{N}_\delta(\overline{X})$. Hence, we can apply the Mean Value Theorem to obtain
\begin{equation}\label{eqnMVT}
 \widetilde{U}_1\widetilde{V}_1^{\mathbb{T}} - \overline{U}_1\overline{V}_1^{\mathbb{T}} = \widehat{F}(\widetilde{X}) - \widehat{F}(\overline{X}) = \int_0^1 \widehat{F}'(\widetilde{X}_t)(\widetilde{X}-\overline{X}) \, {\rm d}t,
\end{equation}
where $\widetilde{X}_t:= \overline{X} + t(\widetilde{X} - \overline{X})$. Clearly, $\widetilde{X}_t \in \mathcal{N}_\delta(\overline{X})$ when $t\in [0,1]$.

Regarding (\ref{eqnMVT}), we need to look into the derivative of $\widehat{F}$ over $\mathcal{N}_\delta(\overline{X})$. Let $X \in \mathcal{N}_\delta(\overline{X})$ be arbitrary and $(U, V) \in \mathbb{O}^{n_1, n_2}(X)$. Without loss of generality, we assume $n_1 \leq n_2$. Let $\chi_1:=\{1, \ldots, r\}$, $\chi_2:=\{r+1, \ldots, n_1\}$ and $\chi_3:= \{n_1+1,\ldots, n_2\}$.  Then, according to \cite[Theorem 3.6]{Din12}, we have that for any $H \in \mathbb{V}^{n_1\times n_2}$,
\begin{equation}\label{eqnspetralderiv}
\widehat{F}'(X)(H) = U\bigg[\mathcal{E}_1(X)\circ \frac{\widetilde{H}_1 + \widetilde{H}_1^{\mathbb{T}}}{2} + \mathcal{E}_2(X) \circ \frac{\widetilde{H}_1 - \widetilde{H}_1^{\mathbb{T}}}{2} \ \ \ \ \Upsilon(X)\circ \widetilde{H}_2 \bigg]V^{\mathbb{T}},
\end{equation}
where $[\widetilde{H}_1 \ \widetilde{H}_2] = \widetilde{H} := UHV^{\mathbb{T}}$ with $\widetilde{H}_1 \in \mathbb{V}^{n_1\times n_1}$, $\widetilde{H}_2 \in \mathbb{V}^{n_1 \times (n_2-n_1)}$,   and $\mathcal{E}_1(X) \in \mathbb{V}^{n_1\times n_1}$, $\mathcal{E}_2(X)\in \mathbb{V}^{n_1\times n_1}$, $\Upsilon(X) \in \mathbb{V}^{n_1\times (n_2-n_1)}$ take the form
\begin{align*}
& \big(\mathcal{E}_1(X)\big)_{ij} = \begin{cases}
 \frac{1}{\sigma_i(X) - \sigma_j(X)} & \text{if} \ \ i \in \chi_1, j\in \chi_2 \ \text{or} \ i \in \chi_2, j \in \chi_1, \\
 0 & \text{otherwise},
\end{cases} \\
& \big(\mathcal{E}_2(X)\big)_{ij} = \begin{cases}
 \frac{2}{\sigma_i(X)+\sigma_j(X)} & \text{if} \  \ i\in \chi_1, j\in \chi_1, \\
 \frac{1}{\sigma_i(X) - \sigma_j(X)} & \text{if} \ \ i \in \chi_1, j\in \chi_2 \ \text{or} \ i \in \chi_2, j \in \chi_1,\\
 0 & \text{otherwise}, \\
\end{cases}\\
& \big(\Upsilon(X)\big)_{ij} = \, \begin{cases}
 \frac{1}{\sigma_i(X)}  & \hspace{1.05cm} \text{if} \ \ i\in \chi_1, j\in \chi_3, \\
 0 & \hspace{1.05cm} \text{otherwise}.
\end{cases}
\end{align*}
Here, ``$\circ$'' stands for the Hadamard product of matrices. Let $\Delta$ denote the matrix in the bracket of (\ref{eqnspetralderiv}). Moreover, let $\Delta_{\chi_i, \chi_j}$ and $\widetilde{H}_{\chi_i, \chi_j}$ denote the submatrices of $\Delta$  and $\widetilde{H}$  with row indices $\chi_i$ and column indices $\chi_j$, respectively. Then, a direct calculation yields
\begin{gather*}
 \|\Delta_{\chi_1, \chi_1} \|_F^2 \leq \frac{\|\widetilde{H}_{\chi_1, \chi_1}\|_F^2}{\sigma_r^2(X)}, \qquad
  \|\Delta_{\chi_1, \chi_2}  \|_F^2 + \|\Delta_{\chi_2, \chi_1}  \|_F^2  \leq \frac{\|\widetilde{H}_{\chi_1, \chi_2}\|_F^2 + \|\widetilde{H}_{\chi_2, \chi_1}\|_F^2}{(\sigma_r(X) - \sigma_{r+1}(X))^2}, \\
\|\Delta_{\chi_2, \chi_2} \|_F^2 =0, \quad \|\Delta_{\chi_1, \chi_3} \|_F^2 \leq \frac{\|\widetilde{H}_{\chi_1, \chi_3}\|_F^2}{\sigma_r^2(X)}  \quad \text{and} \quad  \|\Delta_{\chi_2, \chi_3} \|_F^2 = 0.
\end{gather*}
%\begin{gather*}
%\sum_{i\in \chi_1} \sum_{j \in \chi_1} \Delta_{ij}^2 \leq \frac{1}{\sigma_r^2(X)}\sum_{\scriptstyle i\in \chi_1 \atop \scriptstyle i\neq j} \sum_{j\in \chi} \widetilde{H}_{ij}^2, \qquad \sum_{i\in \chi_2}\sum_{j\in \chi_3} \Delta_{ij}^2 = 0, \\
%\sum_{i\in \chi_1}\sum_{j\in \chi_2} \Delta_{ij}^2 + \sum_{i \in \chi_2}\sum_{j \in \chi_1} \Delta_{ij}^2 \leq \frac{1}{(\sigma_r(X)-\sigma_{r+1}(X))^2} \Bigg(\sum_{i\in \chi_1}\sum_{j\in \chi_2} \widetilde{H}_{ij}^2 + \sum_{i\in \chi_2}\sum_{j\in \chi_1} \widetilde{H}_{ij}^2\Bigg), \\
%\sum_{i\in \chi_3} \sum_{j \in \chi_1} \Delta_{ij}^2 \leq \frac{1}{\sigma_r^2(X)}\sum_{\scriptstyle i\in \chi_1} \sum_{j\in \chi} \widetilde{H}_{ij}^2,
%\sum_{i\in \chi_2}\sum_{j\in \chi_2} \Delta_{ij}^2 = 0.
%\end{gather*}
Note that $\|\widehat{F}'(X)(H)\|_F  = \|\Delta\|_F$ and $\|\widetilde{H}\|_F = \|H\|_F$. By summing up the above inequalities together, we obtain that for any $X\in \mathcal{N}_\delta(\overline{X})$,
\begin{equation}\label{eqndevbnd}
\big\|\widehat{F}'(X)(H)\big\|_F  \leq \frac{\sqrt{\|H_{\chi_1, \chi_1\cup\chi_2\cup\chi_3}\|_F^2+\|H_{\chi_2, \chi_1}\|_F^2}}{\sigma_r(X)-\sigma_{r+1}(X)} \leq \frac{\|H\|_F}{\sigma_r(X)-\sigma_{r+1}(X)}.
\end{equation}
Now, we proceed with the proof by applying (\ref{eqndevbnd}) to (\ref{eqnMVT}). This leads to
\begin{equation}\label{eqnU1V1bound}
\big\|\widetilde{U}_1\widetilde{V}_1^{\mathbb{T}} - \overline{U}_1\overline{V}_1^{\mathbb{T}}\big\|_F \leq \int_0^1 \big\|\widehat{F}'(\widetilde{X}_t)(\widetilde{X}-\overline{X})\big\|_F \, {\rm d} t \leq \int_0^1 \frac{\delta}{\sigma_r(\widetilde{X}_t) - \sigma_{r+1}(\widetilde{X}_t)} \, {\rm d} t.
\end{equation}
Moreover, using \cite[Theorems IV.3.4 \& II.3.1]{Bha97}, we have
$$\big(\sigma_r(\widetilde{X}_t) - \sigma_r(\overline{X})\big)^2 + \sigma_{r+1}^2(\widetilde{X}_t) \leq \big\|\sigma(\widetilde{X}_t) - \sigma(\overline{X})\big\|_F^2\leq \big\|\sigma(\widetilde{X}_t-\overline{X})\big\|_F^2 = \big\|\widetilde{X}_t - \overline{X}\big\|_F^2 \leq t^2\delta^2.$$
This implies that $\sigma_r(\widetilde{X}_t)-\sigma_r(\overline{X}) = \delta_t \cos\theta$ and $\sigma_{r+1}(\widetilde{X}_t) = \delta_t \sin\theta$ for some $\delta_t \leq t\delta$ and  $\theta \in [0, 2\pi)$. Thus,
\begin{equation}\label{eqnsigmadifbnd}
\sigma_r(\widetilde{X}_t) - \sigma_{r+1}(\widetilde{X}_t) = \sigma_r(\overline{X}) + \delta_t\cos\theta - \delta_t\sin\theta \geq \sigma_r(\overline{X}) - \sqrt{2}\delta_t \geq \sigma_r(\overline{X}) - \sqrt{2}t\delta.
\end{equation}
Substituting (\ref{eqnsigmadifbnd}) into (\ref{eqnU1V1bound}), we obtain that
$$ \big\|\widetilde{U}_1\widetilde{V}_1^{\mathbb{T}} - \overline{U}_1\overline{V}_1^{\mathbb{T}}\|_F  \leq \int_0^1 \frac{\delta}{\sigma_r(\overline{X}) - \sqrt{2}t\delta} \, {\rm d} t = -\frac{1}{\sqrt{2}}\log\bigg(1 - \frac{\sqrt{2} \, \delta}{\sigma_r(\overline{X})}\bigg).$$
This, together with (\ref{eqnambndinproof}), completes the proof.

\subsection{Proof of Theorem \ref{thmcons}}

We first prove the following properties of the sample operator $\mathcal{R}_\Omega$ and its adjoint $\mathcal{R}_\Omega^*$.

\begin{lemma}\label{lemoper}
{\bf (i)} For any given $X \in \mathbb{V}^{n_1\times n_2}$, the random matrix $\displaystyle{\frac{1}{m}} \mathcal{R}_\Omega^* \mathcal{R}_\Omega(X) \stackrel{a.s.}{\rightarrow} \mathcal{Q}_\beta(X)$. \\
\noindent
{\bf (ii)} The random vector $\displaystyle{\frac{1}{\sqrt{m}}}\mathcal{R}_{\alpha\cup\beta}\mathcal{R}_\Omega^*(\xi) \stackrel{d}{\rightarrow} N\big(0, {\rm Diag}(p)\big)$, where $p=(p_1,\ldots,p_d)^{\mathbb{T}}$.
\end{lemma}

\begin{proof}
(i) It follows from the definitions of $\mathcal{R}_\Omega$ and its adjoint $\mathcal{R}_\Omega^*$ that
$\frac{1}{m} \mathcal{R}_\Omega^* \mathcal{R}_\Omega(X) = \frac{1}{m} \sum_{i=1}^m \langle \Theta_{\omega_i}, X \rangle \, \Theta_{\omega_i}.$
This is an average value of $m$ i.i.d. random matrices $\langle \Theta_{\omega_i}, X \rangle \Theta_{\omega_i}$. Note that  $\mathbb{E}\big(\langle \Theta_{\omega_i}, X \rangle \Theta_{\omega_i}\big) = \mathcal{Q}_\beta(X) \ \forall \, i=1,\cdots, m$.
Then the result follows directly from the strong law of large numbers.

\noindent
(ii) It directly follows from the definitions of $\mathcal{R}_\Omega^*$ and $\mathcal{R}_{\alpha\cup\beta}$ that
$\frac{1}{\sqrt{m}} \mathcal{R}_{\mathcal{\alpha\cup\beta}} \mathcal{R}_\Omega^*(\xi)
= \frac{1}{\sqrt{m}} \mathcal{R}_{\alpha\cup\beta} \big(\sum_{i=1}^m \xi_i \Theta_{\omega_i}\big)
= \frac{1}{\sqrt{m}}\sum_{i=1}^m \xi_i \mathcal{R}_{\alpha\cup\beta}(\Theta_{\omega_i}).$
Since $\mathbb{E}(\xi_i) = 0$ and $\mathbb{E}(\xi_i^2) =1$, according to the independence of $\xi_i$ and $\mathcal{R}_{\alpha\cup\beta}(\Theta_{\omega_i})$, we obtain
$\mathbb{E}\big(\xi_i \mathcal{R}_{\alpha\cup\beta}(\Theta_{\omega_i})\big) = 0$ and ${\rm cov}\big(\xi_i\mathcal{R}_{\alpha\cup\beta}(\Theta_{\omega_i})\big)={\rm Diag}(p).$
Then, applying the vector-valued central limit theorem yields the result.
\end{proof}

To prove the convergence in distribution of minimizers, the following theorem of Knight \cite[Theorem 1]{Kni99} on epi-convergence in distribution is particularly useful in this regard (see also \cite[Proposition 9]{HanP06}).

\begin{lemma}[{Knight \cite{Kni99}}]\label{lemepicon}
Let $\{\Phi_m\}$ be a sequence of random lower-semicontinuous functions that epi-converges
in distribution to $\Phi$. Assume that
\begin{description}
\setlength{\itemsep}{0pt}
\item[(i)] $\widehat{x}_m$ is an $\varepsilon_m$-minimizer of $\Phi_m$, i.e., $\Phi_m(\widehat{x}_m)\leq \inf \Phi_m(x) +\varepsilon_m$, where $\varepsilon_m \stackrel{p}{\rightarrow} 0$;
\item[(ii)] $\widehat{x}_m= O_p(1)$;
\item[(iii)] the function $\Phi$ has a unique minimizer $\overline{x}$.
\end{description}
Then, $\widehat{x}_m \stackrel{d}{\rightarrow} \overline{x}$. In addition, if $\Phi$ is a deterministic function,
then $\widehat{x}_m \stackrel{p}{\rightarrow} \overline{x}$.
\end{lemma}

It is know from \cite{Gey96} that $\widehat{x}_m$ is guaranteed to be $O_p(1)$ when all $\Phi_m$ are convex functions and $\Phi$ has a unique minimizer.
 For more details on epi-convergence in distribution, one may refer to King and Wets \cite{KinW91}, Geyer \cite{Gey94}, Pflug \cite{Pfl92,Pfl95} and Knight \cite{Kni99}. As Lemma \ref{lemepicon} is only applicable to unconstrained optimization problems, constrained optimization problems need to be equivalently converted to unconstrained ones using the indicator function of feasible set. This leads to the issue of epi-convergence in distribution of the sum of two sequences of random functions; see, e.g., Pflug \cite[Lemma 1]{Pfl95}.

Now we proceed with the proof of Theorem \ref{thmcons}. Let $\Phi_m$ denote the objective function of (\ref{eqnrcs}) and $\mathcal{F}$ denote the feasible set. Then, the problem (\ref{eqnrcs}) can be concisely written as
$$\min_{X \in \mathbb{V}^{n_1\times n_2}} \big\{\Phi_m(X) + \delta_\mathcal{F}(X)\big\}.$$
By Assumptions \ref{asmpfun} and \ref{asmpini} and Lemma \ref{lemoper}, we have that the convex function $\Phi_m$ converges pointwise in probability to the convex function
$\Phi$, where $\Phi(X):=\frac{1}{2}\langle X-\overline{X},  \mathcal{Q}_\beta(X-\overline{X})\rangle$ for any $X\in\!\mathbb{V}^{n_1\times n_2}$.
 As a direct extension of Rockafellar \cite[Theorem 10.8]{Roc70}, Andersen and Gill \cite[Theorem II.1]{AndG82} proved that the pointwise convergence in probability implies the convergence in probability (and thus in distribution) with respect to the topology of uniform convergence on compact subset. Then, according to Pflug \cite[Lemma 1]{Pfl95}, we further obtain that $\Phi_m+\delta_\mathcal{F}$ epi-converges in distribution to $\Phi+\delta_\mathcal{F}$. Note that $\overline{X}$ is the unique minimizer of $\Phi(X) + \delta_\mathcal{F}(X)$ since $\Phi(X)$ is strongly convex over the feasible set $\mathcal{F}$. Thus, we complete the proof by applying Lemma \ref{lemepicon} on epi-convergence in distribution.

\subsection{Proof of Theorem \ref{thmnessuf}}
Theorem \ref{thmcons} actually implies that $\widehat{X}_m$ has a higher rank than $\overline{X}$ with probability converging to $1$ if $\rho_m \rightarrow 0$, due to the  straightforward result:

\begin{lemma}\label{lemrkrhs}
If $X_m \stackrel{p}{\rightarrow} \overline{X}$, then $\lim\limits_{m\rightarrow \infty}{\rm Pr}\big({\rm rank}(X_m) \geq {\rm rank}(\overline{X})\big)=1$.
\end{lemma}
\begin{proof}
   It follows from the Lipschitz continuity of singular values that $$\sigma_k(X_m)\stackrel{p}{\rightarrow} \sigma_k(X) \quad \forall\, 1\leq k\leq n.$$ Thus, for any $\varepsilon > 0$, we have
$$\mathbb{P}\big({\rm rank}(X_m) \geq {\rm rank}(\overline{X})\big)\geq \mathbb{P}\big(|\sigma_r(X_m)-\sigma_r(\overline{X})|\leq \varepsilon \sigma_r(\overline{X})\big)\rightarrow 1 \quad \text{as} \quad m \rightarrow \infty.$$
\end{proof}

Now we take a look at the local property for the rank function for the perturbation.

\begin{lemma}\label{lemlocrk}
Let $\overline{\Delta} \in \mathbb{V}^{n_1\times n_2}$ satisfy $\overline{U}_2^\mathbb{T} \overline{\Delta}\,\overline{V}_2 \neq 0$.
Then, for all $\rho\neq 0$ sufficiently small and $\Delta$ sufficiently close to $\overline{\Delta}$,
${\rm rank}(\overline{X}+\rho \Delta) > {\rm rank}(\overline{X})$.
\end{lemma}
\begin{proof}
Let $\sigma'_i(X;\cdot)$ denote the directional derivative function of the $i$-th largest singular
value function $\sigma_i(\cdot)$ at $X$. Let $r:=\text{rank}(\overline{X})$. Note that $\sigma_{r+1}(\overline{X})=0$. Then, according to \cite[Section 5.1]{Lew05} and \cite[Proposition 6]{DinST10},
for any $\Delta \in \mathbb{V}^{n_1\times n_2}$ and $\rho \rightarrow 0$, we have
$$\sigma_{r+1}(\overline{X} + \rho\Delta)-\sigma'_{r+1} (\overline{X};\rho\Delta) = O(\|\rho\Delta\|_F^2),$$
where $\sigma_{r+1}'(\overline{X}; \rho\Delta) = \|\overline{U}_2^\mathbb{T} (\rho\Delta)\overline{V}_2\|$.
Since $\overline{U}_2^\mathbb{T}\overline{\Delta} \, \overline{V}_2 \neq 0$, from the sign-preserving property of limits,
for any $\rho\neq 0$ sufficiently small and $\Delta$ sufficiently close to $\overline{\Delta}$, we have
\begin{align*}
\frac{\sigma_{r+1}(\overline{X}+\rho\Delta)}{|\rho|}
= & \ \|\overline{U}_2^\mathbb{T} \Delta \overline{V}_2\| + O(|\rho|\|\Delta\|_F^2)>0.
%\geq & \ \frac{1}{2} \|\overline{U}_2^\mathbb{T}\overline{\Delta} \, \overline{V}_2\| >0.
 \end{align*}
 This implies that ${\rm rank}(\overline{X}+\rho \Delta) > {\rm rank}(\overline{X})$.
\end{proof}

Define $\widehat{\Delta}_m: = \rho_m^{-1}(\widehat{X}_m-\overline{X})$.
To guarantee the efficiency of the nuclear semi-norm on encouraging a low-rank solution, the parameter $\rho_m$ should not decay too fast. Then, for a slow decay on $\rho_m$, we can establish the following result.

\begin{lemma}\label{lemdellim}
If $\rho_m \rightarrow 0$ and $\sqrt{m}\rho_m\rightarrow \infty$, then $\widehat{\Delta}_m \stackrel{p}{\rightarrow} \widehat{\Delta}$, where $\widehat{\Delta}$ is the unique optimal solution to the following  convex optimization problem
\begin{equation}\label{eqndellim}
\begin{aligned}
\min_{\Delta \in \mathbb{V}^{n_1\times n_2}} &\ {\displaystyle \frac{1}{2}} \langle \mathcal{Q}_\beta(\Delta), \Delta \rangle + \langle \overline{U}_1 \overline{V}_1^\mathbb{T} - F(\overline{X}), \Delta \rangle + \|\overline{U}_2^\mathbb{T} \Delta \overline{V}_2\|_*\\
{\rm s.t.}\ \ \ &\ \mathcal{R}_\alpha(\Delta) = 0, \ \ \mathcal{R}_{\beta^+}(\Delta) \leq 0, \ \ \mathcal{R}_{\beta^-}(\Delta) \geq 0.
\end{aligned}
\end{equation}
\end{lemma}
\begin{proof}
 Take a variable transformation $\Delta:= \rho_m^{-1}(X-\overline{X})$ in the optimization problem (\ref{eqnrcs}). Then one can easily see that $\widehat{\Delta}_m$ is the optimal solution to
\begin{equation}\label{eqndellimappr}
\begin{aligned}
\min_{\Delta \in \mathbb{V}^{n_1\times n_2}}&\ {\displaystyle \frac{1}{2m}}\|\mathcal{R}_\Omega(\Delta)\|_2^2 - \frac{\nu}{m \rho_m}\langle \mathcal{R}_\Omega^*(\xi), \Delta\rangle + \frac{1}{\rho_m}\big(\|\overline{X}+\rho_m\Delta \|_* - \|\overline{X}\|_*\big) - \langle F(\widetilde{X}_m), \Delta\rangle \\
{\rm s.t.} \quad  & \Delta \in \mathcal{F}_m := \rho_m^{-1}(\mathcal{K} -\overline{X}),
\end{aligned}
\end{equation}
where $\mathcal{K} := \big\{X \in \mathbb{S}^n \mid \mathcal{R}_\alpha(X) = \mathcal{R}_\alpha(\overline{X}), \ \|\mathcal{R}_\beta(X)\|_\infty \leq b\big\}$.
Let $\Phi_m$ and $\Phi$ denote the objective functions of (\ref{eqndellimappr}) and (\ref{eqndellim}), respectively. By the definition of directional derivative and \cite[Theorem 1]{Wat92}, we have
 $$\lim_{\rho_m \rightarrow 0} \frac{1}{\rho_m}\big(\|\overline{X}+\rho_m\Delta \|_* - \|\overline{X}\|_*\big) = \langle \overline{U}_1\overline{V}_1^\mathbb{T}, \Delta \rangle + \|\overline{U}_2^\mathbb{T} \Delta \overline{V}_2\|_*.$$
Then, under Assumptions \ref{asmpfun} and \ref{asmpini}, according to Lemma \ref{lemoper}, we obtain that $\Phi_m$ converges pointwise in probability to $\Phi$. Together with the convexity of $\mathcal{K}$, we know that $\mathcal{F}_m$ converges in the sense of Painlev{\'e}-Kuratowski to the tangent cone $\mathcal{T}_{\mathcal{K}}(\overline{X})$ (see \cite{RocW98, BonS00}), taking the form
\begin{equation}\label{eqntangentK}
\mathcal{T}_{\mathcal{K}}(\overline{X}) = \big\{\Delta \in \mathbb{V}^{n_1\times n_2} \mid \mathcal{R}_\alpha(\Delta) = 0, \ \mathcal{R}_{\beta^+}(\Delta) \leq 0, \ \mathcal{R}_{\beta^-}(\Delta) \geq 0\big\}.
\end{equation}
Since epi-convergence of functions corresponds to set convergence of their epigraphs \cite{RocW98}, we obtain that $\delta_{\mathcal{F}_m}$ epi-converges to $\delta_{\mathcal{T}_{\mathcal{K}}(\overline{X})}$. Then, by using the same argument as in the proof of Theorem \ref{thmcons}, we obtain that $\Phi_m + \delta_{\mathcal{F}_m}$ epi-converges in distribution to $\Phi + \delta_{\mathcal{T}_{\mathcal{K}}(\overline{X})}$. In addition, the optimal solution to (\ref{eqndellim}) is unique due to the strong convexity of $\Phi$ over the feasible set $\mathcal{K}$. Then, applying
Lemma \ref{lemepicon} on the epi-convergence in distribution leads to the desired result.
\end{proof}

Note that $\widehat{X}_m = \overline{X} + \rho_m \widehat{\Delta}_m$. From Lemmas \ref{lemrkrhs}, \ref{lemlocrk} and \ref{lemdellim}, we can see that $\overline{U}_2^\mathbb{T}\widehat{\Delta}\overline{V}_2=0$
is a necessary condition for the rank consistency of $\widehat{X}_m$. Then, we look into an explicit characterization of this condition.

\begin{lemma}\label{lemdelcond}
Let $\widehat{\Delta}$ be the optimal solution to the problem (\ref{eqndellim}).
Then $\overline{U}_2^\mathbb{T} \widehat{\Delta} \overline{V}_2 = 0$ if and only if the linear system (\ref{eqndelcond})
has a solution $\widehat{\Gamma} \in \mathbb{V}^{(n_1-r)\times (n_2-r)}$ with $\|\widehat{\Gamma}\| \leq 1$.
Moreover, in this case, $\widehat{\Delta} = \mathcal{Q}_\beta^\dag\big(\overline{U}_2\widehat{\Gamma}\,\overline{V}_2^\mathbb{T}-\overline{U}_1\overline{V}_1^\mathbb{T} + F(\overline{X})\big)$.
\end{lemma}
\begin{proof}
Assume that $\overline{U}_2^\mathbb{T}\widehat{\Delta}\overline{V}_2=0$. Since $\widehat{\Delta}$ is the optimal solution to (\ref{eqndellim}), from the optimality condition, the subdifferential of  $\|X\|_*$ at $0$,
and \cite[Theorem 23.7]{Roc70}, we obtain that there exist
some $\widehat{\Gamma} \in \mathbb{V}^{(n_1-r)\times (n_2-r)}$ with $\|\widehat{\Gamma}\|\le 1$ and $(\widehat{\eta}^0, \widehat{\eta}^1, \widehat{\eta}^2) \in \mathbb{R}^{|\alpha|} \times \mathbb{R}^{|\beta^+|} \times \mathbb{R}^{|\beta^-|}$ such that
\begin{equation}\label{eqndellimkkt}
\left\{
\begin{aligned}
& \mathcal{Q}_\beta(\widehat{\Delta}) + \overline{U}_1\overline{V}_1^\mathbb{T} - F(\overline{X})
+ \mathcal{R}_\alpha^*(\widehat{\eta}^0)
+ \mathcal{R}_{\beta^+}^*(\widehat{\eta}^1)
+ \mathcal{R}_{\beta^-}^*(\widehat{\eta}^2)
- \overline{U}_2 \widehat{\Gamma} \,\overline{V}_2^\mathbb{T}=0,\\
& \mathcal{R}_\alpha(\widehat{\Delta}) = 0,\\
& \mathcal{R}_{\beta^+}(\widehat{\Delta}) \leq 0, \quad \widehat{\eta}^1 \geq 0, \quad \langle \mathcal{R}_{\beta^+}(\widehat{\Delta}), \widehat{\eta}^1 \rangle = 0, \\
& \mathcal{R}_{\beta^-}(\widehat{\Delta}) \geq 0, \quad \widehat{\eta}^2 \leq 0, \quad \langle \mathcal{R}_{\beta^-}(\widehat{\Delta}), \widehat{\eta}^2 \rangle = 0.
\end{aligned}
\right.
\end{equation}
Note that $\mathcal{R}_{\beta^+}(\widehat{\Delta}) \leq 0$ and $\mathcal{R}_{\beta^-}(\widehat{\Delta}) \geq 0$ implies that $\mathcal{Q}_\beta^\dag\mathcal{Q}_\beta(\widehat{\Delta}) = \mathcal{P}_\beta(\widehat{\Delta})$. Moreover, $\mathcal{Q}_\beta^\dag \mathcal{R}_\alpha^*(\widehat{\eta}^0) = \mathcal{Q}_\beta^\dag\mathcal{R}_{\beta^+}^*(\widehat{\eta}^1) = \mathcal{Q}_\beta^\dag\mathcal{R}_{\beta^-}^*(\widehat{\eta}^2) = 0$. Then, we apply the operator $\mathcal{Q}_\beta^\dag$ to the first equation of (\ref{eqndellimkkt}) and then obtain
\begin{equation}\label{eqndellimkkt2}
\mathcal{P}_\beta(\widehat{\Delta}) + \mathcal{Q}_\beta^{\dag}(\overline{U}_1\overline{V}_1^\mathbb{T} - F(\overline{X}))
- \overline{U}_2 \widehat{\Gamma} \,\overline{V}_2^\mathbb{T})=0.
\end{equation}
Further note that $\mathcal{R}_\alpha(\widehat{\Delta}) = 0$ implies $\mathcal{P}_\alpha(\widehat{\Delta}) = 0$. This leads $\overline{U}_2^\mathbb{T}\mathcal{P}_\beta(\widehat{\Delta})\overline{V}_2=0$ since $\overline{U}_2^\mathbb{T}\widehat{\Delta}\overline{V}_2=0$. Then, together with(\ref{eqndellimkkt2}), we obtain that  $\widehat{\Gamma}$ is a solution to  (\ref{eqndelcond}).

Conversely, if the linear system (\ref{eqndelcond}) has a solution $\widehat{\Gamma}$ with $\|\widehat{\Gamma}\|\leq 1$, then it is easy to check that the KKT conditions (\ref{eqndellimkkt}) are satisfied with $\widehat{\Delta} = \mathcal{Q}_\beta^\dag(\widehat{Z})$
 and $ \widehat{\eta}^0 = \mathcal{R}_\alpha(\widehat{Z})$, $\widehat{\eta}^1 = (\mathcal{R}_{\beta^+}(\widehat{Z}))_+$, $\widehat{\eta}^2 = (\mathcal{R}_{\beta^-}(\widehat{Z}))_-$, where $\widehat{Z} = \overline{U}_2\widehat{\Gamma}\,\overline{V}_2^\mathbb{T}-\overline{U}_1\overline{V}_1^\mathbb{T} + F(\overline{X})$. Then, $\overline{U}_2^\mathbb{T} \widehat{\Delta} \overline{V}_2 = 0$ directly follows
from (\ref{eqndelcond}).
\end{proof}

  With Lemma \ref{lemdelcond}, the necessary part of Theorem \ref{thmnessuf} is immediate due to the necessity of the condition $\overline{U}_2^\mathbb{T}\widehat{\Delta}\overline{V}_2=0$
for rank consistency. Now we proceed with the sufficient part.

Define $\beta_m^+$, $\beta_m^-$, $\beta_m^\circ$ similar to (\ref{defbetadivision}) with $\overline{X}$ replaced by $\widehat{X}_m$. From Theorem \ref{thmcons}, we have $\widehat{X}_m \stackrel{p}{\rightarrow} \overline{X}$ as $m \rightarrow \infty$. The convergence implies that $\beta_m^+ \subseteq \beta^+$ and $\beta_m^- \subseteq \beta^-$ for sufficiently large $m$. In this circumstance,
the estimator $\widehat{X}_m$ is the optimal solution to (\ref{eqnrcs}) with $\mathcal{C} = \mathbb{V}^{n_1\times n_2}$
if and only if there exists a subgradient $\widehat{G}_m$ of the nuclear norm at $\widehat{X}_m$
and $(\widehat{\eta}_m^0, \widehat{\eta}_m^1, \widehat{\eta}_m^2) \in \mathbb{R}^{|\alpha|} \times  \mathbb{R}^{|\beta_m^+|} \times \mathbb{R}^{|\beta_m^-|}$ such that $(\widehat{X}_m, \widehat{\eta}_m^0, \widehat{\eta}_m^1, \widehat{\eta}_m^2)$ satisfies the KKT conditions:
\begin{equation}\label{eqnsufkkt}
\left\{
\begin{aligned}
&\frac{1}{m}\mathcal{R}_\Omega^* \big(\mathcal{R}_\Omega(\widehat{X}_m)\!-\!y\big)\!+\!\rho_m \big(\widehat{G}_m\!-\!F(\widetilde{X}_m)\big)
\!+\! \mathcal{R}_\alpha^*(\widehat{\eta}_m^0)
\!+\! \mathcal{R}_{\beta_m^+}^*(\widehat{\eta}_m^1)
\!+\! \mathcal{R}_{\beta_m^-}^*(\widehat{\eta}_m^2) = 0, \\
& \mathcal{R}_\alpha(\widehat{X}_m) = \mathcal{R}_\alpha(\overline{X}), \\
& \mathcal{R}_{\beta_m^\circ}(\widehat{X}_m) < b, \ \mathcal{R}_{\beta_m^+}(\overline{X}_m) = b, \ \mathcal{R}_{\beta_m^-}(\overline{X}_m) = -b, \ \eta_m^1 \geq 0,  \ \eta_m^2 \leq 0.
\end{aligned}
\right.
\end{equation}
 Let $ (\widehat{U}_m, \widehat{V}_m)\in \mathbb{O}^{n_1,n_2}(\widehat{X}_m)$ with $\widehat{U}_{m,1} \in \mathbb{O}^{n_1\times r}$, $\widehat{U}_{m,2} \in \mathbb{O}^{n_1\times (n_1-r)}$, $\widehat{V}_{m,1} \in \mathbb{O}^{n_2\times r}$ and $\widehat{V}_{m,2} \in \mathbb{O}^{n_2\times (n_2-r)}$. From Theorem \ref{thmcons} and Lemma \ref{lemrkrhs}, we know that $\text{rank}(\widehat{X}_m) \geq r$ with probability tending to one. When $\text{rank}(\widehat{X}_m) \geq r$ holds,
from the characterization of the subdifferential of the nuclear norm \cite{Wat92, Wat93}, we have that
 $\widehat{G}_m = \widehat{U}_{m,1} \widehat{V}_{m,1}^\mathbb{T} + \widehat{U}_{m,2} \widehat{\Gamma}_m\widehat{V}_{m,2}^\mathbb{T}$
 for some $\widehat{\Gamma}_m \in \mathbb{V}^{(n_1-r)\times (n_2-r)}$ satisfying $\|\widehat{\Gamma}_m\| \leq 1$. Now we want to show $\|\widehat{\Gamma}_m\| <1$ so that $\text{rank}(\widehat{X}_m)=r$. Since $\widehat{X}_m \stackrel{p}{\rightarrow} \overline{X}$,
 by \cite[Proposition 8]{DinST10} we have $\widehat{U}_{m,1}\widehat{V}_{m,1}^\mathbb{T} \stackrel{p}{\rightarrow} \overline{U}_1 \overline{V}_1^\mathbb{T}$. As $\widehat{\Gamma}$ is the unique optimal solution to (\ref{eqndelcond}), applying Lemma \ref{lemoper}  with the equation (\ref{eqnobs}) leads to
$$\frac{1}{m\rho_m}\mathcal{R}_\Omega^* \big(\mathcal{R}_\Omega(\widehat{X}_m)-y\big)+\widehat{U}_{m,1} \widehat{V}_{m,1}^\mathbb{T} - F(\widetilde{X}_m) \stackrel{p}{\rightarrow}  \mathcal{Q}_\beta(\widehat{\Delta})+\overline{U}_1 \overline{V}_1^\mathbb{T} - F(\overline{X}),$$
Then, by further applying the operator $\mathcal{Q}_\beta^\dag$ to the above equation, together with (\ref{eqndellimkkt2}) in Lemma \ref{lemdelcond} and  (\ref{eqnsufkkt}), we obtain   that
\begin{equation}\label{eqnsuflsconv}
\overline{U}_2^\mathbb{T}\mathcal{Q}_\beta^\dag(\widehat{U}_{m,2}\widehat{\Gamma}_m\widehat{V}_{m,2}^\mathbb{T})\overline{V}_2 \stackrel{p}{\rightarrow}
\overline{U}_2^\mathbb{T}\mathcal{Q}_\beta^\dag(\overline{U}_2\widehat{\Gamma}\overline{V}_2^\mathbb{T})\overline{V}_2.
\end{equation}
Since $\widehat{X}_m \stackrel{p}{\rightarrow} \overline{X}$, according to \cite[Proposition 7]{DinST10},
there exist two sequences of matrices $Q_{m,U} \in \mathbb{O}^{n_1-r}$ and  $Q_{m,V} \in \mathbb{O}^{n_2-r}$ such that
\begin{equation}\label{eqnsufuvprop}
\widehat{U}_{m,2}Q_{m,U} \stackrel{p}{\rightarrow} \overline{U}_2 \quad \text{and} \quad \widehat{V}_{m,2} Q_{m,V} \stackrel{p}{\rightarrow} \overline{V}_2.
\end{equation}
Moreover, the uniqueness of the solution to the linear system (\ref{eqndelcond}) is equivalent to the non-singularity of its linear operator. By combining (\ref{eqnsuflsconv}) and (\ref{eqnsufuvprop}), we obtain that
$Q_{m,U}^\mathbb{T} \widehat{\Gamma}_m Q_{m,V} \stackrel{p}{\rightarrow} \widehat{\Gamma}.$
Hence, we obtain that $\|\widehat{\Gamma}_m\|<1$ and thus $\text{rank}(\widehat{X}_m) = r$ with probability tending to one since $\|\widehat{\Gamma}\|<1$. Thus, we complete the proof of Theorem \ref{thmnessuf}.

\subsection{Proof of Theorem \ref{thmnessufpos}}

The proof of Theorem \ref{thmnessufpos} is similar to the proof of Theorem \ref{thmnessuf}. Define $\widehat{\Delta}_m: = \rho_m^{-1}(\widehat{X}_m-\overline{X})$.

\begin{lemma}\label{lemdellimpos}
If $\rho_m \rightarrow 0$ and $\sqrt{m}\rho_m\rightarrow \infty$, then $\widehat{\Delta}_m \stackrel{p}{\rightarrow} \widehat{\Delta}$, where $\widehat{\Delta}$ is the unique optimal solution to the following convex optimization problem
\begin{equation}\label{eqndellimpos}
\begin{aligned}
\min_{\Delta \in \mathbb{S}^n} & \ \ {\displaystyle \frac{1}{2}} \langle \mathcal{Q}_\beta(\Delta), \Delta \rangle + \langle I_n - F(\overline{X}), \Delta \rangle\\
{\rm s.t.} &  \ \ \mathcal{R}_\alpha(\Delta) = 0, \ \ \mathcal{R}_{\beta^+}(\Delta) \leq 0, \ \ \mathcal{R}_{\beta^-}(\Delta) \geq 0, \ \ \overline{P}_2^\mathbb{T} \Delta \overline{P}_2 \in \mathbb{S}_+^{n-r}.
\end{aligned}
\end{equation}
\end{lemma}
\begin{proof}
It is easy to verify that $\widehat{\Delta}_m$ is the optimal solution to
\begin{equation}\label{eqndellimapprpos}
\begin{aligned}
\min_{\Delta\in\mathbb{S}^n} & \ \ {\displaystyle \frac{1}{2m}}\|\mathcal{R}_\Omega(\Delta)\|_2^2 - \frac{\nu}{m \rho_m}\langle \mathcal{R}_\Omega^*(\xi), \Delta\rangle + \langle I_n - F(\widetilde{X}_m), \Delta\rangle\\
{\rm s.t.} & \ \ \Delta \in \mathcal{F}_m := \rho_m^{-1}(\mathcal{K} \cap \mathbb{S}_+^n-\overline{X}),
\end{aligned}
\end{equation}
where $\mathcal{K} := \big\{X \in \mathbb{S}^n \mid \mathcal{R}_\alpha(X) = \mathcal{R}_\alpha(\overline{X}), \ \|\mathcal{R}_\beta(X)\|_\infty \leq b\big\}$. Then, $\mathcal{F}_m$ converges in the sense of Painlev{\'e}-Kuratowski to
the tangent cone $\mathcal{T}_{\mathcal{K}\cap \mathbb{S}_+^n}(\overline{X})$ (see \cite{RocW98, BonS00}).
Note that the Slater condition in Assumption \ref{asmpslater} implies that $\mathcal{K}$ and $\mathbb{S}_+^n$ cannot be separated.
Then, from \cite[Theorem 6.42]{RocW98}, we have
$\mathcal{T}_{\mathcal{K}\cap \mathbb{S}_+^n}(\overline{X}) = \mathcal{T}_\mathcal{K}(\overline{X}) \cap \mathcal{T}_{\mathbb{S}_+^n}(\overline{X})$
with $\mathcal{T}_{\mathcal{K}}(\overline{X})$ taking the form of (\ref{eqntangentK}) and
$\mathcal{T}_{\mathbb{S}^n_+}(\overline{X})  = \big\{\Delta \in\mathbb{S}^n \mid \overline{P}_2^\mathbb{T} \Delta \overline{P}_2 \in \mathbb{S}_+^{n-r} \big\}$ according to Arnold \cite{Arn71}. Then, the proof can be completed by using the same argument as in the proof of Lemma \ref{lemdellim}.
\end{proof}

For the case $\mathcal{C} = \mathbb{S}_+^n$, Lemmas \ref{lemrkrhs}, \ref{lemlocrk} and \ref{lemdellimpos} imply that $\overline{P}_2^\mathbb{T}\widehat{\Delta}\overline{P}_2=0$ is a necessary condition
for the rank consistency of $\widehat{X}_m$. Then we look into an explicit characterization of this condition.

\begin{lemma}\label{lemdelcondpos}
Let $\widehat{\Delta}$ be the optimal solution to the problem (\ref{eqndellimpos}). Then $\overline{P}_2^\mathbb{T} \widehat{\Delta} \overline{P}_2 = 0$
if and only if the linear system (\ref{eqndelcondpos})
has a solution $\widehat{\Lambda} \in \mathbb{S}_+^{n-r}$.
Moreover, in this case,
$ \widehat{\Delta} = \mathcal{Q}_\beta^\dag\big(\overline{P}_2\widehat{\Lambda}\,\overline{P}_2^\mathbb{T}-I_n + F(\overline{X})\big).$
\end{lemma}
\begin{proof}
Note that the Slater condition also holds for the problem (\ref{eqndellimpos}). (One may check the point $X^0-\overline{X}$.)
Hence, $\widehat{\Delta}$ is the optimal solution to (\ref{eqndellimpos}) if and only if there exists
$(\widehat{\zeta}^0, \widehat{\zeta}^1, \widehat{\zeta}^2, \widehat{\Lambda}) \in \mathbb{R}^{|\alpha|} \times \mathbb{R}^{|\beta^+|} \times \mathbb{R}^{|\beta^-|} \times  \mathbb{S}^{n-r}$ such that
\begin{equation}\label{eqndellimkktpos}
\left\{
\begin{aligned}
& \mathcal{Q}_\beta(\widehat{\Delta}) + I_n - F(\overline{X})
+ \mathcal{R}_\alpha^*(\widehat{\zeta}^0)
+ \mathcal{R}_{\beta^+}^*(\widehat{\zeta}^1)
+ \mathcal{R}_{\beta^-}^*(\widehat{\zeta}^2)
- \overline{P}_2 \widehat{\Lambda} \overline{P}_2^\mathbb{T} =0,\\
& \mathcal{R}_\alpha(\widehat{\Delta}) = 0, \\
& \mathcal{R}_{\beta^+}(\widehat{\Delta}) \leq 0, \quad \widehat{\zeta}^1 \geq 0, \quad \langle \mathcal{R}_{\beta^+}(\widehat{\Delta}), \widehat{\zeta}^1 \rangle = 0, \\
& \mathcal{R}_{\beta^-}(\widehat{\Delta}) \geq 0, \quad \widehat{\zeta}^2 \leq 0, \quad \langle \mathcal{R}_{\beta^-}(\widehat{\Delta}), \widehat{\zeta}^2 \rangle = 0, \\
& \overline{P}_2^\mathbb{T}\widehat{\Delta} \overline{P}_2 \in \mathbb{S}_+^{n-r},\ \widehat{\Lambda} \in \mathbb{S}_+^{n-r},\ \langle  \overline{P}_2^\mathbb{T} \widehat{\Delta} \overline{P}_2, \widehat{\Lambda}\rangle =0.
\end{aligned}\right.
\end{equation}
Then, applying the operator $\mathcal{Q}_\beta^\dag$ to the first equation of (\ref{eqndellimkktpos}) yields the desired expression of $\widehat{\Delta}$ if $\overline{P}_2^\mathbb{T}\widehat{\Delta}\overline{P}_2=0$. It immediately follows that $\widehat{\Lambda}$ is a solution
 to (\ref{eqndelcondpos}).

Conversely, if the linear system (\ref{eqndelcondpos}) has a solution
$\widehat{\Lambda} \in \mathbb{S}_+^{n-r}$,  it is easy to check that  (\ref{eqndellimkktpos}) is satisfied
with $\widehat{\Delta} = \mathcal{Q}_\beta^\dag(\widehat{Z})$ and $\widehat{\zeta}^0 = \mathcal{R}_\alpha(\widehat{Z})$, $\widehat{\zeta}^1 = (\mathcal{R}_{\beta^+}(\widehat{Z}))_+$, $\widehat{\zeta}^2 = (\mathcal{R}_{\beta^-}(\widehat{Z}))_-$, where $\widehat{Z} = \overline{P}_2\widehat{\Lambda}\,\overline{P}_2^\mathbb{T}-I_n + F(\overline{X})$.
Then, $\overline{P}_2^\mathbb{T} \widehat{\Delta} \overline{P}_2 = 0$  directly follows from (\ref{eqndelcondpos}).
\end{proof}

The necessary part of Theorem \ref{thmnessufpos} is immediate from Lemma \ref{lemdelcondpos} due to the necessity of the condition $\overline{P}_2^\mathbb{T}\widehat{\Delta}\overline{P}_2=0$
for rank consistency. Now we proceed with the sufficient part.

Define $\beta_m^+$, $\beta_m^-$, $\beta_m^\circ$ by (\ref{defbetadivision}) with $\overline{X}$ replaced by $\widehat{X}_m$. From Theorem \ref{thmcons}, we have $\widehat{X}_m \stackrel{p}{\rightarrow} \overline{X}$ as $m \rightarrow \infty$. The convergence implies that $\beta_m^+ \subseteq \beta^+$ and $\beta_m^- \subseteq \beta^-$ for sufficiently large $m$. In this circumstance,
the Slater condition implies that $\widehat{X}_m$ is the optimal solution to (\ref{eqnrcs})
if and only if there exists  multipliers $(\widehat{\zeta}_m^0,
\widehat{\zeta}_m^1,\widehat{\zeta}_m^2, \widehat{S}_m) \in \mathbb{R}^{|\alpha|} \times \mathbb{R}^{|\beta^+|} \times \mathbb{R}^{|\beta^-|} \times \mathbb{S}^{n}$
such that $(\widehat{X}_m, \widehat{\zeta}_m^0, \widehat{\zeta}_m^1, \widehat{\zeta}_m^2, \widehat{S}_m)$ satisfies the KKT conditions:
\begin{equation}\label{eqnsufkktpos}
\left\{
\begin{aligned}
&\frac{1}{m}\mathcal{R}_\Omega^* \big(\mathcal{R}_\Omega(\widehat{X}_m)\!-\!y\big)\!+\!\rho_m \big(I_n\!-\! F(\widetilde{X}_m)\big)
\!+\!\mathcal{R}_\alpha^*(\widehat{\zeta}_m^0)
\!+\!\mathcal{R}_{\beta^+}^*(\widehat{\zeta}_m^1)
\!+\!\mathcal{R}_{\beta^-}^*(\widehat{\zeta}_m^2)
\!-\! \widehat{S}_m =0, \\
&\mathcal{R}_\alpha(\widehat{X}_m) = \mathcal{R}_\alpha(\overline{X}), \\
& \mathcal{R}_{\beta_m^\circ}(\widehat{X}_m) < b, \ \mathcal{R}_{\beta_m^+}(\overline{X}_m) = b, \ \mathcal{R}_{\beta_m^-}(\overline{X}_m) = -b, \ \eta_m^1 \geq 0,  \ \eta_m^2 \leq 0,\\
&\widehat{X}_m \in \mathbb{S}_+^n,\ \widehat{S}_m \in \mathbb{S}_+^n,\ \langle \widehat{X}_m, \widehat{S}_m \rangle =0.
\end{aligned}\right.
\end{equation}
 The last equation in (\ref{eqnsufkktpos}) implies that $\widehat{X}_m$ and $\widehat{S}_m$ can have a
 simultaneous eigenvalue decomposition. Let $\widehat{P}_m \in \mathbb{O}^{n}(\widehat{X}_m)$ with $\widehat{P}_{m,1} \in \mathbb{O}^{n\times r}$ and $\widehat{P}_{m,2} \in \mathbb{O}^{n\times (n-r)}$. From Theorem \ref{thmcons} and Lemma \ref{lemrkrhs}, we know that $\text{rank}(\widehat{X}_m) \geq r$ with probability tending to one. When $\text{rank}(\widehat{X}_m) \geq r$ holds, we can write
 $\widehat{S}_m = \widehat{P}_{m,2} \widehat{\Lambda}_m\widehat{P}_{m,2}^\mathbb{T}$ for some diagonal matrix $\widehat{\Lambda}_m \in \mathbb{S}^{n-r}_+$. In addition, if $\widehat{\Lambda}_m \in \mathbb{S}^{n-r}_{++}$, then $\text{rank}(\widehat{X}_m)=r$. Since $\widehat{X}_m \stackrel{p}{\rightarrow} \overline{X}$, according to \cite[Proposition 1]{DinST10}, there exist a sequence of matrices $Q_m \in \mathbb{O}^{n-r}$ such that
 $\widehat{P}_{m,2}Q_m \stackrel{p}{\rightarrow} \overline{P}_2$.
 Then, using the similar arguments to the proof of Theorem \ref{thmnessuf}, we obtain that
 $Q_m^\mathbb{T} \widehat{\Lambda}_m Q_m \stackrel{p}{\rightarrow} \widehat{\Lambda}$.
 Since $\widehat{\Lambda}\in \mathbb{S}_{++}^n$, we have $\widehat{\Lambda}_m \in \mathbb{S}_{++}^n$
 with probability tending to one. Thus, we complete the proof of Theorem \ref{thmnessufpos}.

\subsection{Proof of Theorem \ref{thmsoluni}}

We first prove for the rectangular case $\mathcal{C} = \mathbb{V}^{n_1\times n_2}$ by contradiction. Assume that there exists some $\mathbb{V}^{(n_1-r)\times (n_2-r)} \ni \overline{\Gamma} \neq 0$ such that
$\mathcal{B}_2(\overline{\Gamma}) = \overline{U}_2^\mathbb{T}\mathcal{Q}_\beta^\dag(\overline{U}_2 \overline{\Gamma} \, \overline{V}_2^\mathbb{T})\overline{V}_2  = 0$. Then $\langle \overline{\Gamma}, \overline{U}_2^\mathbb{T}\mathcal{Q}_\beta^\dag(\overline{U}_2 \overline{\Gamma} \, \overline{V}_2^\mathbb{T})\overline{V}_2  \rangle
=\langle \overline{U}_2 \overline{\Gamma} \, \overline{V}_2^\mathbb{T}, \mathcal{Q}_\beta^\dag(\overline{U}_2 \overline{\Gamma} \, \overline{V}_2^\mathbb{T})\rangle = 0$. This immediately leads to $(\mathcal{Q}_\beta^\dag)^{1/2}(\overline{U}_2 \overline{\Gamma} \, \overline{V}_2^\mathbb{T})= 0$ since $\mathcal{Q}_\beta^\dag$ is a self-adjoint and positive semidefinite operator. It then follows that
$[\mathcal{R}_{\beta^\circ}; (\mathcal{R}_{\beta^+})_-; (\mathcal{R}_{\beta^-})_+](\overline{U}_2 \overline{\Gamma} \, \overline{V}_2^\mathbb{T}) = 0$, where $(\mathcal{R}_\pi)_\pm(\cdot):=(\mathcal{R}_\pi(\cdot))_\pm$ with $\pi = \beta^+$ or $\beta^-$. Then by using this equality, we have that for any $H \in \mathcal{T}(\overline{X})$,
\begin{align*}
 0 = & \ \langle \overline{\Gamma},  \overline{U}_2^{\mathbb{T}} H \overline{V}_2 \rangle = \langle \overline{U}_2 \overline{\Gamma} \, \overline{V}_2^\mathbb{T}, H\rangle = \langle \mathcal{R}_{\alpha\cup\beta}(\overline{U}_2 \overline{\Gamma} \, \overline{V}_2^\mathbb{T}), \mathcal{R}_{\alpha\cup\beta}(H) \rangle \\
= & \ \langle [\mathcal{R}_\alpha; (\mathcal{R}_{\beta^+})_+; (\mathcal{R}_{\beta^-})_-](\overline{U}_2 \overline{\Gamma} \, \overline{V}_2^\mathbb{T}), \mathcal{R}_{\alpha\cup\beta^+\cup\beta^-}(H)\rangle.
\end{align*}
By using the arbitrariness of $\mathcal{R}_{\alpha\cup\beta^+\cup\beta^-}(H)$ over $\mathbb{R}^{|\alpha \cup \beta^+ \cup \beta^-|}$ implied by the constraint nondegeneracy (\ref{eqncndc}), we further have $[\mathcal{R}_\alpha; (\mathcal{R}_{\beta^+})_+; (\mathcal{R}_{\beta^-})_-](\overline{U}_2 \overline{\Gamma} \, \overline{V}_2^\mathbb{T})=0$. Therefore, we obtain $\overline{U}_2 \overline{\Gamma} \, \overline{V}_2^\mathbb{T} = 0$ and thus $\overline{\Gamma} = 0$, which leads to a contradiction. Therefore, the linear operator $\mathcal{B}_2$ is positive definite. The proof for the positive semidefinite case is similar.

\subsection{Proof of Theorem \ref{thmrccordencons}}

We first prove for the constraint nondegeneracy.
\begin{lemma}\label{lemcorcn}
For the matrix completion problems of Classes I and II, the constraint nondegeneracy (\ref{eqncndcpos}) holds at $\overline{X}$.
\end{lemma}
\begin{proof}
For the real covariance matrix case, the proof is given in \cite[Lemma 3.3]{QiS06} and \cite[Proposition 2.1]{QiS11}. For the complex covariance matrix case, one can use the similar arguments to prove the result.

We next consider the density matrix case. Suppose that $\overline{X}$ satisfies the density constraint,
i.e., $\mathcal{R}_{\alpha}(\overline{X}) =\frac{1}{\sqrt{n}}\text{Tr}(\overline{X}) = \frac{1}{\sqrt{n}}$. Note that for any $t\in\mathbb{R}$,
 we have $t \overline{X} \in \text{lin}(\mathcal{T}_{\mathcal{H}_+^n}(\overline{X}))$. This,
along with $\text{Tr}(\overline{X})=1$, implies that
\[\frac{1}{\sqrt{n}}\text{Tr} \big(\text{lin}(\mathcal{T}_{\mathcal{H}_+^n}(\overline{X}))\big)   = \mathcal{R}_{\alpha} \big(\text{lin}(\mathcal{T}_{\mathcal{H}_+^n}(\overline{X}))\big) =\mathbb{R}.\]
This means that the constraint nondegeneracy (\ref{eqncndcpos}) holds.
\end{proof}

From Theorem \ref{thmsoluni} and Lemma \ref{lemcorcn}, for both  Classes I an II, the linear system (\ref{eqndelcondpos}) has a unique solution $\widehat{\Lambda}$. Moreover, for both Classes I and II, uniform sampling yields $\mathcal{Q}_\beta^\dag(Z) = \mathcal{P}_\beta(Z) /d_2$ for any $Z \in \mathbb{S}_+^n$. Thus, from (\ref{eqndelcondpos}), we have
\begin{equation}\label{eqnrccordenconseq}
\widehat{\Lambda} - \overline{P}_2^\mathbb{T} \mathcal{P}_\alpha (\overline{P}_2 \widehat{\Lambda} \overline{P}_2^\mathbb{T}) \overline{P}_2
= \overline{P}_2^\mathbb{T} \mathcal{P}_\beta (\overline{P}_2 \widehat{\Lambda} \overline{P}_2^\mathbb{T}) \overline{P}_2
= \overline{P}_2^\mathbb{T} \mathcal{P}_\beta (I_n - F(\overline{X})) \overline{P}_2.
\end{equation}

Then we first prove for Class I by contradiction. For any $Z\in\mathbb{S}_+^n$, $\mathcal{P}_\alpha(Z)$ is the diagonal matrix whose $i$-th diagonal entries is $X_{ii}$ for all $i \in \pi$ and the other entries are $0$.
 %The left-hand side of (\ref{eqnrccordenconseq}) can be written as
%$$\overline{P}_2^\mathbb{T}\mathcal{P}_\beta(\overline{P}_2\widehat{\Lambda}\overline{P}_2^\mathbb{T})\overline{P}_2 =\widehat{\Lambda} - \overline{P}_2^\mathbb{T}\mathcal{P}_\alpha(\overline{P}_2\widehat{\Lambda}\overline{P}_2^\mathbb{T})\overline{P}_2.$$
Assume that $\widehat{\Lambda} \notin \mathbb{S}^{n-r}_{++}$, i.e., $\lambda_{\rm min}(\widehat{\Lambda})\leq 0$, where $\lambda_{\rm min}(\cdot)$ denotes the smallest eigenvalue. Then, we have
$$\lambda_{\rm min}(\widehat{\Lambda}) = \lambda_{\rm min}(\overline{P}_2\widehat{\Lambda}\overline{P}_2^\mathbb{T}) \leq \lambda_{\rm min}\big(\mathcal{P}_\alpha(\overline{P}_2\widehat{\Lambda}\overline{P}_2^\mathbb{T})\big) \leq \lambda_{\rm min}\big(\overline{P}_2^\mathbb{T}\mathcal{P}_\alpha (\overline{P}_2\widehat{\Lambda}\overline{P}_2^\mathbb{T})\overline{P}_2\big),
$$
where the equality follows from the fact that $\widehat{\Lambda}$ and $\overline{P}_2 \widehat{\Lambda}\overline{P}_2^\mathbb{T}$ have the same nonzero eigenvalues, the first inequality follows from the fact that the vector of eigenvalues is majorized by the vector of diagonal entries, (e.g., see \cite[Theorem 9.B.1]{MarOA10}), and the second inequality follows from the Courant-Fischer minmax theorem, (e.g., see \cite[Theorem 20.A.1]{MarOA10}). As a result, the left-hand side of (\ref{eqnrccordenconseq}) is not positive definite. Notice that $\overline{P}_2^\mathbb{T} F(\overline{X})\overline{P}_2=0$. Thus, the right-hand side of (\ref{eqnrccordenconseq}) can be written as
$$\overline{P}_2^\mathbb{T}\mathcal{P}_\beta(I_n-F(\overline{X}))\overline{P}_2 = \overline{P}_2^\mathbb{T}\mathcal{P}_\beta(I_n)\overline{P}_2 %-\overline{P}_2^\mathbb{T}\big(F(\overline{X})-\mathcal{P}_\alpha(F(\overline{X}))\big)\overline{P}_2 =
+\overline{P}_2^\mathbb{T} \mathcal{P}_\alpha(F(\overline{X}))\overline{P}_2 = \overline{P}_2^\mathbb{T} \big(\mathcal{P}_\beta(I_n)+\mathcal{P}_\alpha(F(\overline{X}))\big)\overline{P}_2.$$
Since $\text{rank}(\overline{X})=r$, with the choice (\ref{eqnfcorden}) of $F$, we have that for any $i \in \pi$,
$$\overline{X}_{ii} = \sum_{j=1}^r\lambda_j(\overline{X}) |\overline{P}_{ij}|^2 >0 \quad \text{implies} \quad
 \big(F(\overline{X})\big)_{ii} = \sum_{j=1}^r f_i\big(\lambda_j(\overline{X})\big) |\overline{P}_{ij}|^2  >0.$$
 Moreover, $\mathcal{P}_\beta(I_n)$ is the diagonal matrix with the last $n-r$ diagonal entries being ones and the other entries being zeros. Thus, $\mathcal{P}_\beta(I_n)+\mathcal{P}_\alpha(F(\overline{X}))$ is a diagonal matrix with all positive diagonal entries. It follows that the right-hand side of (\ref{eqnrccordenconseq}) is positive definite. Thus, we obtain a contradiction. Therefore, we should have $\widehat{\Lambda} \in \mathbb{S}^{n-r}_{++}$. Then, we can obtain the rank consistency according to Theorem \ref{thmnessufpos}.

Next, we prove for Class II. It is easy to see $\mathcal{P}_\alpha(\cdot) = \frac{1}{n}\text{Tr}(\cdot)I_n$. By further using $\overline{P}_2^\mathbb{T} F(\overline{X})\overline{P}_2= 0$ and $\mathcal{P}_\beta(I_n)=0$, we can rewrite (\ref{eqnrccordenconseq})   as
$$  \widehat{\Lambda}-\frac{1}{n}\text{Tr}(\widehat{\Lambda})I_{n-r} = \frac{1}{n}\text{Tr}(F(\overline{X}))I_{n-r}.$$
By taking the trace on both sides, we obtain that
$\widehat{\Lambda} = \frac{1}{r}\text{Tr}(F(\overline{X}))I_{n-r}.$
Since $\overline{X}$ is a density matrix of rank $r$, with the choice (\ref{eqnfcorden}) of $F$, we have that
$$\text{Tr}(\overline{X}) = \sum_{i=1}^n \sum_{j=1}^r\lambda_j(\overline{X}) |\overline{P}_{ij}|^2 = 1 \quad \text{implies} \quad
\text{Tr}\big(F(\overline{X})\big) = \sum_{i=1}^n\sum_{j=1}^r f_i\big(\lambda_j(\overline{X})\big) |\overline{P}_{ij}|^2  >0.$$
It follows that $\widehat{\Lambda} \in \mathbb{S}_{++}^{n-r}$ and thus we obtain the rank consistency.

{\small
\bibliographystyle{plain}
\bibliography{rank_correction}
}

\end{document}